\documentclass[11pt,letterpaper]{amsart}

\usepackage{amssymb,amsfonts}
\usepackage[all,arc]{xy}
\usepackage{enumerate}
\usepackage{mathrsfs}
\usepackage{tikz-cd}
\usepackage{mathtools}
\usepackage{bm}
\usepackage{eucal}
\usepackage{quiver} 
\usepackage{scalerel}
\usepackage{stackengine,wasysym}
\usepackage{csquotes}
\usepackage[foot]{amsaddr}
\usepackage{amsmath,calligra,mathrsfs}
\usepackage{hyperref}
\usepackage[capitalise]{cleveref}
\usepackage[style=alphabetic,
minalphanames=1,
maxalphanames=7,
maxbibnames=99]{biblatex}
\newtheorem{thm}{Theorem}[section]
\newtheorem*{thm*}{Theorem}
\newtheorem{cor}[thm]{Corollary}
\newtheorem{prop}[thm]{Proposition}
\newtheorem{lem}[thm]{Lemma}
\newtheorem{conj}[thm]{Conjecture}
\newtheorem{quest}[thm]{Question}
\newtheorem{claim}[thm]{Claim}

\theoremstyle{definition}
\newtheorem{defn}[thm]{Definition}

\newtheorem{construction}[thm]{Construction}
\newtheorem{exmp}[thm]{Example}

\newtheorem{notn}[thm]{Notation}

\newtheorem{property}[thm]{Property}

\newtheorem{spec}[thm]{Speculation}

\theoremstyle{remark}

\newtheorem{warning}[thm]{Warning}

\newtheorem{disc}[thm]{Discussion}
\newtheorem{rem}[thm]{Remark}

\makeatletter
\let\c@equation\c@thm
\makeatother
\numberwithin{equation}{section}

\makeatletter
\renewcommand{\email}[2][]{%
	\ifx\emails\@empty\relax\else{\g@addto@macro\emails{,\space}}\fi%
	\@ifnotempty{#1}{\g@addto@macro\emails{\textrm{(#1)}\space}}%
	\g@addto@macro\emails{#2}%
}
\makeatother

\usepackage{geometry}
\newcommand{\beqn}{\begin{equation}}
\newcommand{\eeqn}{\end{equation}}
\newcommand{\bclaim}{\begin{claim}}
\newcommand{\eclaim}{\end{claim}}
\newcommand{\blem}{\begin{lem}}
\newcommand{\elem}{\end{lem}}
\newcommand{\bproof}{\begin{proof}}
\newcommand{\eproof}{\end{proof}}

\newcommand{\An}{\mathbb{A}^n}

\newcommand{\Aone}{\mathbb{A}^1}
\newcommand{\const}{\underline{e}}

\newcommand{\hook}{\hookrightarrow}
\newcommand{\F}{\mathcal{F}}
\newcommand{\bdef}{\begin{defn}}
\newcommand{\edefn}{\end{defn}}
\newcommand{\bprop}{\begin{prop}}
\newcommand{\eprop}{\end{prop}}
\newcommand{\bthm}{\begin{thm}}
\newcommand{\ethm}{\end{thm}}

\newcommand{\Four}{\mathscr{F}}
\newcommand{\Z}{\bm{Z}}
\newcommand{\X}{\bm{X}}

\renewcommand{\a}{\alpha}
\renewcommand{\b}{\beta}

\newcommand{\T}{\bm{T}}
\newcommand{\N}{\mathcal{N}}
\renewcommand{\O}{\mathcal{O}}
\renewcommand{\i}{\iota}

\renewcommand{\L}{\mathbb{L}}
\newcommand{\brem}{\begin{rem}}
\newcommand{\erem}{\end{rem}}
\newcommand{\Q}{\mathbb{Q}}
\newcommand{\K}{\mathscr{K}}
\newcommand{\G}{\mathbb{G}}
\newcommand{\longhook}{\lhook\joinrel\longrightarrow}
\renewcommand{\P}{\mathbb{P}}
\newcommand{\Bl}{\mathrm{Bl}}
\newcommand{\D}{\bm{D}}
\newcommand{\0}{\{0\}}
\newcommand{\f}{\bm{f}}
\newcommand{\act}{\mathrm{act}}
\renewcommand{\-}{\textendash}
\newcommand{\Shv}{\mathrm{Shv}}
\newcommand{\too}{\twoheadrightarrow}
\newcommand{\vircodim}{\mathrm{codim.vir}}
\newcommand{\cBl}{\mathrm{Bl}}
\newcommand{\wt}{\widetilde}
\newcommand{\pr}{\mathrm{pr}}
\newcommand{\U}{\mathcal{U}}
\newcommand{\V}{\mathcal{V}}
\newcommand{\W}{\mathcal{W}}

\newcommand{\C}{\mathbb{C}}
\newcommand{\R}{\mathbb{R}}
\newcommand{\crit}{\mathrm{crit}}
\newcommand{\Cee}{\mathscr{C}}

\newcommand{\Fun}{\mathrm{Fun}}

\newcommand{\op}{\mathrm{op}}
\newcommand{\pt}{\mathrm{pt}}
\newcommand{\supp}{\mathrm{supp}}

\newcommand{\Sch}{\mathrm{Sch}}
\newcommand{\oblv}{\mathrm{oblv}}
\newcommand{\PrShv}{\mathrm{PrShv}}
\newcommand{\h}{\mathrm{h}}
\newcommand{\B}{\mathcal{B}}

\newcommand{\dcrit}{\mathbf{crit}}
\newcommand{\codim}{\mathrm{codim}}
\newcommand{\LKE}{\mathrm{LKE}}
\newcommand{\BM}{\mathrm{BM}}
\newcommand{\bbT}{\mathbb{T}}
\newcommand{\wit}{\widetilde}
\newcommand{\con}{\mathrm{con}}
\newcommand{\alg}{\mathrm{alg}}
\newcommand{\an}{\mathrm{an}}

\DeclareMathOperator{\QCoh}{QCoh}

\DeclareMathOperator{\colim}{colim}
\DeclareMathOperator{\Perv}{Perv}
\DeclareMathOperator{\CAlg}{CAlg}
\DeclareMathOperator{\Spc}{Spc}
\DeclareMathOperator{\Mod}{Mod}
\DeclareMathOperator{\cdga}{cdga}
\DeclareMathOperator{\Map}{Map}
\DeclareMathOperator{\Sym}{Sym}
\DeclareMathOperator{\fib}{fib}

\DeclareMathOperator{\Oblv}{Oblv}
\DeclareMathOperator{\Spec}{Spec}
\DeclareMathOperator{\Hom}{Hom}
\DeclareMathOperator{\Sp}{Sp}
\DeclareMathOperator{\id}{id}
\DeclareMathOperator{\DSch}{DSch}
\DeclareMathOperator{\VDiv}{VDiv}
\DeclareMathOperator{\Real}{Re}
\DeclareMathOperator{\sHom}{\mathscr{H}\text{\kern -3pt {\calligra\large om}}\,}
\DeclareMathOperator{\ob}{ob}

\title{Microlocal homology}

\author{Kendric Schefers}

\address{Department of Mathematics, The University of California, Berkeley}

\email{kdschefers@gmail.com}

\begin{document}

\begin{abstract}
Let $Z$ be an l.c.i. scheme over $\C$.
In this paper, we introduce a Kashiwara--Schapira-style functor of derived microlocalization,
which we use to define a perverse sheaf $\mu_{Z}$ on the $-1$-shifted cotangent bundle, $T^*[-1]Z$. 
The sheaf $\mu_{Z}$ is designed to be a refinement of the \emph{microlocal homology} of $Z$: 
a family of invariants introduced by Nadler that
interpolates between the singular cohomology and Borel--Moore homology of $Z$.
	
Our main result is an equivalence between $\mu_{Z}$ and
the DT sheaf $\varphi_{T^*[-1]Z}$ on $T^*[-1]Z$. 
This provides an alternative construction for the DT sheaf 
in the case of a shifted cotangent bundle.
The main step of our argument, which may be of
independent interest, is a local computation---similar to
one obtained recently by Kinjo using different methods---providing 
a description of the classical microlocalization functor in terms
of vanishing cycles.

\end{abstract}

\maketitle

\tableofcontents


\section{Introduction}

\subsection{Overview}
On a smooth, oriented manifold $X$ of dimension $m$, 
Borel--Moore homology and singular cohomology are 
equivalent under Poinar\'e duality:
\[H_{m-\bullet}^{\BM}(X) \simeq H^{\bullet}(X).\]
On a singular topological space, Poincar\'e duality
fails, and Borel--Moore homology diverges from
singular cohomology. That is, rather than an equivalence,
we obtain a map $H^{\bullet}(X) \to H_{m-\bullet}^{\BM}(X)$.
An interesting mathematical question is whether
one can characterize the difference between
between Borel--Moore homology and 
singular cohomology in terms of the singularities
of the space in question.

Local complete intersections in complex algebraic
geometry are very tame examples of singular spaces, 
and a sensible investigation of this question 
might begin by considering the zero locus of a
flat, proper map $f: X \to V$, from a smooth complex scheme $X$
of dimension $m$
to a finite dimensional vector space $V$
of dimensional $n$.
Such an investigaton led David Nadler,
in his course lecture notes (\cite{Y13}), to propose
a definition for a family of objects,
parametrized by conical subsets of $V^{\vee}$,
which interpolates between the singular cohomology
and Borel--Moore homology of $Z := f^{-1}(0)$.
whose precise definition we recall below.

\bdef[cf. \cite{Y13}]
\label{def: microlocal homology}
Let $\Lambda \subset V^{\vee}$
be a closed conical subset. The 
\emph{$\Lambda$-microlocal homology} of 
$Z$ is defined as,
\[H_{\bullet}^{\Lambda}(f) := \Gamma_{\Lambda}
	(V^{\vee}; \mu_0(f_*\underline{\C}_{X}[m])),\]
where $\mu_0$ denotes the Kashiwara--Schapira functor
of microlocalization (see \cite{KS90}),
$\Gamma_{\Lambda}(V^{\vee};-)$ denotes
the functor of derived global sections 
with support in $\Lambda$, and $\underline{\C}_X[m]$
is the perverse constant sheaf on $X$.
\edefn

The important features of the $\Lambda$-microlocal homology
are that
\begin{itemize}
\item $H_{\bullet}^{\{0\}}(f)$ recovers singular cohomology
of $Z$ with a shift by $m-n$;\footnotemark
\footnotetext{Compare the shift here with that in
our formulation of Poincar\'e duality above in the
case when $Z$ is smooth.}
\item $H_{\bullet}^{V^{\vee}}(f)$ recovers Borel--Moore
homology of $Z$;  
\item and that, for conical $\Lambda' \subset \Lambda$, there
exists a map $H_{\bullet}^{\Lambda'}(f) \to H_{\bullet}^{\Lambda}(f)$.
\end{itemize}
These features of the microlocal homology 
allow us to extract a naive \emph{singular support} theory for
classes in Borel--Moore homology: the singular support
of an element $[\sigma] \in H_{\bullet}^{\BM}(Z)$
is the smallest closed conical subset $\Lambda \subset V^{\vee}$
such that $[\sigma]$ is in the image of the
map $H_{\bullet}^{\Lambda}(f) \to H_{\bullet}^{\BM}(Z)$.
This characterizes the obstruction to representing
a class in Borel-Moore homology by a singular cochain in
terms of the (co)directions along the base of the
map $f: X \to V$, which has something to do with
the singularities of $Z$.

The $\Lambda$-microlocal homology 
of $Z$, however, depends a priori
on the auxiliary choice of a map $f$
presenting $Z$ as a zero locus. The choice of
such a map, while not totally unrelated to
the singularities of $Z$, is not intrinsic to $Z$.
A characterization of the difference between Borel--Moore
homology and cohomology in terms of the singularities
of $Z$ morally should not depend on the baggage of such choices.
Our goal in this paper, simply put, is to show 
that the microlocal homology is more intrinsic to $Z$ 
than it initially appears.

We accomplish this goal by defining a microlocal
homology \emph{sheaf} on the $-1$ shifted cotangent
bundle of $Z$, $T^*[-1]Z$,\footnotemark
\footnotetext{To clarify, $T^*[-1]Z$ for us denotes
the underlying classical scheme of the derived
scheme $\Spec \Sym_{\O_Z}(\bbT_Z[1])$,
where $\bbT_Z$ is the tangent complex of $Z$.
Alternatively, $T^*[-1]Z := \Spec \Sym_{\O_Z}(H^1(\bbT_Z))$.} 
also known as \textit{the scheme of singularities
of $Z$}. This is a scheme over $Z$ which, in a sense, records
information about the singularities of $Z$: 
$T^*[-1]Z$ has fiber zero over each point
in the regular locus of $Z$. Importantly,
$T^*[-1]Z$ is \emph{intrinsic} to $Z$ (as a scheme).
We show that sections of the microlocal homology sheaf recover
Nadler's microlocal homology, thereby accomplishing our goal.
The following theorem may be seen as a distillation
of our main results.

\begin{thm*}
Given any proper local complete intersection $Z$ over $\C$ that
admits an embedding $Z \hook X$ for $X$ smooth, there
exists a perverse sheaf $\mu_Z$ on $T^*[-1]Z$ such that
$\Gamma_{\Lambda}(T^*[-1]Z; \mu_Z) \simeq H_{\bullet}^{\Lambda}(f)$,
for any $f$ presenting $Z$ as a zero fiber.
\end{thm*}

In order to construct $\mu_Z$ we introduce a
Kashiwara--Schapira-style functor of microlocalization 
in the setting of derived algebraic geometry. 
The proof of the properties of $\mu_Z$ 
uses as its main input a highly non-trivial calculation
which reveals a surprising
and beautiful connection with cohomological
Donaldson--Thomas theory, which may be of
independent interest to the reader.

Although the discussion so far
has required $f$ be flat and $Z$ be
a local complete intersection, in the body of
the paper we allow arbitrary maps $f: X \to V$ 
by using the formalism of derived algebraic geometry
and considering $Z$ as a quasi-smooth derived
scheme.

\subsection{Main results}
Let $\Z$ be a quasi-smooth derived
scheme (``derived l.c.i.")
over $\C$, and let $\Z \hook X$ be any
closed immersion into a smooth scheme $X$.
Then there is a well-defined notion
of derived normal bundle to $\Z$ inside
of $X$, which we denote $\N_{\Z/X}$.
Recent developments in derived
algebraic geometry by Khan and Rydh (\cite{KR19})
allow one to define a derived deformation to the
normal bundle: a derived scheme over $\Aone$
whose special fiber is $\N_{\Z/X}$ with a map
to $X$ and all the
expected properties of an object called
``deformation to the normal bundle."

Using the input of this derived deformation
to the normal bundle, we define the functors
of \emph{quasi-smooth} specialization and
microlocalization,
\begin{align*}
\Sp_{\Z/X}: \Shv(X) \to \Shv(\N_{\Z/X}) \\
\mu_{\Z/X}: \Shv(X) \to \Shv(\N^{\vee}_{\Z/X})
\end{align*} 
by imitating the definitions of their classical
counterparts found in \cite{KS90} and \cite{V83b}.
Our functors are extensions of the classical
ones to the derived setting in the following sense.
In the case when $\Z$ is classical, and therefore 
a local complete intersection, then
our functors agree with the classical ones.

\brem
\label{rem: Adeel's work}
A forthcoming paper of Adeel Khan (\cite{AAK})
independently develops the 
notion of derived microlocalization
for \'etale sheaves along arbitrary morphisms
(i.e.\ not necessarily a quasi-smooth 
closed immersion).
\erem

\subsubsection{First main theorem}
Our first main theorem is the following.

\bthm[= \cref{theorem: local statement}]
\label{intro theorem: local statement}
Suppose that $X$ is a smooth scheme, and $f:
X \to V$ is a regular function whose target
is an $n$-dimensional complex vector space. 
Let $\Z(f) \hook X$ denote the 
derived zero locus of $f$, and let
$\widetilde{f}: X \times V^{\vee} \to \Aone$
denote the map $(x, \lambda) \mapsto 
\langle \lambda, f(x) \rangle$. Then there
exists a natural isomorphism of functors,
	\[\varphi_{\widetilde{f}}(\pr_1^*-) \simeq \mu_{\Z(f)/X}(-),\]
where $\varphi_{\widetilde{f}}$ denotes the functor
of vanishing cycles with respect to $\widetilde{f}$.
 \ethm
The proof of this theorem starts with the
observation that the claimed isomorphism
is trivial in the case when $X$ is bundle on
$\Z(f)$ and the input sheaf is monodromic.
This observation reduces the problem to 
a highly non-trivial computation showing
that vanishing cycles ``factors through"
the specialization functor along $\Z(f)$.

\brem
\label{rem: Kinjo's work}
A result closely related to
\cref{intro theorem: local statement}
was recently shown by Kinjo (\cite[Theorem 4.1]{K21b}) using
different methods. In contrast to our result, however, Kinjo's is purely classical.
Our introduction of the derived microlocalization functor $\mu_{\Z(f)/X}$
clarifies the nature of Kinjo's and our equivalences, and
allows for a clean conceptual statement.
 
When $\Z(f)$ is classical, hence a
local complete intersection variety, 
our quasi-smooth microlocalization
agrees with the classical microlocalization
functor, and Kinjo's statement is
essentially equivalent to our own.

Interestingly, 
Kinjo's proof strategy and our own, while
very different, both similarly work 
only in the Betti setting
(both using the complex analytic
topology in an essential way).
\erem

\subsubsection{Background on Donaldson--Thomas theory}
The derived critical locus of the function $\wit{f}$
from \cref{intro theorem: local statement} is isomorphic
to the derived scheme $\T^*[-1]\Z$. In fact, more is true: 
both are canonically oriented $-1$-shifted symplectic, and
they are isomorphic as such.

Oriented $-1$-shifted geometry is the setting for the modern
study of Donaldson--Thomas theory, as the moduli
spaces of sheaves on Calabi-Yau threefolds have
the structure of oriented $-1$-shifted symplectic
stacks of virtual dimension $0$. The virtual counts
of point in these moduli spaces are known as
``Donaldson--Thomas invariants." 
In the simplest example, the moduli
space admits a description as a derived
critical locus of a regular function $f: U \to \Aone$, 
the Donaldson--Thomas invariants are given by 
the Euler characteristic of 
the perverse vanishing cycles sheaf $\varphi_f(\underline{\C}_U[\dim U])$,
so $\varphi_f(\underline{\C}_U[\dim U])$ in a sense
``categorifies" Donaldson--Thomas invariants.

Derived critical loci are the
prototypical examples of $-1$-shifted
symplectic spaces. In fact, any $-1$-shifted symplectic
scheme $(\X, \omega_{\X})$ admits an open cover in which
it is given locally by derived critical loci.  
In \cite{BBDJSS15}, the authors show that, 
if $(\X, \omega_{\X})$
is oriented, the collection vanishing
cycles given by any such cover of critical charts
glues to a perverse sheaf on $\X$ which is
\emph{independent} of the choice of cover. 
More precisely, they show that
to any oriented $-1$-shifted 
symplectic scheme, $(\X, \omega_{\X}, o_{\X})$, 
there is a canonical\footnotemark
\footnotetext{up to unique isomorphism
determined by choice of square root
of the determinant line 
bundle, $\det(\L_{\X})|_{X}$.} 
perverse sheaf on the underlying 
classical scheme, $X$, denoted
$\varphi_{\X}$
such that if
$\dcrit(f) \simeq \bm{U} \subset \X$
is an open critical chart of $\X$,
then ${\varphi_{\X}}|_{\bm{U}}$
is isomorphic to the perverse vanishing cycles sheaf with
respect to $f$---up to a twist
by a $\mathbb{Z}/2\mathbb{Z}$-bundle.
This sheaf plays a central role in the modern
study of Donaldson--Thomas theory, as it
categorifies Donaldson--Thomas invariants
and their generalizations.

\subsubsection{Second main theorem}
Our second main theorem globalizes the
isomorphism proven in our first by relating it
to the perverse sheaf of twisted vanishing
cycles on $\T^*[-1]\Z$ for a quasi-smooth
derived scheme $\Z$, with its canonical oriented
$-1$-shifted symplectic structure.

\bthm[= {\cref{theorem: global statement}}]
\label{intro theorem: global statement}
Given any closed immersion $\Z \hook X$ 
of a quasi-smooth derived scheme $\Z$ into 
a smooth scheme $X$, there exists an isomorphism,
\[ \varphi_{\T^*[-1]\Z} \simeq \mu_{\Z/X}[\vircodim(\Z,X)](\const_{X}[\dim X]),\]
of perverse sheaves on $\T^*[-1]\Z \subset \bm{\N}^{\vee}_{\Z/X}$, where $e$ is any
well-behaved commutative ring\footnotemark.
\ethm

\footnotetext{e.g.\ of finite homological
dimension. The reader who wishes to
will lose nothing by taking $e := \mathbb{C}$
throughout this entire paper.}

\Cref{intro theorem: global statement}
gives a new construction of the perverse
sheaf $\varphi_{\T^*[-1]\Z}$ as the microlocalization
along any embedding of $\Z$ into a smooth
scheme. 
The original construction detailed in \cite{BBDJSS15} 
is difficult and highly technical, even in
the case of a shifted cotangent bundle,
as we have here. In contrast, the new method
of construction suggested by our results is
much simpler and more functorial.
From its definition, the sheaf $\mu_{\Z/X}$, 
however, is not obviously independent 
of the choice of embedding $\Z \hook X$.
On the other hand, once it is defined,
the intrinsicality of $\varphi_{\T^*[-1]\Z}$
to $\Z$ is automatically seen.
In this way, these independent 
constructions of the same perverse 
sheaf are mutually complementary.

\subsubsection{Singular support for Borel--Moore chains}
Using the yoga of the six functor formalism for sheaves,
it is not difficult to deduce from \cref{intro theorem: global
statement} that global sections of $\varphi_{\T^*[-1]\Z}$
with prescribed support recovers Nadler's 
$\Lambda$-microlocal homology whenever
$\Z$ is the derived zero fiber of a proper map $f: X \to V$.
More precisely, we have the following theorem.

\bthm[= {\cref{theorem: equivalence of microlocal homologies}}]
\label{intro theorem: equivalence of microlocal homologies}
Let $\Lambda \subset V^{\vee}$ be
a closed conic subset, and suppose $f: X \to V$
is a proper map from a smooth scheme $X$
to a finite dimensional vector space $V$. Let
$\wit{\Lambda}$ denote the intersection of
$Z(f) \times \Lambda \subset Z(f) \times V^{\vee}
\simeq \N_{\Z(f)/X}$ with $T^*[-1]\Z(f) \subset
\N_{\Z(f)/X}$.
Then there exists an equivalence-up-to-shifts,
\[\Gamma_{\wit{\Lambda}}\left(\varphi_{\T^*[-1]\Z(f)}\right) 
	\simeq H_{\bullet}^{\Lambda}(f),\]
where ``\,$\Gamma_{\wit{\Lambda}}(-)$" denotes derived
global sections with support contained in $\wit{\Lambda}$.
\ethm

This theorem allows us finally to extract from the
microlocal homology a plausible candidate for
a singular support theory of classes in 
Borel--Moore homology. Namely, for $f$
proper, the singular support
of an element $[\sigma] \in H_{\bullet}^{\BM}(Z(f))$
is the minimal subset $\Lambda \subset T^*[-1]\Z(f)$ such
that $[\sigma]$ is the image of the natural map,
\begin{align*}
\Gamma_{\Lambda}\left(\varphi_{\T^*[-1]\Z(f)}\right) &\to \Gamma\left(\varphi_{\T^*[-1]\Z(f)}\right) \\
											&\simeq H_{\bullet}^{\BM}(Z(f)).
\end{align*}


\subsection{Additional significance of results}
Nadler's microlocal homology is a general construction,
of interest to anyone working on microlocalization
and the topology of singular spaces, and we believe that the motivation
for our results on these invariants was sufficiently explained earlier
in the introduction. We wish to indicate the
significance of our results beyond
the points we have already mentioned.

\subsubsection{Vanishing cycles and microlocalization}
In the classical 
theory of $\mathscr{D}$-modules
and microlocal sheaf theory, the functor
of vanishing cycles and the functors of
microlocalization and Fourier transform
are intimately related in well-known ways. 
For example, if $f: X \to \C$ 
is a holomorphic function on a complex manifold $X$
and the zero locus of $f$ is non-singular, 
$\varphi_f$ is entirely determined by
the microlocalization along $f^{-1}(0)$ 
(\cite[Proposition 8.6.3]{KS90}).
On the other hand, the stalk of the
Fourier transform at a generic covector
is vanishing cycles with respect
to the functional determined by that covector
(\cite{V83a}). 

Yet, these known results stop
short of providing a \emph{full} description
of the microlocalization and Fourier--Sato 
transform functors as a vanishing 
cycles---and vice versa. \Cref{intro theorem: local statement}
can seen as furnishing such a description,
writing the \emph{entire} microlocalization functor in
terms of vanishing cycles. 

\brem
\label{rem: local statement}
We believe that this equivalence
between microlocalization and vanishing cycles
is a fundamental one that should be true in
any sheaf theory context where there are
formal notions of nearby cycles, vanishing cycles,
and Fourier transform---such as $\mathscr{D}$-modules,
$\ell$-adic sheaves, mixed Hodge modules, etc.---though
we in this paper (and Kinjo in \cite{K21b}) only show it 
for $\C$-constructible sheaves (cf. \cref{rem: Kinjo's work}).
\erem

\subsubsection{Donaldson--Thomas theory}
While the class of $-1$-shifted symplectic schemes
which we consider in \cref{intro theorem: global statement}
may seem trivially narrow, shifted cotangent bundles of
this type are already an important and interesting class
of spaces in the study of cohomological and
motivic Donaldson--Thomas theory.

For example, let $\widetilde{S}$ be a 
local surface (i.e.\ the total space of 
the canonical line bundle
on a smooth quasi-projective surface $S$). 
Denote by $\bm{M}_{\widetilde{S}}$ and $\bm{M}_S$ 
the derived moduli stacks
of coherent sheaves with 
proper supports on $\widetilde{S}$ and $S$, 
respectively. Then $\bm{M}_{\widetilde{S}}$
carries a canonical $-1$-shifted symplectic structure,
and $\bm{M}_{S}$ is quasi-smooth. 
In fact, \cite[Theorem 5.1]{K21a} shows
that $\bm{M}_{\widetilde{S}} \simeq \T^*[-1]\bm{M}_S$
as $-1$-shifted symplectic stacks with the latter's
canonical shifted symplectic structure.
When, for example, the underlying coarse moduli
scheme $\pi_0(M_S)$ is quasi-projective,
$\pi_0(\bm{M}_S)$ admits an
embedding into a smooth scheme,
and \cref{intro theorem: global statement}
gives a microlocal interpretation of
the categorification of the classical 
Donaldson--Thomas invariants
of the Calabi-Yau 3-fold, $\widetilde{S}$\footnotemark.
\footnotetext{See \cite[Remark 6.14]{BBDJSS15}.}

Similarly, we might take the category
of representations of a quiver with
superpotential. This furnishes a CY3
category whose moduli stack of
objects can locally be written as
the derived critical locus of a holomorphic
function $f: U \to \mathbb{C}$, for smooth
$U$. Simple examples of such stacks can
sometimes globally be expressed
as the critical locus of such a function,
in which case, \cref{intro theorem: global statement}
(in fact, even \cref{intro theorem: local statement})
gives an alternative construction for
its canonical twisted vanishing cycles sheaf.


\subsection{Structure of the paper}

The content of this paper is organized as follows.

In \cref{sec: notation and conventions} we establish
the notation and conventions used in this work.

In \cref{sec: quasi-smooth specialization and microlocalization},
we define the functors of quasi-smooth specialization
quasi-smooth microlocalization. We establish some basic 
properties expected of the quasi-smooth specialization,
such as base change and a lemma 
on restriction to the zero section.
The classical notions of Verdier specialization and the real
analytic specialization of \cite{KS90} are also discussed
and compared.

In \cref{sec: the local equivalence}, we state and prove
our first main result, \cref{theorem:
local statement}.

In \cref{sec: the global equivalence}, we
state and prove our second main result,
\cref{theorem: global statement}, after recalling
the basic concepts involved in its statement.

In \cref{sec: microlocal chains}, we
introduce the $\Lambda$-microlocal homology
discussed earlier in the introduction and prove its
compatibility with the $\Lambda$-microlocal homology sheaf
we define using \cref{theorem: global statement}. 
For the $\Lambda$-microlocal homology sheaves
we extract a well-behaved singular support theory for classes
in Borel--Moore homology. In this section, we also speculate
about possible descriptions for the $\Lambda$-microlocal homology
spaces in terms of currents with prescribed wavefront sets.

In \cref{sec: sheaf theory}, we give an account of
sheaf theory in the modern setting. Care is taken
to justify the use of classical results in the theory
of sheaves on locally compact spaces. The reader
is reminded of various elements of the sheaf theory
needed for the body of the paper. The theory of
constructible sheaves in the $\infty$-categorical setting
is well-known to experts, but the author could not find
a systematic account of it elsewhere in the literature.

In \cref{sec: deformation to the normal bundle}, 
we recall from \cite{KR19} the construction of 
deformation to the normal bundle along a 
quasi-smooth closed immersion. We prove 
several properties about the construction
in preparation for next section on the functors
of quasi-smooth specialization and microlocalization.


\subsection{Acknowledgments}
I am indebted to my
advisors, David Ben-Zvi and Sam Raskin:
this paper would not exist without them.
Their guidance and influence is felt throughout this work.
Both of them have been extremely generous
with their time and energy over the years,
and I am truly grateful for them.

I give special thanks to Rok Gregoric for
teaching me most of my knowledge of derived
geometry and $\infty$-categories, 
and to Tom Gannon for listening to
my ideas and answering my questions over the years. 

Additionally, I would like to
thank Germ\'an Stefanich for many insightful 
conversations; Ben Davison for taking the
time to answer my questions about his work;
Brian Hepler for explaining aspects
of singular topology to me; 
and Adeel Khan for his thoughtful comments
and suggestions on an earlier draft of this work.


\section{Notation and conventions}
\label{sec: notation and conventions}

\subsubsection{Base ring and sheaf coefficients}
For the remainder of this paper,
$k := \C$ and $e$ is as stated in \cref{intro
theorem: global statement}. This 
choice of notation was made to 
suggest that results analogous to those we 
prove should be true for
sheaf theories other than constructible
sheaves with $e$ coefficients---such as $\ell$\-adic
sheaves, $\mathscr{D}$-modules, mixed
Hodge modules, etc.---on schemes
over a base other than $(k:=)\C$
(see \cref{rem: local statement}
above).

\subsubsection{$\infty$-categories}
We freely use the language of
$\infty$-categories throughout
the paper, and take Joyal's
quasi-categories (i.e. weak
Kan complexes) as our particular model.
Unless otherwise specified, our notation
and conventions for the theory of
$\infty$-categories and related
subjects are taken to be those
appearing in the canonical references,
\cite{HTT} and \cite{HA}.

\subsubsection{Derived algebraic geometry}
The theory of 
derived algebraic geometry
is indispensable for the formulation
of the results of this paper.
There are three choices for a
model of affine derived schemes
over a commutative ring $R$:
(the opposite category of)
simplicial commutative rings, 
(connective) cdgas, and (connective) $\mathbb{E}_{\infty}$-algebras.
If $R$ contains the field of
rational numbers, then all three
choices give rise to the same theory
of derived algebraic geometry\footnotemark.
We will always work over the base
ring $\C$ in this paper, so our
choice of local model is immaterial. 
Nonetheless, in order to access
the results of \cite{SAG} more
directly, we take (the opposite
category of) $\mathbb{E
}_{\infty}$-algebras (i.e. $\mathbb{E}
_{\infty}$-algebra objects in spectra) 
as our model for affine derived
schemes. In light of the above
commentary, we use the terms
``spectral scheme" and
``derived scheme" interchangeably.
\footnotetext{See the introduction of
\cite[ch. 25]{SAG} for a detailed
comparison of the three choices.}

\subsubsection{Miscellanea}
The following is a list of miscellaneous
notation and conventions used in this
paper, but which have no natural home
in the main body of the text.

\begin{itemize}
\item Boldface letters such as
$\X, \bm{Y}, \Z$ always
denote derived schemes or stacks;
likewise, boldface letters such 
as $\bm{f}, \bm{g}, \bm{h}$ 
always denote morphisms
of derived schemes or stacks.
The corresponding plain letters,
$X, Y, Z$ and $f, g, h$,
denote the underlying classical
schemes and morphisms of
underlying classical schemes, respectively,
of their boldface counterparts,
when the meaning is clear from context.

\item We denote by $pt_X: X \to \ast$ the terminal 
map to a point from the topological space $X$. We
omit the subscript when $X$ is clear from context.

\item For $\ast = \{>0, \geq 0, <0, \leq 0\}$,
$\Aone_{\ast}$ denotes the
subset of $\{z \in \C| \Real(z) \, \ast \} \in \C$,
consisting of complex numbers
whose real parts are $\ast$.

\item For $\ast = \{>0, \geq 0, <0, \leq 0\}$, 
$\ell_{\ast}$ denotes the subset $\R_{\ast} 
\subset \C = \Aone$.

\item $\R^+$ denotes the multiplicative 
group of strictly positive real numbers.

\item If $\Cee$ is a small, ordinary
category, it can be viewed as
an $\infty$-category by taking the
nerve, $N(\Cee)$, which
gives a simplicial set. By abuse 
of notation, we often use 
$\Cee$ to refer to the 
simplicial set, $N(\Cee)$.
\end{itemize}


\section{Quasi-smooth specialization and microlocalization}
\label{sec: quasi-smooth specialization and microlocalization}

Recall the functor of specialization
introduced by Verdier in \cite[\S8]{V83b}.
Let $X$ be a separated scheme of
finite type over $\Spec \mathbb{C}$ and let $Z$ be a 
closed subscheme. Then Verdier 
specialization along $Z$ of sheaves 
on $X$ is a functor
\[ \Sp_{Z/X}: \Shv(X) \to \Shv_{\G_m}(C_{Z/X}) \]
where $C_{Z/X}$ is the normal cone of
$Z$ inside $X$.

It is defined using the deformation to 
the normal cone of $Z$ inside $X$, 
a scheme defined as $D_{Z/X} :=
\cBl_{Z \times \{0\}}(X \times \Aone) \setminus
\cBl_Z(X)$. The scheme $D_{Z/X}$ fits into
the following commutative diagram,
\begin{equation}
\begin{tikzcd}
\label{diagram: deformation 1}
	{C_{Z/X}} & {D_{Z/X}} & {X \times \G_m} & X \\
	{\{0\}} & \Aone & {\G_m}
	\arrow[hook, from=2-1, to=2-2]
	\arrow[from=2-3, to=2-2]
	\arrow["j"', from=1-3, to=1-2]
	\arrow["t"', from=1-2, to=2-2]
	\arrow[from=1-1, to=2-1]
	\arrow[from=1-3, to=2-3]
	\arrow["{i_{C_{Z/X}}}", hook, from=1-1, to=1-2]
	\arrow["{\pr_1}", from=1-3, to=1-4]
\end{tikzcd}
\end{equation}

\bdef
The functor of Verdier specialization is
defined by the formula:
\[\Sp_{Z/X}(\F) := \psi_t j_! \pr_1^*(F), \]
where $\psi_t$ is the nearby cycles functor
with respect to the map $t$ (see \cref{definition: nearby cycles}).
\edefn

\brem
A full account of the properties of
the specialization functor is given
in its original source \cite[\S8]{V83b}. 
\erem 


In this section
we define a functor of quasi-smooth 
specialization generalizing
the above construction. 
It is defined in the same 
manner as Verdier's
specialization functor, only
replacing the deformation
to the normal cone with
an appropriate notion of deformation
to the normal \emph{bundle} along
a quasi-smooth closed immersion $\Z \to \X$. 
The appropriate
such notion for us is that
appearing in \cite{KR19}.

The construction of
deformation to the normal
bundle introduced 
in \cite{KR19} is recalled
in \cref{sec: deformation to the normal
bundle} after the definition of
quasi-smooth closed immersion is given. 
In that appendix, we also prove some
general results about the construction which
are needed in order to set up 
the theory of specialization.

\brem
The specialization functor we define coincides 
with Verdier specialization
in the case when $\Z$ and
$\X$ are classical, meaning the
map $\Z \hook \X$
is a regular embedding.
\erem

At the end of this section,
we define quasi-smooth
microlocalization as the
Fourier--Sato transform of
quasi-smooth specialization
and briefly mention one of its
properties.


\subsection{Quasi-smooth specialization}
As mentioned above, \cite{KR19} gives a deformation to
the derived normal bundle for the quasi-smooth
closed immersion $\Z \to \X$ furnishing a diagram as follows:
\begin{equation}
\label{eqn: deformation diagram main 2}
\begin{tikzcd}
	{\bm{N}_{\Z/\X}} & {\bm{D}_{\Z/\X}} & {\X \times \G_m} & {\X} \\
	{\{0\}} & \Aone & {\G_m}
	\arrow[hook, from=2-1, to=2-2]
	\arrow[from=2-3, to=2-2]
	\arrow[from=1-1, to=2-1]
	\arrow[from=1-3, to=2-3]
	\arrow["{\bm{t}}", from=1-2, to=2-2]
	\arrow["{\bm{pr}_1}", from=1-3, to=1-4]
	\arrow["{\bm{j}_{\neq 0}}"', from=1-3, to=1-2]
	\arrow["{\bm{s}}", hook, from=1-1, to=1-2].
\end{tikzcd}
\end{equation}
See \cref{ssec: definition of deformation to the normal bundle} for a review
of this construction.

The following definition
is now very natural.

\bdef
Suppose that $i: \Z \hook \X$ is a 
quasi-smooth closed immersion of derived
schemes. Then the specialization functor
along $i$ is given by the formula:
\[ \Sp_{\Z/\X}(\F) := \psi_{t} j_{\neq 0 !} \pr_1^*(\F), \]
where $\F \in \Shv(X)$ and $\psi_t$ is again the nearby
cycles functor with respect to $t$\footnotemark in the above diagram
\labelcref{eqn: deformation diagram main 2}.
\footnotetext{where, by our notational convention, $t$ denotes the
the map induced by $\bm{t}$ on the underlying classical varieties.}
\edefn

Like its classical counterpart,
$\Sp_{\Z/\X}(\F)$ is monodromic,
lying in $\Shv_{\G_m}(N_{\Z/\X})$,
by virtue of the following well-known lemma.

\blem
\label{lemma: monodromic}
Suppose that $X$ is a scheme 
with a $\G_m$\-action, and that 
$f: X \to \Aone$ is $\G_m$\-equivariant.
Then if $\F \in \Shv(X \setminus Z(f))$ is equivariant,
$\psi_f(\F)$ is monodromic with respect
to the induced $\G_m$\-action on $Z(f)$.
\elem

The structure map 
$t: D_{\Z/\X} \to \Aone$ 
is obviously $\G_m$\-equivariant, 
and the choice of an equivalence
$\pr_1^*\F \simeq \F \boxtimes \const_{\G_m}$
gives $\pr_1^*\F$ the structure of
an equivariant sheaf.
Therefore, $\Sp_{\Z/\X}(\F)$ 
is monodromic.


\subsection{Specialization and algebraic structure}
Deformation to the normal
bundle, and by extension
quasi-smooth specialization, is
sensitive to the non-reduced,
algebro-geometric 
information of the closed immersion
$\Z \hook \X$, even for regular
embeddings of classical schemes.

While the category $\Shv(\bm{Z})
:= \Shv(Z)$ 
by definition only depends on
the associated analytic space
$Z^{\an}$, which itself only depends
on the underlying reduced, classical 
subscheme $Z^{red}$, 
performing deformation 
to the normal bundle along
two closed immersions
possessing the same underlying
reduced subscheme can produce
schemes with genuinely
different reduced subschemes
and analytifications thereof.
 
\subsubsection{}
The dependence of
quasi-smooth specialization on even the
classical nilpotent fuzz of $\Z$ is readily 
exhibited by looking at
the case when $\Z \hook \X$ is given
by the quotient of an
discrete algebra by a regular ideal
	\[A \too A/I.\] 
In this case, deformation to 
the normal bundle is given by 
the Rees algebra of $I$, $R(I)$,
which we may explicitly compute
in the following simple example. 

\begin{exmp}
Consider the closed immersions
given by the following maps of $k$-algebras:
\begin{align*}
	k[x] &\too k[x]/(x) \simeq k \\
	k[x] &\too k[x]/(x^2) =: k[\varepsilon].
\end{align*}
The first is the usual embedding, 
	\[\0 \hook \Aone,\] 
and the second is an embedding of the
dual numbers at zero, 
	\[\Spec (k[\varepsilon]) \hook \Aone.\] 
Both are regular closed embeddings, since
neither $x$ nor $x^2$ is a zero divisor,
and they share the same reduced
subscheme; yet specialization along 
each of these two immersions yields
different sheaves, as shown below.
\end{exmp}

\subsubsection{} 
We compute $\Sp_{\0/\Aone}(\const)$. 
The Rees algebra $R((x))$ is given by the 
$k$-algebra $k[x,t,\hbar]/(t\hbar - x)$, 
and the canonical map from deformation 
to the normal cone to $\Aone$ is given by 
the obvious map $k[t] \to k[x,t,\hbar]/(t\hbar - x)$. 
The analytification, $\Spec(R((x)))^{\an}$, 
is the affine variety $\{t\hbar - x = 0\} \subset \mathbb{C}^3$. 
The fibers of the map $p$ in 
the definition of the specialization 
above are subvarieties given by 
intersections of this variety with 
the hyperplanes $\{t = t_0\}$ for 
various values of $t_0$. The fiber 
at $t_0 \neq 0$ is a linear subspace 
defined by the equation $t_0\hbar = x$; 
the fiber at $t_0 = 0$ is the $\hbar$-axis. 
One can solve a system of equations to 
the find the region given by the intersection 
of a small ball in centered at a point $(0,0, a)$ 
in the zero fiber of $p$ with a nearby fiber at 
$|t_0| \ll 1$ to find that the local Milnor fiber 
at any point is contractible. The stalk of 
nearby cycles is given by the cohomology 
of the local Milnor fiber, the stalk of 
$\Sp_{\0/\Aone}(\const)$ at any point is 
the one dimension vector space $e$.

\subsubsection{} 
For simplicity of notation, denote 
$D := \Spec(k[\varepsilon])$. 
We compute $\Sp_{D/\Aone}(\const)$. 
The Rees algebra $R((x^2))$ is given by 
the $k$-algebra $k[x,t,\hbar]/(t\hbar - x^2)$, 
so the associated complex analytic space 
is a singular quadric surface with singularity 
at $(0,0,0)$. The fibers of the map $p$ are 
parabolas (over $\mathbb{C}$) 
defined at $t=t_0 \neq 0$ by the equation 
$\hbar = x^2/t_{0}$; the fiber at $t_0 = 0$ 
(of the analytification) 
is again the $\hbar$-axis. One can, as 
above, solve a system of equations to the 
find the region given by the intersection 
of a small ball centered at a non-zero point 
$(0,0, a \neq 0)$ in the zero fiber of $p$ with 
a nearby fiber at $|t_0|\ll 1$ to find that 
the Milnor fiber at any non-zero point is 
now a disjoint union of two contractible 
subspaces. The intuition one should have 
(coming from the real picture) is that as 
we approach the zero fiber, the nearby 
parabolas become narrower and narrower--while 
still isomorphic to the complex line--and 
in the limit merge into a complex 
line over $t_0 =0$. Any sufficiently small 
ball at a non-zero point on the $\hbar$-axis 
will then intersect with a nearby fiber in 
disjoint segments of each of the two ``legs" 
of the parabola. Consequently, the stalk of 
$\Sp_{D/\Aone}(\const)$ at 
any such point is the two(!) dimensional 
vector space $e \oplus e$.



\subsection{Properties of quasi-smooth specialization}

The functor $\Sp_{\Z/\X}$,
which we call \emph{quasi-smooth
specialization},
possesses most of the
properties enjoyed by its
classical counterpart.
We list and prove a 
few such properties.

\subsubsection{Proper and smooth base change}

The following two lemmas verify
the base change properties expected
of quasi-smooth specialization.


\blem
\label{lemma: proper base change for specialization}
Suppose that $\bm{f}: \X' \to \X$ is a
map of derived schemes, $\Z \hook \X$
and $\Z' \hook \X'$ are quasi-smooth
closed immersions such that
$\Z' \simeq \X' \times_{\X} \Z$.
Then there exists a natural morphism
\[\epsilon_{\bm{f}}: \Sp_{\Z/\X}(f_*\F) \to
	N_{\bm{f}*}(\Sp_{\Z'/\X'}(\F)) \]
for $\F \in \Shv(X')$.
Moreover, $\epsilon_{\bm{f}}$ is
an isomorphism when $\bm{f}$
is proper.
\elem

\bproof
Using diagram \labelcref{diagram:
classical map of deformations},
we obtain the chain of base change 
morphisms:
\begin{align*}
\Sp_{\Z/\X}(f_*\F) &:= \psi_{t}(j_{\neq 0 !}\pr_1^*f_*\F) \\
					&\to \psi_{t}(j_{\neq 0 !}(f \times \id_{\G_m})_* \pr'^*_1\F) \\
					&\to \psi_{t}(D_{\bm{f}*}j'_{\neq 0 !}\pr'^*_1\F) \\
				        &\to N_{\bm{f}*}\psi_{t'}(j'_{\neq 0 !}\pr'^*_1\F) \\
					&=: N_{\bm{f}*}\Sp_{\Z'/\X'}(\F).
\end{align*}
The composition of all these
morphisms is, by definition,
$\epsilon_{\bm{f}}$.
Since all the morphisms in the
above chain are natural in $\F$,
$\epsilon_{\bm{f}}$ is as well, proving
the first part of our claim.

Clearly, all the arrows in the
above chain of morphisms
are isomorphisms when
$f$ and $N_{\bm{f}}$ are
proper. Assume that $\bm{f}$
is proper as a morphism of
derived schemes \footnotemark. 
\footnotetext{This notion
is defined and discussed at length 
in \cite[ch. 5]{SAG}. We use its
results freely and without
precise reference, but they
may be all be found in there.}
We proceed to show that
$f$ and $N_{\bm{f}}$ are
proper as maps of classical schemes.

By \cite[Theorem 4.1.5(ii)]{KR19}
or \cref{lemma: blowup map}
we have the following Cartesian
diagram:
\[\begin{tikzcd}
	{\Bl_{\Z' \times \0}(\X' \times \Aone)} & {\Bl_{\Z \times \0}(\X \times \Aone)} \\
	{\X' \times \Aone} & {\X \times \Aone}
	\arrow["{\bm{\pi}'}"', from=1-1, to=2-1]
	\arrow["{\bm{f} \times \id_{\Aone}}", from=2-1, to=2-2]
	\arrow["{\bm{\pi}}"', from=1-2, to=2-2]
	\arrow["{\Bl_{\bm{f} \times \id_{\Aone}}}", from=1-1, to=1-2].
\end{tikzcd}\]
Since $\bm{f} \times \id_{\Aone}$
is proper, $\Bl_{\bm{f} \times \id_{\Aone}}$
is proper, as properness
is preserved by base change.
By the proof of
\cref{cor: deformation map}
above, we also have the Cartesian
diagram,
\[\begin{tikzcd}
	{\D_{\Z'/\X'}} & {\Bl_{\Z' \times \0}(\X' \times \Aone)} \\
	{\D_{\Z/\X}} & {\Bl_{\Z \times \0}(\X \times \Aone)}
	\arrow["{\D_{\bm{f}}}", from=1-1, to=2-1]
	\arrow["{\Bl_{\bm{f} \times \id_{\Aone}}}", from=1-2, to=2-2]
	\arrow[from=2-1, to=2-2]
	\arrow[from=1-1, to=1-2],
\end{tikzcd}\]
so $\D_{\bm{f}}$ is proper as well.

A proper map of derived schemes
over $k$ ($= \mathbb{C}$)
is easily seen to induce a proper
map on underlying classical
schemes, so both
$D_{\bm{f}}$ and $f$
are proper. 
The classical notion 
of properness is also preserved
under base change, so we deduce from
diagram \labelcref{diagram: classical
map of deformations} that
$N_{\bm{f}}$ is indeed proper.
\eproof


\blem
\label{lemma: smooth base change for specialization}
Suppose that $\bm{f}: \X' \to \X$ is a
map of derived schemes, $\Z \hook \X$
and $\Z' \hook \X'$ are quasi-smooth
closed immersions such that
$\Z' \simeq \X' \times_{\X} \Z$.
Then there exists a natural morphism
\[\eta_{\bm{f}}: N_{\bm{f}}^*(\Sp_{\Z/\X}(\F)) \to
	\Sp_{\Z'/\X'}(f^*\F) \]
for $\F \in \Shv(X)$.
Moreover, $\eta_{\bm{f}}$ is
an isomorphism on the locus where
$N_{\bm{f}}$ is smooth.
\elem

\bproof
Using diagram \labelcref{diagram:
classical map of deformations},
we obtain the chain of morphisms:
\begin{align*}
N_{\bm{f}}^*\Sp_{\Z/\X}(\F) &:= N_{\bm{f}}^*\psi_t(j_{\neq 0 !}\pr_1^*\F) \\
					&\to \psi_{t'}(D_{\bm{f}}^*j_{\neq 0 !}\pr_1^*\F) \\
					&\simeq \psi_{t'}(j'_{\neq 0 !}(f \times \id_{\G_m})^*\pr_1^*\F) \\
					&\simeq \psi_{t'}(j'_{\neq 0 !}\pr'^*_1f^*\F) \\
					&=: \Sp_{\Z'/\X'}(f^*\F).
\end{align*}
The composition of all these
morphisms is, by definition,
$\eta_{\bm{f}}$.
Since all the morphisms in the
above chain are natural in $\F$,
$\eta_{\bm{f}}$ is as well, proving
the first part of our claim.

Clearly $\eta_{\bm{f}}$ is an isomorphism
at points in $N_{\Z'/\X'}$ where
the base change morphism
	\[N_{\bm{f}}^*\psi_t(-) \to \psi_{t'}(D_{\bm{f}}^*(-)) \]
is an isomorphism. This is true at points
where the maps $N_{\bm{f}}$ and
$D_{\bm{f}}$ are smooth, so to prove the second
part of our claim it suffices to show
that $D_{\bm{f}}$ is smooth at those
points where $N_{\bm{f}}$ is
smooth. Recall that a map of
schemes is smooth at a
point if it is finitely
presented at that point
and the stalk of the relative cotangent
complex at that point is free.

We've shown earlier that $D_{\bm{f}}$
is a morphism of finite presentation.
Now choose $x \in N_{\Z'/\X'}$ where
$N_{\bm{f}}$ is smooth, and consider
the stalk at $x$ of the fiber sequence:
	\[(s'^*\L_{D_{\bm{f}}})_x \to (\L_{D_{\bm{f}} \circ s'})_x \to (\L_{s'})_x. \]
By assumption, $(\L_{N_{\bm{f}}})_x$
is free. The fiber sequence,
	\[N_{\bm{f}}^*\L_s \to \L_{s \circ N_{\bm{f}}} \to \L_{N_{\bm{f}}},\]
exhibits an equivalence $\L_{N_{\bm{f}}} \simeq
\L_{s \circ N_{\bm{f}}} = \L_{D_{\bm{f}} \circ s'}$
since $s$ is unramified as a closed
immersion. Thus $(\L_{D_{\bm{f}} \circ s'})_x$
is free. As a closed immersion,
$s'$ is also unramified, 
so $\L_{s'} \simeq 0$. 
Therefore, $(s'^*\L_{D_{\bm{f}}})_x \simeq
(\L_{D_{\bm{f}}})_x$ is free as well.
\eproof

The following lemma verifies
that quasi-smooth specialization
commutes with Verdier duality.

\blem
\label{lemma: dual of specialization}
For any $\F \in \Shv(X)$,
\[\mathbb{D}(\Sp_{\Z/\X}\F) \simeq \Sp_{\Z/\X}(\mathbb{D}\F). \]
\elem

\bproof
The standard properties of Verdier
duality and nearby cycles yield the
following chain of isomorphisms:
\begin{align*}
\mathbb{D}(\Sp_{\Z/\X}\F) &:= \mathbb{D}(\psi_t(j_{\neq 0!}\pr_1^*\F)) \\
					& \simeq \psi_t(\mathbb{D}(j_{\neq 0!}\pr_1^*\F[1]))[-1] \\
					& \simeq \psi_t(j_{\neq 0*}\pr_1^!(\mathbb{D}\F)[-1])[-1] \\
					& \simeq \psi_t(j_{\neq 0*}\pr_1^*[2](\mathbb{D}\F)[-1])[-1] \\
					& \simeq \psi_t(j_{\neq 0*}\pr_1^*(\mathbb{D}\F)) \\
					& \simeq \psi_t(j_{\neq 0!}\pr_1^*(\mathbb{D}\F)) \\
					& =: \Sp_{\Z/\X}(\mathbb{D}\F).
\end{align*}
\eproof

\subsubsection{Restriction to the base}
The following lemma computing
the restriction of the specialization
to the zero section, is indispensable
in the proof of
\cref{theorem: local statement}.


\blem
\label{lemma: restriction to base}
Suppose that $\X/k$ is a derived scheme 
locally of finite presentation (equivalently, locally of finite type). Let $\Z \hook \X$ be a
quasi-smooth, closed immersion. 
For any $\F \in \Shv(X)$,
\begin{align}
\label{sproperty1} \F|^*_{Z} &\xrightarrow{\simeq} (\Sp_{\Z/\X}\F)|^*_{Z} \\
\label{sproperty2} \F|^!_{Z} &\xleftarrow{\simeq} \Sp_{\Z/\X}(\F)|^!_{Z}.
\end{align}
\elem

\bproof
By \cref{lemma: dual of specialization}
it suffices to verify just either 
of \labelcref{sproperty1} or \labelcref{sproperty2}. We
verify \labelcref{sproperty1}.

Let $\widetilde{\pi}: \widetilde{D^*}_{\Z/\X}
\to D^*_{\Z/\X} \simeq \X \times \G_m$
denote the morphism indicated in diagram
\labelcref{diagram: definition of nearby cycles}.
Then,
\[\Sp_{\Z/\X}\F := \psi_t(j_{\neq 0!}\pr_1^*\F 
	\simeq s^*(\widetilde{\pi} \circ j_{\neq 0})_*
		(\widetilde{\pi} \circ j_{\neq 0})^*p^*\F,\]
where $p$ is the morphism indicated
in diagram \labelcref{diagram: classical
deformation 2}.
From this formula, we obtain
the following comparison map,
natural in $\F$,
\begin{align*}
\F|^*_{Z} 	&\simeq (s^*p^*\F)|^*_{Z} \\
			&\to (s^*(\widetilde{\pi} \circ j_{\neq 0})_*(\widetilde{\pi} \circ j_{\neq 0})^*p^*\F)|^*_{Z} \\
			&\simeq (\Sp_{\Z/\X}\F)|^*_{Z},
\end{align*}
where the first equivalence
follows from the fact that
the composition $Z \to 
Df_{\Z/\X} \to X$
is the inclusion $Z \hook X$.

We show that
this map is an isomorphism.
Note that \labelcref{sproperty1} is local
on $X^{\an}$. That is, it suffices to show
that $\Sp_{\Z/\X}(\F)|^*_{Z} \xrightarrow{\simeq} 
\F|^*_{Z}$ on an open cover of $X^{\an}$.
Since $\Z \hook \X$ is a quasi-smooth
closed immersion, by
\cref{lemma: local model} it
has the form
\[\begin{tikzcd}
		\Z & \X \\
		\0 & {\mathbb{A}^n}
		\arrow[from=2-1, to=2-2]
		\arrow["{\bm{f}}", from=1-2, to=2-2]
		\arrow[from=1-1, to=2-1]
		\arrow[from=1-1, to=1-2]
		\arrow["\lrcorner"{anchor=center, pos=0.125}, draw=none, from=1-1, to=2-2]
\end{tikzcd}\]
Zariski-locally on $\X$. Pick a
Zariski open cover of $\X$ on
each member of which $\Z \hook \X$
has the above form; then consider
the pullback of this cover along
the map $X \to \X$. Zariski
open covers are in particular
covers in the analytic topology
on $X$, so we may assume
without loss of generality
that $\Z \hook \X$ is the zero
locus of a finite
set of regular functions on $\X$,
$\bm{f} = (\bm{f}_1, \ldots, \bm{f}_n)$.

The derived scheme $\X$ 
embeds as a
closed subscheme of 
$\X \times \An$
via the graph map 
$\Gamma_{\bm{f}} :=
\id_{\X} \times \bm{f}: \X \hook
\X \times \An$. This map is,
in particular, proper, since
$\X \times \An$ is locally
Noetherian. Therefore,
\[\Sp_{(\X \times \0)/(\X \times \An)}
	(\Gamma_{f*}\F) \simeq
		N_{\Gamma_{\bm{f}}*}\Sp_{\Z/\X}(\F).\]
Restricting to $X \times \0$, we obtain
\begin{equation}
\label{equation: restriction to base}
\Sp_{(\X \times \0)/(\X \times \An)}
	(\Gamma_{f*}\F)|^*_{(X \times \0)} \simeq
		(Z \hook X)_*(\Sp_{\Z/\X}(\F)|^*_{Z}).
\end{equation}

\begin{claim}
There is an equivalence of functors,
	\[\Sp_{(\X \times \0)/(\X \times \An)} \simeq \Sp_{(X \times \0)/(X \times \An)},\]
where the 
left-hand side 
is quasi-smooth specialization
and the right-hand side is 
Verdier specialization.
\end{claim}

\bproof
By \cref{theorem: deformation} (i),
deformation to the normal bundle is stable
under arbitrary base change, meaning
\[\D_{(\X \times \0)/(\X \times \An)} \simeq 
	\D_{\0/\An} \times_{\An} (\X \times \An) \simeq
		D_{\0/\An} \times \X.\]
From this, we see clearly that
\[D_{(\X \times \0)/(\X \times \An)} 
	\simeq D_{(X \times \0)/(X \times \An)},\]
where the left-hand side is
the underlying classical scheme
of the derived deformation to the normal
bundle and the right-hand
side is the ordinary deformation
to the normal cone along a
closed subscheme. This suffices
to prove the claim, since both
quasi-smooth specialization
and Verdier specialization are
defined by the same formula.
\eproof

The claim in hand, we have
\[\Sp_{(\X \times \0)/(\X \times \An)}
	(\Gamma_{f*}\F)|^*_{(X \times \0)}
		\xrightarrow{\simeq} (\Gamma_{f*}\F)|^*_{(X \times \0)}\]
by \cite[\S8, Property (SP5)]{V83b}.
The base change isomorphism yields
\[(\Gamma_{f*}\F)|^*_{(X \times \0)} \simeq (Z \hook X)_*(\F|^*_{Z}),\]
which, combined with 
\labelcref{equation: restriction to base},
gives the result, since
$(Z \hook X)^*(Z \hook X)_* \simeq \id_{\Shv(Z)}$.
\eproof

\subsubsection{Specialization preserves perversity}

Finally, we show that quasi-smooth
specialization preserves perversity.


\blem
\label{lemma: sp and perversity}
Suppose that $\F \in \Perv(X)$. Then
$\Sp_{\Z/\X}(\F) \in \Perv(N_{\Z/\X})$.
\elem

\bproof
Recall that, with respect
to the perverse t-structure, 
\begin{enumerate}[(i)]
	\item $j_{\neq 0!}$ is right t-exact
		\cite[Theorem 5.2.4(iii)]{Dim04}.
	\item $j_{\neq 0*}$ is left t-exact
		\cite[Theorem 5.2.4(v)]{Dim04}.
	\item $\pr_1^*[1]$ is t-exact 
		\cite[Theorem 2.6(d)]{BBDJSS15}.   
\end{enumerate}
Items (i) and (iii) imply:
\begin{align*}
^p\psi_t(j_{\neq 0!}\pr_1^*[1](\F)) &= \psi_t(j_{\neq 0!}\pr_1^*[1](\F))[-1] \\
						&\simeq \psi_t(j_{\neq 0!}\pr_1^*\F) \\
						&=: \Sp_{\Z/\X}(\F) \in \Shv(N_{\Z/\X})^{\leq 0}.
\end{align*}
On the other hand, $\psi_t(j_{\neq 0!}\pr_1^*\F)
\simeq \psi_t(j_{\neq 0*}\pr_1^*\F)$, and the latter
lives in $\Shv(N_{\Z/\X})^{\geq 0}$ by items (ii)
and (iii). We conclude that 
$\Sp_{\Z/\X}(\F) \in \Perv(N_{\Z/\X})$.
\eproof


\subsection{Comparing classical notions
of specialization}

Earlier in this section we recalled
the notion of Verdier specialization,
a functor which takes a sheaf on
the analytic topology of a scheme $X$
and produces a sheaf on the normal
cone along a designated subscheme $Z$.

\subsubsection{}
When $X$ and $Z$ are smooth, their
analytifications have the structure of smooth
manifolds, and there is an alternative
notion of specialization available to us as
described in \cite[\S4.1 and \S4.2]{KS90}.
In \cref{sec: the local equivalence}, we
use a result from \cite[\S4]{KS90}
on the functor of specialization.
In order to access these
results, we sketch an argument 
that, in the case when
$Z \subset X$ is given as the zero
locus of a section of a vector bundle
on $X$, the notions from \cite{KS90}
and \cite{V83b} agree\footnotemark.

\footnotetext{It is folklore that these two
notions always agree, but we only verify
this agreement in the current setting.}

\begin{disc}

Suppose that $Z \subset X$ 
is given as the zero locus of 
the section $s: X \to E$ of 
a vector bundle $E \to X$.
By \cite[Remark 5.1.1]{F98}, 
$D_{Z/X}$ is then described 
by the following construction
attributed to MacPherson. 

Consider the map,
\[\begin{tikzcd}[row sep=0]
	{X \times \Aone_{\neq 0}} & {E \times \Aone_{\neq 0}} \\
	\rotatebox{90}{$\in$} & \rotatebox{90}{$\in$} \\ 
	{(x, t)} & {(\frac{s(x)}{t}, t)}
	\arrow["{S}", from=1-1, to=1-2]
	\arrow[maps to, from=3-1, to=3-2].
\end{tikzcd}\]
Then $D_{0_Z/V_Z}$ is 
given by the closure
of the image of $S$ 
inside of $E \times \Aone$.
Intuitively: $D_{0_Z/V_Z}$ 
is given by adding the ``limit 
point" at $0$ to the smooth 
$\Aone_{\neq 0}$\-family 
of $X$s whose slice at a point $t$ is
given by the embedded copy of 
$X$ inside $E$ given by $s$, 
scaled by $\frac{1}{t}$.

The diligent reader may now
check that the following
constructions are the identical: 
\begin{enumerate}
\item
\label{item: real def construction}
The restriction of the above 
construction to those values 
$t \in \R \subset \mathbb{C}$---
i.e. consider the closure of
the image of ${s|}_{V_Z \times \R_{\neq 0}}$
inside of $V_Z \times \R \times V$.
\item
\label{item: KS def construction}
The construction of deformation
to the normal bundle in \cite[\S4.1, pg. 186]{KS90}.
\end{enumerate}
In an open subset of $X$
over which $E$ is trivialized,
the correspondence 
goes as follows. In the 
notation of \cite[\S4.1]{KS90}:
$X \leftrightarrow U_i$, $E 
\leftrightarrow U_i \times \R^l$,
$s \leftrightarrow (\id_{U_i},
(\R^n \too \R^l \times \0) \circ \phi_i)$,
and $D_{Z/X} \leftrightarrow V_i$.

Now denote by $t: D_{Z/X}
\to \Aone$ and $pr_{\neq 0}: 
(X \times \Aone_{\neq 0} \simeq) 
D^*_{Z/X} \to X$ 
the restrictions of the 
obvious projections.

The diligent reader may
again check that $t^{-1}(\ell_{>0})$
coincides with the fiber
over $\R^+$ of
the object outlined
in (\labelcref{item: real def construction}),
as well as the subset
in \cite[\S4.1]{KS90} that
the authors denote by $\Omega$. 
Clearly, $\Sp_{Z/X}\F
\simeq \psi_t(\pr_{\neq 0}^*\F)$.
At the same time, formula 
\labelcref{nearby cycles at fixed angle no.2}
for nearby cycles
at fixed angle $\theta = 0$ yields:
\[\psi_t(\pr_{\neq 0}^*\F) \simeq
	(t^{-1}(0) \hook D_{Z/X})^*
		(t^{-1}(\ell_{>0}) \hook D_{Z/X})_*
			(t^{-1}(\ell_{>0}) \hook D^*_{Z/X})^*
				\pr_{\neq 0}^*\F,\]
which coincides precisely
with the formula for the
specialization in \cite[\S4.2]{KS90},
taking into account the agreement
of constructions 
(\labelcref{item: real def construction})
and (\labelcref{item: KS def construction})
sketched above.
\end{disc}

\subsubsection{Specialization of monodromic sheaves}
MacPherson's description
of deformation to the normal
cone allows for a nice proof
of the following
lemma, which is
well-known but hard to
find in the literature.

\blem
\label{lem: specialization of monodromic sheaf is identity}
Suppose that $\bm{E} \to \X$ is a
finite dimensional complex 
vector bundle on the 
derived scheme $\X$. 
If $\F \in \Shv_{\G_m}(\bm{E})$, then
\begin{align*}
\Sp_{\X/\bm{E}}\F 	&\simeq \F \\ 
\Sp_{X/E}\F  	&\simeq \F,
\end{align*}
where $\X$ (resp. $X$) is seen as the
zero section of $\bm{E}$
(resp.\ $E$).
\elem

\bproof
The following equivalence
of schemes is easily seen
from the properties of
quasi-smooth blow-up:
\[D_{\X/\bm{E}} \simeq D_{X/E}.\]
As such, $\Sp_{\X/\bm{E}} \simeq
\Sp_{X/E}$, so it suffices to prove
that $\Sp_{X/E}\F \simeq \F$.
Without loss of generality,
we may assume $X = \ast$
and $E$ is a vector space $V$.

MacPherson's description
of deformation to the normal
cone, outlined in 
\cite[Remark 5.1.1]{F98},
identifies $D_{0/V}$ with the
closure of the image of
the map,
\[\begin{tikzcd}[row sep=0]
	{V \times \Aone_{\neq 0}} & {V \times V \times \Aone} \\
	\rotatebox{90}{$\in$} & \rotatebox{90}{$\in$} \\
	{(v, t)} & {(v, \frac{v}{t}, t)}
	\arrow["{S}", from=1-1, to=1-2]
	\arrow[maps to, from=3-1, to=3-2],
\end{tikzcd}\]
and identifies the maps 
$t: D_{0/V} \to \Aone$
and $p: D_{0/V} \to V$ with
the restrictions to 
$D_{0/V}$ of $\pr_3: V
\times V \times \Aone \to \Aone$
and $\pr_1: V \times V \times \Aone \to V$,
respectively.
Let $U := D_{0/V} \setminus  
(0,0) \times \Aone$.
By smooth base change,
\[(\Sp_{0/V}\F|)_{V \setminus 0} 
	\simeq \psi_{t_U}(p_U^*\F)|_{V \setminus 0},\]
where $t_U$ and $p_U$ are
the restriction to $U$ of $t$
and $p$, respectively. We show
that $\psi_{t_U}(p_U^*\F)|_{
V \setminus 0} \simeq \F_{V
\setminus 0}$.

Let $\pi: V \setminus 0 \to
(V \setminus 0)/{\R^+} = 
S(V)$ be the map to the spherization of $V$.
The adjoint functors 
\[\pi^*: \Shv(S(V) \leftrightarrows 
	\Shv_{\R^+}(V \setminus 0): \pi_*\]
induce an equivalence of categories
because the counit
map is an isomorphism.
While the map $S$ doesn't
extend to a continuous map
on $V \times \Aone$, the induced map,
\[\begin{tikzcd}[row sep=0]
	{S(V) \times \Aone_{\neq 0}} & {S(V \times V) \times \Aone} \\
	\rotatebox{90}{$\in$} & \rotatebox{90}{$\in$} \\
	{([v], t)} & {([v, \frac{v}{t}], t)} 
	\arrow["{\widetilde{S}}", from=1-1, to=1-2]
	\arrow[maps to, from=3-1, to=3-2]
\end{tikzcd}\]
\emph{does} extend to a continuous
map on $S(V) \times \Aone$.
This extension is a 
homeomorphism
onto its image, which
we can identify with
$S(D_{0/V})$, the fiberwise 
spherization of $D_{0/V}$.
In all, we obtain the following
commutative diagram, whose
maps are induced 
from maps on $D_{0/V}$:
\[\begin{tikzcd}
	& {S(V)} \\
	{S(V) \times \Aone} & {S(D_{0/V})} \\
	& \Aone
	\arrow["{\widetilde{t}}", from=2-2, to=3-2]
	\arrow["{\widetilde{p}}"', from=2-2, to=1-2]
	\arrow["{\pr_2}"', from=2-1, to=3-2]
	\arrow["{\pr_1}", from=2-1, to=1-2]
	\arrow["{\widetilde{S}}", from=2-1, to=2-2]
\end{tikzcd}\]
Denote $\widetilde{F} :=
\pi_*(\F|_{V \setminus 0})$
From the above diagram:
\begin{align*}
\widetilde{\F}			&\simeq \psi_{\pr_2}(\pr_1^*\pi_*\widetilde{\F}) \\
					&\simeq \psi_{\widetilde{t} \circ \widetilde{S}}(\widetilde{S}^*\widetilde{p}^*\pi_*\widetilde{\F}) \\
					&\xleftarrow{\simeq} \widetilde{S}|_{\widetilde{t}^{-1}(0)}^*\psi_{\widetilde{t}}(p^*\pi_*\widetilde{\F}) \\
					&\simeq \widetilde{\Sp_{0/V}\F}
\end{align*}
Because $\pi_*$ induces
an equivalence of categories,
we obtain an equivalence,
$\F|_{V \setminus 0} \simeq 
(\Sp_{0/V}\F)|_{V \setminus 0}$. 
Combining this equivalence
with that obtained in
\cref{lemma: restriction to base},
we use an adjunction triangle
for the complementary pair
of subspaces
$0, V \setminus 0 \subset V$
to obtain the desired equivalence,
$\F \simeq \Sp_{0/V}\F$.
\eproof

The following is an immediate
corollary of \cref{lem:
specialization of monodromic sheaf
is identity} which is also
well-known but hard to find in
the literature.

\begin{cor}
\label{cor: specialization is idempotent}
Both Verdier and quasi-smooth
specialization are idempotent 
in the following sense. 
Suppose $Z \hook X$ 
an arbitrary closed
immersion of ordinary schemes. 
Then, for
$\F \in \Shv(X)$, 
\[\Sp_{Z/\N_{Z/X}}(\Sp_{Z/X}\F) 
	\simeq \Sp_{Z/X}\F,\]
after making the identification
$\N_{Z/N_{Z/X}} \simeq \N_{Z/X}$.
If $\Z \hook \X$ is a
quasi-smooth embedding
of derived schemes, the
corresponding statement
is obtained by substituting
$\bm{\N}, \Z,$ and $\X$
for $\N, Z$, and $X$.
\end{cor}

\subsubsection{}
Lastly, we record a result from \cite{KS90} that
will be used in \cref{sec: the local equivalence}, 
as well as a definition needed to state it.

\bdef(\cite[Definition 4.1.1(i)]{KS90})
\label{def: normal cone along Z}
Suppose that $p: D_{Z/X} \to X$
is the canonical map, and $p_{>0}:
t^{-1}(\Aone_{>0}) =: D_{Z/X, >0} 
\to X$ is its restriction. Suppose given
a subset $S \subset X$. Then we
define the normal cone
of $S$ along $Z$ as:
\[C_Z(S) := \overline{p_{>0}^{-1}(S)} \cap C_{Z/X}.\]
\edefn 

\bthm(\cite[Theorem 4.2.3(ii)]{KS90})
\label{thm: sections of sp on opens}
Suppose $M$ is a smooth
submanifold of a smooth manifold $X$.
Let $\F \in \Shv(X)$ and $\V$ be an 
open conic subset of $N_{M/X}$.
Then:
\begin{equation}
\label{equation: sections of sp on opens}
\varinjlim_{\U} \Gamma(\U; \F) 
	\xrightarrow{\simeq} \Gamma(\V;\Sp_{M/X}\F),
\end{equation}
where $\U$ ranges through the
family of open subsets of $X$
such that $C_M(X \setminus \U) 
\cap \V = \varnothing$, which
we denote by $\mathfrak{F}_{\V}$.
\ethm 


\subsection{Quasi-smooth microlocalization}
\subsubsection{}
We make the following 
obvious definition.

\bdef
Suppose that $\Z \hook \X$ is a quasi-smooth
closed immersion, and that $\F \in \Shv(X)$. 
The \emph{microlocalization} of $\F$ along $\Z
\hook \X$ is defined as
	\[\mu_{\Z/\X}\F := \Four(\Sp_{\Z/\X}\F)\]
\edefn

Clearly, $\mu_{\Z/\X}$ is a functor from
$\Shv(X)$ to $\Shv_{\G_m}(N^{\vee}_{\Z/\X})$.
We also call $\mu_{\Z/\X}$ the functor of
\emph{quasi-smooth microlocalization}.

\subsubsection{}
Unlike quasi-smooth specialization, $\mu_{\Z/\X}$
preserves perversity only up to a shift. The following
is an immediate corollary of
\cref{lemma: sp and perversity} and 
\cref{lemma: FST and perversity}.

\blem
\label{lemma: mu and perversity}
$\mu_{\Z/\X}[\vircodim \, \Z]: \Perv(X) \to \Perv(N^{\vee}_{\Z/\X})$.
\elem

\brem
Just as quasi-smooth specialization
enjoys most of the formal properties
of its classical counterpart, quasi-smooth
microlocalization enjoys many more
properties than the one stated
in \cref{lemma: mu and perversity}
above. Listing and proving them
is beyond the scope of this paper.
\erem


\section{The local equivalence}
\label{sec: the local equivalence}


\subsection{Stating the equivalence}
For the remainder of this
section, we use the following
set-up. Let $X$ be a 
smooth $k$-scheme, 
and let $V$ be a $k$-vector 
space of dimension $n$, 
viewed as a $k$-scheme. 
Given a map $f: X \to V$, we consider
the derived zero locus,
\[\begin{tikzcd}
	{\Z(f)} & X \\
	\0 & V
	\arrow[hook, from=2-1, to=2-2]
	\arrow[hook, from=1-1, to=1-2]
	\arrow["f", from=1-2, to=2-2]
	\arrow[from=1-1, to=2-1]
	\arrow["\lrcorner"{anchor=center, pos=0.125}, draw=none, from=1-1, to=2-2].
\end{tikzcd}\]
The closed immersion $\Z(f) 
\hook X$ is obviously quasi-smooth,
and by application of base change 
for cotangent complexes, we find that
	\[\bm{N}_{\Z(f)/X} \simeq \Z(f) \times V.\]
For such a map, $f$, we consider 
the map $\widetilde{f}: X \times V^{\vee} 
\to V$ whose definition is indicated 
by the following diagram:
	\[\begin{tikzcd}
		{X \times V^{\vee}} & {V\times V^{\vee}} \\
		{} & {\mathbb{A}^1}
		\arrow["{\langle \cdot \,, \cdot \rangle}", from=1-2, to=2-2]
		\arrow["{f \times \id_{V^{\vee}}}", from=1-1, to=1-2]
		\arrow["{\widetilde{f}}"', from=1-1, to=2-2],
	\end{tikzcd}\]
where $\langle \cdot \,, \cdot \rangle$ is the natural 
pairing between a vector space and its dual.

\subsubsection{}
\label{sssec: identify conormal with product}
The conormal bundle $N^{\vee}_{\Z(f)/X}$
embeds into $X \times V^{\vee}$ as the 
closed subscheme $Z(f) \times V^{\vee}$,
using the identification $N^{\vee}_{\Z(f)/X}
\simeq Z(f) \times V^{\vee}$ given above.  
Note in particular that, under
the inclusion $Z(f) \times V^{\vee}
\hook Z(\widetilde{f})$, $N^{\vee}_{\Z(f)/X} 
\subset Z(\widetilde{f})$, the (underived)
zero locus of $\widetilde{f}$.

\subsubsection{}
Denote by $\pr_1$ and $\pr_2$ the canonical 
projections from $X \times V$ onto its factors.
In this section, we prove the following theorem.

\begin{thm}
\label{theorem: local statement}
 There exists a natural isomorphism of functors,
	\[\varphi_{\widetilde{f}}(\pr_1^*-) \simeq \mu_{\Z(f)/X}(-).\]
\end{thm}

\subsubsection{}
Prima facie, \cref{theorem: 
local statement} is ill-posed: for a given
$\F \in \Shv(X)$, the left-hand side 
and right-hand side produce sheaves 
that live on different spaces. However, 
$\varphi_{\widetilde{f}}(\pr_1^*\F)$ has support contained
in $Z(f) \times V^{\vee}$ by \cite[Lemma A.4]{Dav17}, 
so it can be viewed as a sheaf on such. At the same
time, the identification discussed in \cref{sssec:
identify conormal with product} allows us
to view $\mu_{\Z(f)/X}(\F)$ as an element
of $\Shv(Z(f) \times V^{\vee})$, as well, so the two
functors of \cref{theorem: local statement}
are indeed comparable.


\subsection{The monodromic case}

The following lemma is the
basis for our strategy of proof.

\blem
\label{lem: phi is four}
Suppose that $p: E \to S$ is a complex vector bundle
of finite rank. If $\F \in \Shv_{\G_m}(E)$, then 
\[\mathscr{F}(\F) \simeq \varphi_m(\pr_1^*\F),\]
where $\pr_1: E \oplus E^{\vee} \to E$ is the obvious
projection and $m: E \oplus E^{\vee} \to \Aone$ is the natural
pairing map.  
\elem

\bproof
The key observation in this proof is
the following identity:
\begin{align*}
P &:= \{(v,w) \in E \oplus E^{\vee}| \Real w(v) \geq 0\} \\
   &=: \{\Real m \geq 0\}
\end{align*}
Recalling the expression
for vanishing cycles at the fixed
angle $\theta = \pi$,
\labelcref{vanishing cycles at 
fixed angle no.1}, we obtain 
the following chain of isomorphisms:
\begin{align*}
\Four(\F) &:= \pr_{2*}\Gamma_{P}(\pr_1^*\F) \\
		&:= \pr_{2*}\Gamma_{\{\Real m \geq 0\}}(\pr_1^*\F) \\
		&\simeq (0 \oplus E^{\vee} \hook E \oplus E^{\vee})^*\Gamma_{\{\Real m \geq 0\}}(\pr_1^*\F) \\
		&\simeq (0 \oplus E^{\vee} \hook m^{-1}(0))^*(m^{-1}(0) \hook E \oplus E^{\vee})^*\Gamma_{\{\Real m \geq 0\}}(\pr_1^*\F) \\
		&\simeq (0 \oplus E^{\vee} \hook m^{-1}(0))^*\varphi_m(\pr_1^*\F),
\end{align*}
where the third equivalence holds by
monodromicity with respect to the $\G_m$\-action
given by scaling on the first factor alone.

Now observe that $\varphi_m(\pr_1^*\F)$
is supported on the closed
subset $0 \oplus E^{\vee}$ by 
\cite[Lemma A.4]{Dav17}, 
so we may identify
$\varphi_m(\pr_1^*\F)$ with its pullback
$(0 \oplus E^{\vee} \hook m^{-1}(0))^*\varphi_m(\pr_1^*\F)$,
which completes the proof.
\eproof

\subsubsection{}
Suppose that in the 
statement of 
\cref{lem: phi is four} we
took $E \to S$ to be the normal bundle
$N_{\Z(f)/X} \to Z(f)$ and $\F$ to be
the specialization $\Sp_{\Z(f)/X}\F$.
Then we obtain,
\[\mu_{\Z(f)/X}(\F) := \Four(\Sp_{\Z(f)/X}\F)
	\simeq \varphi_m(\pr_1^*\Sp_{\Z(f)/X}\F).\]
Therefore, it suffices to prove
that the functor $\varphi_m(\pr_1^*-)$
factors through the specialization
along $\Z(f)$ in the sense of the
following lemma.
\blem
\label{lemma: phi factors through sp}
\begin{equation*}
\varphi_{\widetilde{f}}(\pr_1^*\F) \simeq 
	\varphi_m(\pr_1^*\Sp_{\Z(f)/X}\F)
\end{equation*}
\elem

\subsubsection{Reducing to the case of a product}
\label{sssec: graph embedding}
In order to prove
\cref{lemma: phi factors
through sp}, it suffices to 
take $X := Z \times V$,
for $Z$ smooth, V a finite dimensional
vector space over $k$, and $f := \pr_V$.
Demonstrating this is
a matter of simple, repeated
application of proper base 
change; we write it out 
for completeness.

Consider the following diagram,
\begin{equation}
\label{diagram: graph embedding}
\begin{tikzcd}
	& V \\
	X & {X \times V} \\
	{X \times V^{\vee}} & {X \times V \times V^{\vee}} \\
	& \Aone
	\arrow["{\Gamma_f \times \id_{V^{\vee}}}", from=3-1, to=3-2]
	\arrow["{\widetilde{\pr}_V}", from=3-2, to=4-2]
	\arrow["{\widetilde{f}}"', from=3-1, to=4-2]
	\arrow["{\pr_1}"', from=3-2, to=2-2]
	\arrow["{\pr_1}"', from=3-1, to=2-1]
	\arrow["{\Gamma_f}", from=2-1, to=2-2]
	\arrow["{\pr_V}"', from=2-2, to=1-2]
	\arrow["f", from=2-1, to=1-2],
\end{tikzcd}
\end{equation}
where the meaning of the
overloaded notation ``$\pr_1$" 
is inferred from context in
what follows.

Using the graph map
$\Gamma_f: X \hook X \times V$, we apply proper base 
change for vanishing cycles to obtain:
\begin{align*}
(Z(\wt{f}) \to Z(\wt{\pr}_V))_*\varphi_{\wt{f}}(\pr_1^*\F) &\simeq \varphi_{\wt{\pr}_V}((\Gamma_f \times \id_{V^{\vee}})_*\pr_1^*\F) \\
											& \simeq \varphi_{\wt{\pr}_V}(\pr_1^*\Gamma_{f*}\F)
\end{align*}
One the other hand, since $f = \pr_V \circ \Gamma_f$,
the square,
\[\begin{tikzcd}
	{\Z(f)} & {X \times \0 = \Z(\pr_V)} \\
	X & {X \times V}
	\arrow[from=1-1, to=2-1]
	\arrow[from=1-2, to=2-2]
	\arrow[hook, from=1-1, to=1-2]
	\arrow["{\Gamma_f}", hook, from=2-1, to=2-2]
	\arrow["\lrcorner"{anchor=center, pos=0.125}, draw=none, from=1-1, to=2-2],
\end{tikzcd}\]
is derived Cartesian. Note that
$N_{(X \times \0)/(X \times V)}
\simeq X \times V$, and $N_{\pr_V}
= \pr_V$. Meanwhile, $N_{\Z(f)/X}
\simeq Z(f) \times V$, and $N_{f}:
Z(f) \times V \to V$ is the projection.
The following diagram
is induced from diagram
\labelcref{diagram: graph embedding}:
\[\begin{tikzcd}
	\Aone \\
	{N_{\Z(f)/X} \times V^{\vee}} & {N_{(X \times\0)/(X \times V)} \times V^{\vee}} \\
	{Z(f) \times V \simeq N_{\Z(f)/X}} & {N_{(X \times\0)/(X \times V)} \simeq X \times V}
	\arrow["{\pr_1}", from=2-1, to=3-1]
	\arrow["{\pr_1}", from=2-2, to=3-2]
	\arrow["{N_{\Gamma_f} \simeq i_{Z(f)} \times \id_V}"{pos=0.7}, from=3-1, to=3-2]
	\arrow["{N_{\Gamma_f} \times \id_{V^{\vee}}}", from=2-1, to=2-2]
	\arrow["m"', from=2-1, to=1-1]
	\arrow["{\widetilde{N_{\pr_V}} = \widetilde{pr}_V}"', from=2-2, to=1-1],
\end{tikzcd}\]
where $m$ is as follows:
\[\begin{tikzcd}[row sep=0]
	{N_{\Z(f)/X} \times V^{\vee}} & {\Aone} \\
	\rotatebox{90}{$\in$} & \rotatebox{90}{$\in$} \\
	{(z, v, w)} & {w(v)} 
	\arrow["{m}", from=1-1, to=1-2]
	\arrow[maps to, from=3-1, to=3-2].
\end{tikzcd}\]
Proper base change
for quasi-smooth specialization 
(\cref{lemma: proper base change for 
specialization}) and vanishing cycles now yields:
\begin{align*}
(Z(m) \to Z(\wt{\pr}_V))_*\varphi_m(\pr_1^*\Sp_{\Z(f)/X}\F) &\simeq \varphi_{\wt{\pr}_V}((N_{\Gamma_f} \times \id_{V^{\vee}})_*\pr_1^*\Sp_{\Z(f)/X}\F) \\
											&\simeq \varphi_{\wt{\pr}_V}(\pr_1^*(N_{\Gamma_{f}*}\Sp_{\Z(f)/X}\F) \\
											&\simeq \varphi_{\wt{\pr}_V}(\pr_1^*\Sp_{(X \times \0)/(X \times V)}(\Gamma_{f*}\F)).
\end{align*}
Hence, if the equivalence
\[\varphi_{\wt{\pr}_V}(\pr_1^*\mathcal{G}) 
	\simeq \varphi_{\wt{\pr}_V}(\pr_1^*
		\Sp_{(X \times \0)/(X \times V)}(\mathcal{G})),\]
holds for all $\mathcal{G} \in \Shv(X \times V)$,
the equivalence
\[(Z(\wt{f}) \to Z(\wt{\pr}_V))_*
	\varphi_{\wt{f}}(\pr_1^*\F) 
		\simeq (Z(m) \to Z(\wt{\pr}_V))_*
			\varphi_m(\pr_1^*\Sp_{\Z(f)/X}\F),\]
also holds. The latter then 
implies
\cref{lemma: phi factors through sp}
once we take into account
that each side is supported
on $Z \times \0 \times V^{\vee}
\subset X \times V \times V^{\vee}$.

\subsubsection{}
For the remainder
of this section,
take $Z$ smooth, 
$X := Z \times V$, and $f := \pr_V$,
per \cref{sssec: graph embedding}
above.

\subsubsection{}
To ease the exposition, we
denote $Z \times V$ 
by $V_Z$, viewing $X$ as a scheme 
over $Z$ with structure map given 
by $\pr_1$. Thus, products are taken
in the overcategory $\Sch_{/Z}$, 
meaning $V_Z \times_Z V^{\vee}_Z$ denotes 
$Z \times V \times V^{\vee}= X \times V^{\vee}$,
$0_Z \subset V_Z$ denotes
$Z \times \0 \subset Z \times V$, etc.
By abuse of notation, we denote the 
analytifications of all these schemes
by the same. From hereon out, 
we denote $\widetilde{\pr_V}$ by $m$.
$\pr_1$ will always denote the projection
$V_Z \times_Z V^{\vee}_Z \to V_Z$.


\subsection{Proof strategy for the non-monodromic case}

Our proof of \cref{lemma:
phi factors through sp}
proceeds indirectly by first
proving the equivalence
\[\psi_m(\pr_1^*\F)|_{V^{\vee}_Z (= (0 \times V^{\vee})_Z)} 
\xrightarrow{\simeq} \psi_m(\pr_1^*\Sp_{0_Z/V_Z}\F)|_{V^{\vee}_Z},\]
and then using the standard fiber
sequence relating nearby and vanishing cycles. 
In order to prove this
latter equivalence, we
define a natural comparison map, 
$\psi: \psi_m(\pr_1^*\F)|_{V^{\vee}_Z} 
\to \psi_m(\pr_1^*\Sp_{0_Z/V_Z}\F)|_{V^{\vee}_Z}$,
which we then show is an equivalence by checking
on stalks. Since obtaining the comparison
map is somewhat involved,
we briefly outline its construction below.

\subsubsection{Outline of comparison map construction}
Observe that that 
$\varphi_m(\pr_1^*\F)$
is conic by
\cref{lemma: monodromic},
so $\psi_m(\pr_1^*\F)
|_{V^{\vee}_Z}$ is conic via
the canonical fiber sequence
\labelcref{fiber sequence: 
nearby/vanishing cycles}.
Let $\B_{\R^+}(V^{\vee}_Z) 
\subset \U_{\R^+}(V^{\vee}_Z)$
denote the partially ordered
subset consisting of
convex, conic open subsets. 
These are a basis for the conic
topology on $V^{\vee}_Z$.

We begin the construction of 
the comparison map, $\psi$,
by finding convenient expressions
for the sections of
$\psi_m(\pr_1^*\F)|_{V^{\vee}_Z}$
and $\psi_m(\pr_1^*\Sp_{0_Z/V_Z}\F)|_{V^{\vee}_Z}$
over each element
of $\B_{\R^+}(V^{\vee}_Z)$.
These convenient expressions
are colimits of sections of $\F$
over certain families of open subsets in $V_Z$. 
From these descriptions, we argue
that the restrictions of $\psi_m(\pr_1^*\F)|_{V^{\vee}_Z}$
and $\psi_m(\pr_1^*\Sp_{0_Z/V_Z}\F)|_{V^{\vee}_Z}$
to $\B_{\R^+}(V^{\vee}_Z)$
are actually left Kan extensions of the
functors (to be defined), $Q \circ \alpha$ 
and $Q$, respectively.
Standard abstract nonsense gives
a map between these two Kan extensions,
which are sheaves on
$\B_{\R^+}(V^{\vee}_Z)$. 
We then verify that the category 
of conic sheaves on 
$V^{\vee}_Z$ is equivalent to
the category of sheaves on
the basis $\B_{\R^+}(V^{\vee}_Z)$,
thereby obtaining the desired 
comparison map.

\begin{notn}
To ease the exposition, we sometimes
choose to denote
$\psi_m(\pr_1^*\F)|_{V^{\vee}_Z}$
and $\psi_m(\pr_1^*\Sp_{0_Z/V_Z}\F)|_{V^{\vee}_Z}$
by $\psi_m$ and $\psi_m^{\Sp}$, respectively,
for the remainder of this section.
\end{notn}


\subsection{Obtaining expressions for sections of $\psi_m$ and $\psi_m^{\Sp}$}
\label{ssec: obtaining expressions for sections}
\subsubsection{Expression for sections of $\psi_m$}

Let $j_{>0}: m^{-1}(\Aone_{>0}) 
\hook V_Z \times_Z V^{\vee}_Z$
and $i_{V^{\vee}_Z}: V^{\vee}_Z 
\hook V_Z \times_Z V^{\vee}_Z$
denote the obvious inclusions.
Using the expression 
\labelcref{nearby cycles at fixed angle no.1}
for nearby cycles at the 
fixed angle $\theta = \pi$,
we obtain:
\[\psi_m(\pr_1^*\F)|_{V^{\vee}_Z} \simeq 
	i_{V^{\vee}_Z}^*j_{>0*}j_{>0}^*\pr_1^*\F. \]
Now choose a convex, conic open subset
$\V \subset V^{\vee}_Z$ (see
\cref{rem: clarifying definition of convex}).
We compute the following
chain of isomorphisms:
\begin{align}
\Gamma(\V; \psi_m(\pr_1^*\F)|_{V^{\vee}_Z}) 	&\xleftarrow{\simeq} \varinjlim_{0_Z \times \V \subset \U} \Gamma(\U; j_{>0*}j_{>0}^*\pr_1^*\F) \\
									&\xleftarrow{\simeq} \varinjlim_{\epsilon>0} \Gamma(B_{\epsilon}(0_Z) \times_Z \V; j_{>0*}j_{>0}^*\pr_1^*\F) \\
									&\xrightarrow{\simeq} \varinjlim_{\epsilon>0} \Gamma(B_{\epsilon}(0_Z) \times_Z \V \cap m^{-1}(\Aone_{>0}); j_{>0}^*\pr_1^*\F), \label{last term in align}
\end{align}
where, in the second line,
$B_{\epsilon}(0_Z)$ is the
subset $Z \times B_{\epsilon}(0)$,
determined by fixing an isomorphism
of $V$ with $\R^{\dim_{\R}V}$.

Now consider the diagram,
\[\begin{tikzcd}
	{(B_{\epsilon}(0_Z) \times_Z \V) \cap m^{-1}(\Aone_{>0})} & {V_Z \times_Z V^{\vee}_Z} \\
	{\U_{\epsilon, \V} :=\pr_1((B_{\epsilon}(0_Z) \times_Z \V) \cap m^{-1}(\Aone_{>0}))} & {V_Z} \\
	\ast
	\arrow["{j_{\epsilon, \V}}", from=1-1, to=1-2]
	\arrow["{\pr_1}", from=1-2, to=2-2]
	\arrow["{j'_{\epsilon, \V}}"', from=2-1, to=2-2]
	\arrow["{p_{\epsilon, \V}}"', from=1-1, to=2-1]
	\arrow["pt"', from=2-1, to=3-1].
\end{tikzcd}\]

In the notation of
this diagram, we have:
\begin{align}
\Gamma(B_{\epsilon}(0_Z) \times_Z \V \cap m^{-1}(\Aone_{>0}); j_{>0}^*\pr_1^*\F) 	&\simeq pt_*j_{\epsilon, \V}^*\pr_1^*\F \\
																&\simeq pt_*p_{\epsilon, \V *}p_{\epsilon, \V}^*{j'_{\epsilon, \V}}^*\F \\
																&\xleftarrow{} pt_*{j'_{\epsilon, \V}}^*\F \\
																&=: \Gamma(\U_{\epsilon, \V}; \F).
\end{align}
We use the following lemma
to prove that the arrow
in the third line is an isomorphism.

\blem
\label{lemma: p_*p^* = id} 
Suppose $E \xrightarrow{p} X$ 
is a continuous map
of topological spaces
such that $E \subset W$ 
where $W$ is a topological
finite dimensional,
real vector bundle over $X$.
Suppose also that the fibers
of $p$ are nonempty, convex subsets
of $W$ and that there exist continuous
local sections of $p$ around each
$x \in X$. Then $\id 
\xrightarrow{\simeq} p_*p^*$.
\elem

\bproof
The statement that $\id \to 
p_*p^*$ is an isomorphism
can be checked on stalks,
so without loss of generality,
we assume $p$ actually has a
global section, $X \xrightarrow{s} E$.

By \cite[Corollary 2.7.7(ii)]{KS90} it
suffices to show that $s \circ p$ is
homotopic to $\id_E: E \to E$ over
$X$. That is, we exhibit a continuous
map $H: E \times [0,1] \to E$ such that
$H(-,1) = s \circ p$, $H(-,0) = \id_E$, 
and $p \circ \pr_1 = p \circ H$.
Such a map is given by
\[H(e,t) := t \cdot s(p(e)) + (1-t) \cdot e,\]
where the scalar
multiplication, $\cdot$, is 
performed in $W$.
Since the fiber $E_x$ is
convex, for all $x \in X$,
$H(t,e) \in E_{p(e)}$ for
all $t \in [0,1]$.
\eproof

\begin{claim}
$p_{\epsilon, \V}$ satisfies the hypotheses of \cref{lemma: p_*p^* = id}.
\end{claim}

\bproof
Fix a point $(z_0, v_0) \in \U_{\epsilon, \V}$.
The fiber $p_{\epsilon, \V}^{-1}(z_0, v_0)$
is the intersection
\[\{(v_0, \lambda) \in V_Z \times_{z_0} V^{\vee}_Z| 
	\Real \lambda(v_0) > 0\} \cap (\V|_{z_0}).\]
It is nonempty by the definition of
$\U_{\epsilon, \V}$.
Note that $\Real \lambda(z_0, v_0): V^{\vee}_Z \to \R_Z$
is map of real vector bundles
over $Z$, so $\{(v_0, \lambda) \in 
V_Z \times_{z_0} V^{\vee}_Z| 
\Real \lambda(v_0) > 0\}$ is the
product of an open half-plane
in $V^{\vee}$ with $Z$; in particular,
it is a convex, conic open subset
of $V^{\vee}_Z$. On the other hand,
$\V$ is a convex, conic open subset
of $V^{\vee}_Z$ by assumption.
By elementary convex geometry,
therefore, their intersection
is also a convex, conic open
subset of $\{(z_0, v_0)\} \times V^{\vee}$.

Now, choosing an arbitrary 
$\lambda_0 \in 
p_{\epsilon, \V}^{-1}(z_0, v_0)$,
it is clear from the continuity
of $m$ that $(z, v, \lambda_0) 
\in p_{\epsilon, \V}^{-1}(z, v)$
for all $(z,v)$ in a sufficiently 
small neighborhood of $z_0, v_0$.
This shows that $p_{\epsilon, \V}$
has local sections (i.e.\ $s(z,v) = 
\lambda_0)$), so the hypotheses of
\cref{lemma: p_*p^* = id} are satisfied.
\eproof

Thus, we obtain the equivalence,
\begin{align*}
\labelcref{last term in align} &\xleftarrow{\simeq} \varinjlim_{\epsilon>0} \Gamma(\U_{\epsilon, \V}; \F).
\end{align*}

\subsubsection{Expression for sections of $\psi_m^{\Sp}$}

Observe that $\U_{\epsilon, \V}$ is
itself open (as $\pr_1$ is an open
map) and, not quite conic, but rather the 
intersection of the conic
open subset $\U_{\V} := \pr_1(V \times \V \cap 
m^{-1}(\Aone_{>0}))$ 
and the open ball
$B_{\epsilon}(0_Z)$. Thus, in light
of the above discussion,
\begin{align}
\Gamma(\V; \psi_m(\pr_1^*\Sp_{0_Z/V_Z}\F)|_{V^{\vee}_Z}) 	&\simeq \varinjlim_{\epsilon>0} \Gamma(\U_{\epsilon, \V}; \Sp_{0_Z/V_Z}\F) \\
											&\xleftarrow{\simeq} \varinjlim_{\epsilon>0} \Gamma(\U_{\V}; \Sp_{0_Z/V_Z}\F) \\
											&\simeq \Gamma(\U_{\V}; \Sp_{0_Z/V_Z}\F),
\end{align}
where we have used that
the restriction of sections
along $\U_{\epsilon, \V} \hook
\U_{\V}$ is an equivalence
because $\Sp_{0_Z/V_Z}\F$
is conic.

Finally, we take advantage of
\labelcref{equation: sections of sp
on opens} to write:
\[\Gamma(\U_{\V}; \Sp_{0_Z/V_Z}\F) 
	\xleftarrow{\simeq} \varinjlim_{\U \in 
		\mathfrak{F}_{\U_{\V}}} \Gamma(\U; \F).\]

\brem
\label{rem: defining psi_V}
Note that $\{\U_{\epsilon, \V}\}_{\epsilon>0}$
is a subcollection of $\mathfrak{F}_{\U_{\V}}$,
so there is a natural map of
sections over $\V$,
which we denote by $\psi_{\V}$,
coming from the composition of the following
chain of morphisms:
\begin{align*}
\Gamma(\V; \psi_m(\pr_1^*\F)|_{V^{\vee}_Z}) 	&\simeq \varinjlim_{\epsilon>0} \Gamma(\U_{\epsilon, \V}; \F) \\
									&\to \varinjlim_{\U \in \mathfrak{F}_{\U_{\V}}} \Gamma(\U; \F) \\
									&\simeq \Gamma(\U_{\V}; \Sp_{0_Z/V_Z}\F) \\
									&\simeq \Gamma(\V; \psi_m(\pr_1^*\Sp_{0_Z/V_Z}\F)|_{V^{\vee}_Z}).
\end{align*}
\erem


\subsection{Obtaining the comparison map}
\label{ssec: obtaining the comparison map}

\subsubsection{Left Kan extensions}

We briefly recall the notion of left Kan extension
found in \cite[\S4.3]{HTT}. First suppose
that $i: \mathcal{C} \to \widetilde{\mathcal{C}}$ is
a fully faithful embedding, and let $F: 
\mathcal{C} \to \mathcal{D}$ be an
arbitrary functor. In the definition
below, for $\widetilde{x} \in \ob(\widetilde{\mathcal{C}})$:
$\mathcal{C}_{/\widetilde{x}} := 
\mathcal{C} \times_{\widetilde{\mathcal{C}}} \widetilde{\mathcal{C}}_{/\widetilde{x}}$
and $(\mathcal{C}_{/\widetilde{x}})^{\vartriangleright} :=
(\mathcal{C})^{\vartriangleright} \times_{\widetilde{\mathcal{C}}}
\widetilde{\mathcal{C}}_{/\widetilde{x}}$, where the final object in 
$(\mathcal{C})^{\vartriangleright}$ is mapped to $\widetilde{x}$.

\bdef
\label{def: LKE along inclusion}
A left Kan extension of $F$ along $i$
is a functor, $\widetilde{F}: \widetilde{\mathcal{C}} \to \mathcal{D}$,
making the following diagram commute,
\[\begin{tikzcd}
	& {\mathcal{C}} \\
	{\widetilde{\mathcal{C}}} & {\mathcal{D}}
	\arrow["i", from=2-1, to=1-2]
	\arrow["{\widetilde{F}}", from=1-2, to=2-2]
	\arrow["F"', from=2-1, to=2-2],
\end{tikzcd}\]
such that for every object $C \in \widetilde{\mathcal{C}}$,
the induced diagram,
\[\begin{tikzcd}
	& {(\mathcal{C}_{/\widetilde{x}})^{\vartriangleright}} \\
	{\mathcal{C}_{/\widetilde{x}}} & {\mathcal{D}}
	\arrow["i", from=2-1, to=1-2]
	\arrow[from=1-2, to=2-2]
	\arrow["F_{\widetilde{x}}"', from=2-1, to=2-2],
\end{tikzcd}\]
exhibits $\widetilde{F}(\widetilde{x})$ as the colimit
of $F_{\widetilde{x}}$.
\edefn

\bdef
\label{def: LKE}
A left Kan extension of $F$ along
an arbitrary functor $\delta: \mathcal{C}
\to \widetilde{\mathcal{C}}$ is defined to be a left
Kan extension of $F$ along the 
inclusion of $\mathcal{C}$
into the mapping cylinder 
of $F$, $\mathcal{C} \times \{0\} 
\hook (\mathcal{C} \times [1])
\coprod_{\mathcal{C} \times \{1\}}
\widetilde{\mathcal{C}} =: M_F$.
\edefn

\brem
Left Kan extensions are
unique up to a contractible space
of choices. We denote the left Kan
extension of $F$ along $\delta$ by
$\LKE_{\delta}(F)$.
\erem

\brem
Lurie shows that \cref{def: LKE} 
coincides with \cref{def: LKE along 
inclusion} in the case when
$\delta$ is a fully faithful embedding.
\erem

\brem
Note that \cref{def: LKE along inclusion}
says the left Kan extension of $F$
along $i$ is any functor
whose on an object $\widetilde{x} \in \widetilde{\mathcal{C}}$ is
the colimit of the diagram,
$\mathcal{C}_{/\widetilde{x}} \to \mathcal{D}$.
Similarly, the value of $\LKE_{\delta}(F)$ on an
object $\widetilde{x} \in \widetilde{\mathcal{C}}$ is the colimit
of the diagram $\mathcal{C} \times_{M_F} {M_F}_{/\widetilde{x}}
\to \mathcal{D}$. But the natural map $\mathcal{C}
\times_{M_F} {M_F}_{/\widetilde{x}} \to \mathcal{C}_{/\widetilde{x}}$ is
an equivalence since $\widetilde{\mathcal{C}} \simeq M_F$,
so, in fact, \[\LKE_{\delta}(F)(\widetilde{x}) = 
\varinjlim_{x \in \mathcal{C}_{/\widetilde{x}}} F(x).\]  
\erem

\subsubsection{}
Let $I := \B_{\R^+}(V^{\vee}_Z)^{op} \times (\R^+)^{op}$
be the product category, 
where $\R^+$, being a totally ordered set,
is considered as a category.
Similarly, let 
$J$ denote the full subcategory of
$\B_{\R^+}(V^{\vee}_Z)^{op} \times \U(V_Z)^{op}$
on pairs $\{(\V, \U)| \V \in \ob(\B_{\R^+}(V^{\vee}_Z)), 
\U \in \ob(\mathfrak{F}_{\U_{\V}})\}$. 
There are obvious functors, $F: I \to \B_{\R^+}(V^{\vee}_Z)^{op}$
and $G: J \to \B_{\R^+}(V^{\vee}_Z)^{op}$ given by the
respective projection functors, as well as a functor $\alpha: I \to J$
defined by $(\V, \epsilon) \mapsto (\V, \U_{\epsilon, \V})$. 
Now define a functor $Q: J \to \Mod^b_e$ 
by the assignment $(\V, \U)$ to 
$\Gamma(\U; \mathcal{F}) \in \Mod^b_e$;
in other words, define $Q$ to be the composition
of the projection $J \to \U(V_Z)$ and the
sheaf $\F$. Altogether, we have
the following commutative diagram,
\begin{equation}
\label{diagram: extending map 1}
\begin{tikzcd}
	I & J & {\Mod^b_e} \\
	& {\B_{\R^+}(V^{\vee}_Z)^{op}}
	\arrow["\alpha", from=1-1, to=1-2]
	\arrow["Q", from=1-2, to=1-3]
	\arrow["F"', from=1-1, to=2-2]
	\arrow["G", from=1-2, to=2-2].
\end{tikzcd}
\end{equation}

\subsubsection{}
The category 
$\B_{\R^+}(V^{\vee}_Z)$
has a natural Grothendieck
topology that we use to
consider ``basis" sheaves, 
$\Shv(\B_{\R^+}(V^{\vee}_Z))$.
There is a canonical 
restriction functor,
$\theta: \Shv({(V^{\vee}_{\con})}_Z) 
\to \Shv(\B_{\R^+}
(V^{\vee}_Z))$, given
by restriction along the
inclusion $\B_{\R^+}(V^{\vee}_Z)
\hook \U({V^{\vee}_{\con}}_Z)$.
 
\blem
\label{lemma: psi_m and psi_m^Sp are Kan extensions}
$\theta(\psi_m)$ and $\theta(\psi_m^{\Sp})$ are the
left Kan extensions
of $Q \circ \alpha$ along $F$
and of $Q$ along $G$, respectively.
\elem

\bproof
We begin with the first claim.
It was shown in the previous section
that $\psi_m(\V)$ is the colimit of the
following $(\R^+)^{op}$-shaped diagram, 
\[\begin{tikzcd}[row sep=0]
	{(\R^+)^{op}} & {\Mod^b_e} \\
	\rotatebox{90}{$\in$} & \rotatebox{90}{$\in$} \\
	\epsilon & {\Gamma(\U_{\epsilon, \V}; \F)}
	\arrow[maps to, from=3-1, to=3-2]
	\arrow[from=1-1, to=1-2]
\end{tikzcd}\]
This is a diagram with the same
terms as the diagram, which we denote $D_{\V}$, whose colimit computes
$\LKE_F(Q \circ \alpha)(\V)$ but it is indexed 
by a different category. Nonetheless, there
is an obvious cofinal map $r_{\V}: (\R^+)^{op} \to I_{/\V}$
given by the assignment $\epsilon \mapsto (\V, \epsilon)$.
Indeed, given any $(\V', \epsilon) \in I_{/\V}$, there
is a unique morphism $(\V', \epsilon) \to (\V, \epsilon)$
coming from the fact that $\V \xrightarrow{id} \V$ is
final in $(\B_{\R^+}(V^{\vee}_Z)^{op})_{/\V}$.
Thus, 
\[\colim(D_{\V} \circ r_{\V}) \simeq \colim(D_{\V}),\]
and we have shown that $\psi_m$ is the left Kan extension,
$\LKE_F(Q \circ \alpha)$. 

Likewise, in the previous section we
established that $\psi_m^{\Sp}(\V)$ is the colimit of
the $\mathfrak{F}_{\U_{\V}}$-shaped diagram,
\[\begin{tikzcd}[row sep=0]
	{\mathfrak{F}_{\U_{\V}}} & {\Mod^b_e} \\
	\rotatebox{90}{$\in$} & \rotatebox{90}{$\in$} \\
	\U & {\Gamma(\U; \F)}
	\arrow[maps to, from=3-1, to=3-2]
	\arrow[from=1-1, to=1-2]
\end{tikzcd}\]
On the other hand, $\LKE_G(Q)(\V)$ is the
colimit of the diagram $D'_{\V}: J_{/\V} \to \Mod^b_e$.
In this case, as well, there is an
obvious map $\mathfrak{F}_{\U_{\V}} \to
J_{/\V}$ given by $\U \mapsto (\V, \U)$,
which is cofinal for the same reason as above.
Thus, $\psi_m^{\Sp} \simeq \LKE_G(Q)$, and
the claim is proven. 
\eproof

From the commutativity of
\labelcref{diagram: extending map 1} 
we obtain a map of Kan extensions,
$\psi: \LKE_F(Q \circ \alpha) \to \LKE_G(Q)$.
Unwinding the definitions,
it is clear that the map $\psi(\V): \colim(D_{\V}) 
\to \colim(D'_{\V})$ is necessarily
unique morphism coming from
the universal property of the first
colimit. Upon identification
of $\LKE_F(Q \circ \alpha)$ and
$\LKE_G(Q)$ with $\psi_m$ and $\psi_m^{\Sp}$,
respectively, then, we see that $\psi(\V)$
is precisely the map $\psi_{\V}$ constructed
in \cref{rem: defining psi_V}.
We thus obtain 
a map of basis sheaves,
$\theta(\psi_m) \xrightarrow{\psi} 
\theta(\psi_m^{\Sp})$.

\subsubsection{Extension of basis sheaves}
The following lemma
allows us to lift $\psi$
to a map of conic sheaves
on $V^{\vee}_Z$.

\blem
The canonical restriction map
\[\theta: \Shv_{\R^+}(V^{\vee}_Z) \to \Shv(\B_{\R^+}(V^{\vee}_Z))\]
is an equivalence of
categories.
\elem

\bproof
The following argument
is based on the discussion
in \cite{HTT} right before
Warning 7.1.1.4.
Observe that the
category $\B_{\R^+}(V^{\vee}_Z)$
has all finite products since
the intersection of finitely
many convex, conic open subsets
is again convex, conic, and open.
Moreover, there exists a final
object: $V^{\vee}_Z$. By 
\cite[Lemma 002O]{Stacks}, 
$\B_{\R^+}(V^{\vee}_Z)$ 
therefore admits all finite limits.
Moreover, it is a small
$0$\-category. The proof
of \cite[Proposition
6.4.5.7]{HA} shows that
the $\infty$-topos, 
$\Shv(\B_{\R^+}(V^{\vee}_Z); 
\mathcal{S})$,
is a $0$-localic. As such, it
is determined by the locale\footnotemark
\footnotetext{\cite[Definition 6.4.2.3]{HTT}}
of subobjects of the final
object $\underline{\ast} \in
\Shv(\B_{\R^+}(V^{\vee}_Z); 
\mathcal{S}),$ which in
turn is equivalent to
the locale given by
the partially ordered set
$\B_{\R^+}(V^{\vee}_Z)$ itself.
For the same reasons, the topos 
$\Shv_{\R^+}(V^{\vee}_Z; \mathcal{S})$
is determined by the partially
ordered set $\U_{\R^+}(V^{\vee}_Z)$.

Now observe that $\B_{\R^+}(V^{\vee}_Z)$
being a basis means precisely
that the inclusion of partially
ordered sets, $\B_{\R^+}(V^{\vee}_Z)
\hook \U_{\R^+}(V^{\vee}_Z)$ induces
an isomorphism of locales.
Thus, the restriction map
$\theta_{\mathcal{S}}:
\Shv((V^{\vee}_{\con})_Z; \mathcal{S})
\to \Shv(\B_{\R^+}(V^{\vee}_Z); 
\mathcal{S})$ of $\infty$-topoi
is an isomorphism.

Now recall \cite[Proposition
1.3.1.7]{SAG}, which states that,
for a small $\infty$-category
$\mathcal{T}$ equipped with
a Grothendieck topology and $\Cee$
an arbitrary $\infty$-category
admitting small limits,
\[\Shv_{\Cee}(\Shv(\mathcal{T}; 
\mathcal{S})) \xrightarrow{\simeq} \Shv(\mathcal{T}; \Cee),\]
where $\Shv_{\Cee}(-)$ denotes
the category of small limit-preserving
functors $(-)^{\op} \to \Cee$.
Then, in light of the equivalence
$\theta_{\mathcal{S}}$, the restriction map,
$\Shv_{\R^+}(V^{\vee}_Z; \Mod^b_e)
\xrightarrow{\theta} \Shv(\B_{\R^+}(V^{\vee}_Z); 
\Mod^b_e)$, 
is an equivalence. This suffices
to prove the lemma, as $\Shv(-)$
is a full subcategory of $\Shv(-; \Mod^b_e)$.
\eproof

Thus, $\psi$ lifts
to a map of conic
sheaves (which we
also call $\psi$ by abuse of
notation),
\begin{equation}
\label{equation: comparison map}
\psi_m(\pr_1^*\F)|_{V^{\vee}_Z} \xrightarrow{\psi}
	\psi_m(\pr_1^*\Sp_{0_Z/V_Z}\F)|_{V^{\vee}_Z},
\end{equation}
such that for
$\V \in \B_{\R^+}(V^{\vee}_Z)$,
$\psi(\V) = \psi_{\V}$, as desired.


\subsection{Proof of \cref{lemma: phi factors through sp}}

We proceed to show that the comparison
morphism, \labelcref{equation: 
comparison map}, is an isomorphism
by computing the induced map of
stalks. The map of stalks
is shown to be an 
isomorphism at each point 
by finding a common 
cofinal subsequence in each
of the colimits computing the stalks
of each the sheaves, 
$\psi_m(\pr_1^*\F)|_{V^{\vee}_Z}$ and
$\psi_m(\pr_1^*\Sp_{0_Z/V_Z})|_{V^{\vee}_Z}$.

\subsubsection{Proving the induced map on stalks is an isomorphism}
Choose a covector,
$(z, \lambda) \in V^{\vee}_Z$.
Clearly, the comparison map
\labelcref{equation: comparison map}
induces the following map of
stalks at $(z, \lambda)$:
\begin{equation}
\label{equation: map of stalks}
\varinjlim_{(z, \lambda) \in \V} \left(\varinjlim_{\epsilon>0} 
	\Gamma(\U_{\epsilon, \V}; \F)\right) \to \varinjlim_{(z, \lambda) \in \V} 
		\left(\varinjlim_{\U \in \mathfrak{F}_{\U_{\V}}} \Gamma(\U; \F)\right),
\end{equation}
where $\V$ ranges over
the family of convex, conic
open subsets of $V^{\vee}_Z$
containing $(z, \lambda)$.

Let $\{\W_m\}_{m \in \mathbb{N}}$ be 
a countable neighborhood
basis of $z \in Z$. If $\lambda
\neq 0$, let $\gamma_n$,
$n \in \mathbb{N}$,
be a countable neighborhood
basis of $\lambda$ in 
the conic topology on $V^{\vee}$
such that $\gamma_{n'} \subsetneq
\gamma_{n}$ for all $n'>n$.
If $\lambda =0$, let $\gamma_n
= V^{\vee}$ for all $n$.

Clearly, the subcollection $\{\V_{m, n} :=
\W_m \times \gamma_n\}$
is cofinal for both the outer left- and
right-hand colimits above. Furthermore,
\[\varinjlim_{k,m,n} \Gamma(\U_{k, \V_{m,n}}; \F) 
	\xrightarrow{\simeq} \varinjlim_{(z, \lambda) \in \V} \left(\varinjlim_{\epsilon>0} 
		\Gamma(\U_{\epsilon, \V}; \F)\right),\]
where $\{\U_{k, \V_{m,n}}\}_{k \in \mathbb{N}}$ is the
subcollection of $\{\U_{\epsilon, \V_{m,n}}\}$
where $\epsilon = \frac{1}{k}$.

We now show that the subsequence
$\{\U_{k, \V_{m,n}}\} \subset
\{\mathfrak{F}_{\U_{\V_{m,n}}}\}$
is cofinal for the right-hand side
as well by demonstrating that
any elements of the latter
contains an element of the former.
In order to ease the exposition,
we assume without loss of
generality that $Z = \ast$ (so that
we can drop the index $m$), and
denote $\U_{\V_n}$ and
$\mathfrak{F}_{\U_{\V_n}}$
by $\U_n$ and $\mathfrak{F}_n$,
respectively.

First suppose that $\lambda \neq 0$.
Fix $n$, and suppose 
that $\U \in \mathfrak{F}_n$. Recall that this
means $C_{0}(V \setminus \U)
\cap \U_n = \varnothing$, where
$\U_n = \pr_1(V \times
\gamma_n \cap m^{-1}(\Aone_{>0}))$. 

\begin{claim} 
There exists an $\epsilon_n>0$
such that $B_{\epsilon_n}(0) \cap \U_{\ell}
\subset \U$ for all $\ell>n$.
\end{claim}

\bproof
Suppose that no such $\epsilon_n$
existed. Then for each $\epsilon>0$,
$(B_{\epsilon}(0) \cap U_m) \cap
\partial \U \neq \varnothing$. To see
this, note that $(B_{\epsilon}(0) 
\cap \U_{\ell}) \cap \U \neq \varnothing$,
for otherwise,
\begin{align*}
C_0(V \setminus \U) \cap \U_n 	&\supset C_0(V \setminus \U) \cap \U_{\ell} \\
						&= C_0(V \setminus \U) \cap C_0(B{\epsilon}(0) \cap \U_{\ell}) \\
						&= C_0(B_{\epsilon}(0) \setminus (B_{\epsilon}(0) \cap \U)) \cap C_0(B_{\epsilon}(0) \cap \U_{\ell}) \\
						&= T_0V \cap \overline{p_{>0}^{-1}(B_{\epsilon}(0) \setminus (B_{\epsilon}(0) \cap \U))} \cap \overline{p_{>0}^{-1}(B_{\epsilon}(0) \cap \U_{\ell})} \\
						&\supset T_0V \cap \overline{p_{>0}^{-1}(B_{\epsilon}(0) \setminus (B_{\epsilon}(0) \cap \U)) \cap p_{>0}^{-1}(B_{\epsilon}(0) \cap \U_{\ell})} \\
						&= T_0V \cap \overline{p_{>0}^{-1}(B_{\epsilon}(0) \setminus (B_{\epsilon}(0) \cap \U) \cap B_{\epsilon}(0) \cap \U_{\ell})} \\
						&= T_0V \cap \overline{p_{>0}^{-1}(B_{\epsilon}(0) \cap \U_{\ell})} \\
						&= \U_{\ell} \\
						&\neq \varnothing,
\end{align*}
contradicting our assumption
that $\U \in \mathfrak{F}_n$.

Now suppose that  $(B_{\epsilon}(0)
\cap \U_{\ell}) \cap \partial \U = \varnothing$.
Observe that $\partial((B_{\epsilon}(0) 
\cap \U_{\ell}) \cap \U) \neq \varnothing$ because
$(B_{\epsilon}(0) \cap \U_{\ell}) \cap \U$ is a
nonempty, open proper subset of $V
\simeq \R^{\dim V}$. Thus, given a
point $x \in \partial((B_{\epsilon}(0) \cap \U_{\ell})
\cap \U)$, there is a sequence $x_i \in
(B_{\epsilon}(0) \cap \U_{\ell}) \cap \U$
such that $x_i \to x$. If $x \notin \partial(B_{
\epsilon}(0) \cap \U_{\ell}) \cap \partial \U$,
then $x \in (B_{\epsilon}(0) \cap \U_{\ell}) \cap \U$
necessarily since $x$ belongs to the
closures of both $B_{\epsilon}(0) \cap \U_{\ell}$
and $\U$---a contradiction. Since $B_{\epsilon}(0) \cap \U_{\ell}
\neq \U$ by assumption, there exists a point
in $(B_{\epsilon}(0) \cap \U_{\ell}) \cap \partial \U$;
simply choose $x \in (B_{\epsilon}(0) \cap \U_{\ell}) 
\cap \partial((B_{\epsilon}(0) 
\cap \U_{\ell}) \cap \U)$.

Now choose a point $x_{\epsilon} \in 
(B_{\epsilon}(0) \cap \U_{\ell}) 
\cap \partial \U$. Doing
so for each value $\epsilon = \frac{1}{i}$
yields a sequence in $V \setminus
\U$, $x_i \to 0$. Normalize the terms
of this sequence as $x'_i := \frac{x_i}{|x_i|}$.
This is a bounded sequence in
$\U_{\ell} \subset V$, and therefore contains
a subsequence $\{x'_{i_j}\}_j$
converging to $v \in \overline{\U_{\ell}}$.
By \cite[Proposition 4.1.2(ii)]{KS90}, $v
\in C_0(V \setminus \U)$, so
$C_0(V \setminus \U) \cap \U_n \supset
C_0(V \setminus \U) \cap \overline{\U_{\ell}}
\neq \varnothing$, a contradiction.
\eproof

Choose $\frac{1}{k} < \epsilon_n$.
Then $\U_{n,k} \subset \U$, which
proves the equivalence of stalks
at $\lambda \neq 0$.

Finally, suppose that $\lambda = 0$. 
A neighborhood basis for $0 \in V^{\vee}$
in the conic topology was given
by the single, nonempty conic
open containing: $\gamma_n = V^{\vee}$.
Consequently, equivalence of stalks reduces
to the claim that,
\[\varinjlim_{\epsilon>0} \Gamma(B_{\epsilon}(0); \F) 
	\xrightarrow{\simeq} \varinjlim_{\U \in \mathfrak{F}_V} \Gamma(\U; \F).\]
But this is identical to
the claim that $\F_0 \xrightarrow{\simeq}
(\Sp_{0/V}\F)_0$ which was
shown earlier, in
\cref{lemma: restriction to base}
(alternatively, in the proof
of \cite[Theorem 4.2.1(iv)]{KS90}
provided there).
Thus, the comparison map,
\labelcref{equation: comparison map},
is an equivalence.

\subsubsection{Finishing the proof of \cref{lemma: phi factors through sp}}
Observe that, by smooth base change,
the naive $*$-restriction of $\pr_1^*\F$
and $\pr_1^*\Sp_{0_Z/V_Z}$ to the
subspace $V^{\vee}_Z$ coincide:
\begin{align*}
(\pr_1^*\F)|_{V^{\vee}_Z} 	&\simeq \F|_{0_Z} \boxtimes \const_{V^{\vee}_Z} \\
						&\simeq (\Sp_{0_Z/V_Z}\F)|_{0_Z} \boxtimes \const_{V^{\vee}_Z} \\
						&\simeq (\pr_1^*\Sp_{0_Z/V_Z}\F)|_{V^{\vee}_Z}.
\end{align*}
As such, we obtain a
commutative square
which extends to the 
map of fiber sequences,
\[\begin{tikzcd}
	{\varphi_m(\pr_1^*\F)} & {(\pr_1^*\F)|_{V^{\vee}_Z}} & {\psi_m(\pr_1^*\F)|_{V^{\vee}_Z}} \\
	{\varphi_m(\pr_1^*\Sp_{0_Z/V_Z}\F)} & {(\pr_1^*\Sp_{0_Z/V_Z}\F)|_{V^{\vee}_Z}} & {\psi_m(\pr_1^*\Sp_{0_Z/V_Z}\F)|_{V^{\vee}_Z}}
	\arrow["\simeq", from=1-2, to=2-2]
	\arrow["\simeq", from=1-3, to=2-3]
	\arrow[from=2-2, to=2-3]
	\arrow[from=2-1, to=2-2]
	\arrow[from=1-1, to=1-2]
	\arrow["\varphi", from=1-1, to=2-1]
	\arrow[from=1-2, to=1-3]
\end{tikzcd}\]
via the canonical fiber 
sequences \labelcref{fiber sequence: 
nearby/vanishing cycles}.
We deduce that the induced
map $\varphi$ is also an equivalence
by the three lemma for stable
$\infty$-categories.
This completes the proof of
\cref{lemma: phi 
factors through sp}.


\section{The global equivalence}
\label{sec: the global equivalence}

In this section we prove the main result of this
paper, which is a globalized version of
\cref{theorem: local statement}, 
proven in the previous section.

\bthm
\label{theorem: global statement}
Given any closed immersion $\bm{Z} \hook X$ 
of a quasi-smooth derived scheme $\Z$ into 
a smooth scheme $X$, there exists an isomorphism,
\[ \varphi_{\T^*[-1]\Z} \simeq \mu_{\Z/X}[\vircodim(\Z,X)](\const_{X}[\dim X]),\]
of perverse sheaves on $T^*[-1]\Z \subset N^{\vee}_{\Z/X}$.
\ethm
The left-hand side of this equivalence
is the perverse sheaf associated to a
$-1$-shifted symplectic scheme,
as defined in \cite{BBDJSS15}.
We briefly review its construction,
and state some relevant results.


\subsection{d-critical loci and the canonical perverse sheaf}
\label{ssec: d-critical loci and the canonical perverse sheaf}
\subsubsection{}
We assume familiarity with
the basic concepts and definitions
of shifted symplectic geometry,
as appears in e.g. \cite{PTVV13}.

\subsubsection{d-critical loci}
An algebraic d-critical locus,
defined in \cite{J15}, is a classical model
for the structure of a $-1$-shifted
symplectic derived scheme. 
Its definition consists 
of a classical scheme
equipped with extra data. We state
the precise definition after introducing
some notation.

\begin{notn}
The exposition we present
here is largely derived
from \cite{K21a}, which
in turn is a summary of
the results and ideas of
\cite{J15}.

To any $\C$-scheme
$X$, locally of finite type,
Joyce associates a
$\C$-sheaf $\mathcal{S}_X$,
whose sections, loosely speaking,
parametrize the various
ways of writing $X$ as the
critical locus of a regular function
on a smooth scheme.

The sheaf $\mathcal{S}_X$
enjoys the following (more-or-less
defining) property (\cite{J15}):
for an open $R \subset X$
and any closed embedding
$i: R \hook U$ into a smooth scheme
$U$, there is an exact sequence
of sheaves,
\[0 \to {\mathcal{S}_X}|_R \xrightarrow{\iota_{R,U}} 
	\frac{i^{-1}\O_U}{I^2_{R, U}} \xrightarrow{d} 
		\frac{i^{-1}\Omega_U}{I_{R, U} \cdot i^{-1}\Omega_U} \to 0,\]
which identifies $\mathcal{S}_X$ as
the kernel of the de Rham differential $d$;
here $I_{R,U} \subset i^{-1}(\O_U)$
is the ideal sheaf of $i^{-1}\O_U$
corresponding to the closed embedding
$i: R \hook U$.

The compositions,
\[{\mathcal{S}_X}|_R \to 
	\frac{i^{-1}\O_U}{I^2_{R,U}} 
		\to \frac{i^{-1}\O_U}{I_{U,R}} \simeq \O_R,\]
for various $R \subset X$
are compatible and glue to
a map $\mathcal{S}_X 
\xrightarrow{\beta_X} \O_X$.
Joyce defines a subsheaf $\mathcal{S}_X^0
\subset \mathcal{S}_X$ as the
kernel of the composition,
\[\mathcal{S}_X \xrightarrow{\beta_X} \O_X \to \O_{X^{red}},\]
and shows that there is
decomposition $\mathcal{S}_X
= \underline{\C}_X \oplus \mathcal{S}_X^0$.
\end{notn}

\bdef{\cite[Definition 2.5]{J15}}
An algebraic d-critical locus
is a pair $(X, s)$, where $X$ is a 
$\C$-scheme (locally of finite type),
and $s \in \Gamma(X; \mathcal{S}_X^0)$
satisfying the condition that $X$ may
be covered by Zariski open sets
$R \subset X$ with a closed embedding
$i: R \hook U$ into a smooth scheme $U$, 
and a regular function
$f$ on $U$ such that $i(R) = \crit(f)$ and 
$f+I^2_{R,U}={s|}_{R}$. 
The tuple $\mathscr{R} = 
(R, U, f, i)$ as above is called a 
critical chart for $(X, s)$. 
\edefn

\brem
There is an appropriate notion
of embedding of critical charts\footnotemark
\footnotetext{\cite[Definition 2.18]{J15}}
which allows us to compare two
different, overlapping critical charts.
Indeed, by \cite[Theorem 2.20]{J15}, 
around each point in the overlap of
critical charts $\mathscr{R}_1$ and 
$\mathscr{R}_2$, we may find a third
critical chart into which $\mathscr{R}_1$
and $\mathscr{R}_2$ commonly embed.
\erem

\brem
Joyce also defines the notion
of a complex analytic 
d-critical locus in a very similar fashion.
For our purposes, the notion of
algebraic d-critical locus suffices,
but the results of this section are
equally valid should we consider
d-critical structures
on the analytifications of the
schemes in consideration, rather
than algebraic d-critical structures
on the schemes themselves.
\erem

\subsubsection{Oriented d-critical loci}
Joyce proves the existence of
a canonical line bundle on any
algebraic d-critical locus.

\bthm[{\cite[Theorem 2.28]{J15}}]
Let $(X,s)$ be an algebraic
d-critical locus. There is a natural
line bundle $K_{X,s}$ on $X^{red}$
called the canonical line bundle
characterized by the following properties:
\begin{enumerate}[(i)]
\item If $\mathscr{R} = (R,U,f,i)$ is
a critical chart of $(X,s)$, we have
an isomorphism,
\[\iota_{\mathscr{R}}: {K_{X,s}|}_{R^{red}} 
	\xrightarrow{\simeq} {(i^*K_U)^{\otimes^2}|}_{R^{red}},\]
where $K_U$ is the usual canonical
bundle on $U$.
\item If $\Phi: \mathscr{R}_1 \hook \mathscr{R}_2$ is
an embedding of critical charts, we have
\[{\iota_{\mathscr{R}_2}|}_{R_1^{red}} = 
	J_{\Phi} \circ \iota_{\mathscr{R}_1},\]
where $J_{\Phi}$ is a natural isomorphism
$(i_1^*K_{U_1})^{\otimes^2}|_{R_1^{red}} \simeq 
(i_2^*K_{U_2})^{\otimes^2}|_{R_1^{red}}$\footnotemark.
\end{enumerate}
\ethm
\footnotetext{\cite[Definition 2.26]{J15}}

An orientation on a d-critical locus
is then defined in terms of this 
canonical bundle.

\bdef
An orientation
on a d-critical locus $(X,s)$
is a choice of square root $K_{x,s}^{1/2}$ for
the canonical bundle on $X^{red}$,
\[o: (K_{X,s}^{1/2})^{\otimes^2} \xrightarrow{\simeq} K_{X,s}.\]
\edefn

A d-critical locus $(X,s)$ equipped with
an orientation $o$ is called an oriented
d-critical locus, which we denote
by the tuple $(X,s,o)$.

\subsubsection{The canonical perverse sheaf on an oriented d-critical locus}

To any oriented d-critical locus $(X,s,o)$,
Joyce et al.\ introduce a canonical
perverse sheaf, $\varphi_{X,s,o}$, on $X$ 
which is a kind of globalized vanishing cycles (see
\cref{rem: vanishing cycles name convention}) sheaf; 
namely, if $(X,s)$ is locally modeled on
the critical locus $\crit(f: U \to \Aone)$,
$\varphi_{X,s,o}$ is locally modeled
on the vanishing cycles along $f$, up
to a twist by a $\mathbb{Z}/2 \mathbb{Z}$-bundle 
determined by the orientation
$o$. The precise statement follows.

\bthm[{\cite[Theorem 6.9(i)]{BBDJSS15}}]
Let $(X,s, o)$ be an oriented d-critical
locus, with orientation given by square
root $K_{X,s}^{1/2}$. Then there exists
a perverse sheaf $\varphi_{X,s,o} \in \Perv(X)$
which is natural up to canonical isomorphism,
such that if $\mathscr{R} = (R,U,f,i)$ is a critical chart
on $(X,s)$, there is a natural isomorphism,
\[\omega_{\mathscr{R}}: {\varphi_{X,s,o}|}_{R}
	\to i^*(\varphi_f) \otimes_{\mathbb{Z}/2\mathbb{Z}} Q_{\mathscr{R}},\]
where $Q_{\mathscr{R}} \to R$ is the principal
$\mathbb{Z}/2\mathbb{Z}$-bundle parametrizing
local isomorphisms $\alpha: K_{X,s}^{1/2} \to
(i^*K_U)|_{R^{red}}$ such that $\alpha \otimes
\alpha = \iota_{\mathscr{R}} \circ o$.
\ethm

\subsubsection{The canonical perverse sheaf on $T^*[-1]\Z$}

We remarked above that d-critical loci
are classical models for $-1$-shifted
symplectic schemes. This remark is 
made precise by the following theorem.

\bthm[{\cite[Theorem 6.6]{BBJ19}}]
Let $(\X, \omega_{\X})$ be a $-1$-shifted
symplectic derived scheme. Then the underlying
classical scheme $X$ carries a canonical
d-critical structure $s_X$ characterized by the
following property: if $(\bm{R} = \dcrit(f), 
\omega_{\bm{R}})$ is a $-1$-shifted 
symplectic scheme, where $f$ is a regular
function on a smooth scheme $U$, 
$\omega_{\bm{R}}$ is the usual 
$-1$-shifted symplectic form
on a derived critical locus\footnotemark,
and an open inclusion $\bm{\iota}: \bm{R} \hook \X$
such that $\bm{\iota}^*\omega_{\X} ~ 
\omega_{\bm{R}}$,
the tuple $(R,U,f,i)$ gives a critical chart
for $(X, s_X)$, where $i: R \hook U$ is
the natural closed embedding.
\ethm
\footnotetext{Constructed
in e.g. \cite[Example 2.5]{K21a} or
\cite[Example 5.8 and 5.15]{BBJ19}}

In particular, the $-1$-shifted cotangent
bundle, $\T^*[-1]\Z$, of any 
derived scheme $\Z$ carries
a canonical $-1$-symplectic structure,
and therefore the underlying classical 
scheme $T^*[-1]\Z$ carries a canonical 
d-critical structure, the section of 
which is denoted by $s_{\T^*[-1]\Z}$. 
When $\Z$ is additionally 
quasi-smooth, the d-critical
locus $(T^*[-1]\Z, s_{\T^*[-1]\Z})$
carries a canonical orientation\footnotemark\ 
denoted by $o_{\T^*[-1]\Z}$.   
\footnotetext{See
\cite[Example 2.15]{K21a} or 
\cite[Lemma 3.3.3]{T21}.}

\bdef
Given a quasi-smooth derived scheme
$\Z$, we denote the canonical perverse sheaf, \\
$\varphi_{(T^*[-1]\Z, s_{\T^*[-1]\Z}, 
o_{\T^*[-1]\Z})}$, by $\varphi_{\T^*[-1]\Z}$.
\edefn

In general, $-1$-shifted symplectic
schemes are Zariski locally modeled on
a derived critical locus\footnotemark.
	\footnotetext{\cite[Theorem 5.18]{BBJ19}}
On the other hand, quasi-smooth 
derived schemes are locally 
modeled on the derived zero
locus of a section $s \in \Gamma(U; 
\mathscr{E})$ of a vector bundle $\mathscr{E}$
on a smooth affine scheme $U$\footnotemark.
	\footnotetext{\cite[Proposition 2.3]{K21a}}
When the quasi-smooth scheme $\Z$
is locally modeled on the zero locus,
$\Z(s)$, $\T^*[-1]\Z$
is locally modeled on the derived
critical locus $\dcrit(\widetilde{s})$,
where $\widetilde{s}$ is the associated
regular function on $E^{\vee} :=
\Spec_U(\Sym^{\bullet}_{\O_U}
(\mathscr{E}^{\vee}))$. The following
lemma from \cite{K21a}
establishes the compatibility between
the $-1$-shifted symplectic structures
on $\T^*[-1]\Z(s)$ coming from its
structure as a derived critical locus
and from its being a shifted
cotangent bundle.

\blem[{\cite[Lemma 2.7]{K21a}}]
\label{lemma: Kinjo 2.7}
Let $U = \Spec A$ be a smooth 
affine scheme admitting a global
\'etale coordinate, $\mathscr{E}$ 
be a trivial vector bundle on $U$, and
$s \in \Gamma(U; \mathscr{E})$ be 
a section. Denote by $\widetilde{s}$
the regular function on $E^{\vee}$ 
corresponding to $s$. Then we have 
an equivalence of $-1$-shifted symplectic 
schemes,
\[(\dcrit(\widetilde{s}), \omega_{\dcrit(\widetilde{s})})
	\simeq (\T^*[-1]\Z(s), \omega_{\T^*[-1]\Z(s)}).\]
\elem

\subsubsection{Triviality of $Q_{\mathscr{R}}$}
Later in this section, we expend a
good deal of effort proving that the
equivalences of
\cref{theorem: local statement} 
that are obtained for each member of a
particular critical chart cover of the 
oriented d-critical locus, $(T^*[-1]\Z,
s_{\T^*[-1]\Z}, o_{\T^*[-1]\Z})$,
glue. 
Fortunately, our task is made much
simpler by the apparent triviality
of $Q_{\mathscr{R}}$
for critical charts of 
$(T^*[-1]\Z, s_{\T^*[-1]\Z},
 o_{\T^*[-1]\Z})$ obtained
from expressing $\Z$ locally
as a derived zero locus $\Z(s)$.

\blem[{\cite[Lemma 2.19]{K21a}}]
Let $U,s$ be as in \cref{lemma: Kinjo 2.7}. 
Denote by $o_{\T^*[-1]\Z(s)}$ 
the canonical orientation and $\mathscr{Z} =
(\crit(\widetilde{s}), E^{\vee}, \widetilde{s}, i)$
the critical chart induced by the equivalence
of \cref{lemma: Kinjo 2.7}.
Then $Q^{o_{\T^*[-1]\Z}}_{\mathscr{Z}}$
is a trivial $\mathbb{Z}/2\mathbb{Z}$-bundle.
\elem

As such, $\omega_{\mathscr{Z}}: 
{\varphi_{\T^*[-1]\Z(s)}|}_{\crit(\widetilde{s})}
\xrightarrow{\simeq} i^*(\varphi_{\widetilde{s}})$.


\subsection{Local models of quasi-smooth schemes}

In order to prove
\cref{theorem: global statement}, 
we will some particular results
on the local structure of 
quasi-smooth closed immersions,
which we now establish.

\subsubsection{} 
The following construction 
is a very slight variant
of one that appears in the proof
of \cite[Theorem 8.4.3.18]{HA}, 
and virtually identical---apart 
from presentation---to 
the one that appears in 
\cite[Example 2.8]{BBJ19}
of what the authors call ``a standard form
cdga." We use the conventions
and notation of \cite{HA}.

\begin{construction} 
\label{construction: standard form}
Let $B$ be a connective $\mathbb{E}_\infty$-algebra locally
of finite presentation over a connective $\mathbb{E}_\infty$-ring $A$ 
such that $\pi_0(B)$
is finitely generated as a $\pi_0(A)$-algebra, 
and suppose that $\mathbb{L}_{B/A}$ is perfect
of Tor-amplitude $[0,1]$. Assume
additionally that $\pi_0 A$ is Noetherian, 
and choose generators $\{x_i\}_{i=1}^n$
for $\pi_0 B$ over $\pi_0 A$.

Choose $B(0)$ a smooth $A$-algebra 
(i.e. $\L_{B(0)/A}$ is projective) and
$f: B(0) \to B$ for which
the induced map
$\pi_0 f: \pi_0 B(0) \to \pi_0 B$
is surjective\footnotemark.

\footnotetext{Such a choice might be the
free $A$-algebra generated on $n$ 
elements with $f$ given by the choose of generators
$\{x_i\} \subset \pi_0 B$, by which we mean the following.
The ``forgetful" functor $\Omega^\infty: \Mod_A \to \Spc$ admits 
a left adjoint denoted $\Sigma^\infty: \Spc \to \Mod_A$, and 
$M$ is $\Sigma^\infty$ applied to the disjoint union of $n$ points.
The choice of elements $\{x_i\}_{i=1}^n$ in the presentation of $\pi_0(B)$
lifts to a map of $A$-modules $f: M \to B$ as follows.
By adjunction, $\pi_0 \Map(M, B) = \pi_0 \Map(\coprod \ast, \Omega^\infty(B))
= \pi_0 (\Map(*, \Omega^\infty(B))^n) =  \pi_0 \Omega^\infty(B)^n) =
(\pi_0 \Omega^\infty(B))^n = (\pi_0 B)^n$. That is, an $n$-tuple of elements
in $\pi_0 B$ determine a map $M \to B$.
The forgetful functor $\Oblv: \CAlg_A \to \Mod_A$ has a left adjoint
$\Sym_A: \Mod_A \to \CAlg_A$.
By abuse of notation, we denote the map in $\CAlg_A$ adjoint to $f$ by
$f: \Sym_A(M) \to B$.} 

The map $f$ endows $B$ with the structure 
of a $B(0)$-module, and so we can in fact view $f$
as a map in $\CAlg_{B(0)}$.
Since $B(0)$ is smooth over $A$, $\pi_0 B(0)$ 
is a finitely generated $\pi_0 A$-algebra.
Since $\pi_0 A$ is Noetherian, 
$\pi_0 B(0)$ is Noetherian. Thus
$I := \ker \pi_0 f$ is finitely generated.
Choose generators $\{g_j\}_{j=1}^m \subset
\ker \pi_0 f$ for $I$.  

From the long exact sequence associated
to a fiber sequence, the image of 
$\pi_0 \fib(f) \to \pi_0 B(0)$ is identified with
the ideal $I$. By abuse of notation
we use $\{g_j\}$ to denote both the elements in
$\pi_0 B(0)$ and a choice of their lifts to $\pi_0 \fib(f)$.
Choose a finitely generated, projective $B(0)$-module $N$
of rank $m$ and a map
$N \to \fib(f)$ sending the generators
of $N$ to the elements $\{g_j\}$ (for example, the free 
module on $m$ elements and the map generated by the choice
of $\{g_j\} \subset \pi_0 \fib(f)$). We
thereby obtain a map $g: N \to \fib(f) \to B(0)$. Hence, $g$ classifies
a diagram
	\[\begin{tikzcd}
		N & 0 \\
		{B(0)} & B
		\arrow["g"', from=1-1, to=2-1]
		\arrow["f", from=2-1, to=2-2]
		\arrow[from=1-2, to=2-2]
		\arrow[from=1-1, to=1-2]
	\end{tikzcd}\]
in $\Mod_{B(0)}$. 
Adjoint to this is a diagram in $\CAlg_{B(0)}$,
	\[\begin{tikzcd}
		\label{diagram: pushout}
		{\Sym_{B(0)}(N)} & {B(0)} \\
		{B(0)} & B
		\arrow["g"', from=1-1, to=2-1]
		\arrow["f", from=2-1, to=2-2]
		\arrow["0", from=1-1, to=1-2]
		\arrow[from=1-2, to=2-2]
		\tag{A}
	\end{tikzcd},\]
which classifies a map $h_B: B(1) := B(0) \underset{\Sym_{B(0)}(N)}{\otimes} B(0) \to B$.
\end{construction}

\begin{notn}
To prevent possible confusion, 
we sometimes denote $N$ in the
preceding construction by $N_B$.
\end{notn}

\blem
\label{lemma: local form}
	Take $B$ as in \cref{construction: standard form}.
	There exists an affine open neighborhood $\Spec C$ of 
	each point $x \in \Spec B$, such the map $h_C: C(1) \to C$
	is an isomorphism of $\mathbb{E}_\infty$-algebras over $A$.
\elem



\bproof

First observe that the hypotheses of
\cref{construction: standard form}
are stable under localization, so the statement
 of the lemma makes sense.

Fix $x \in \Spec B$, and take $\Spec B$ as an
ansatz for such a neighborhood. The proof proceeds by 
first trying to show $h_B$ is
an isomorphism, then localizing when necessary.

Observe that $\pi_0 h_B: \pi_0 B(1) \to \pi_0 B$ is an isomorphism by construction.
By \cite[Corollary 7.4.3.2]{HA}, which relates the connectivity
of the cofiber of a map to the connectivity of
its cotangent complex, 
it therefore suffices to show that the relative cotangent
complex $\L_{B/B(1)} = \L_{h_B}$ vanishes.

Consider the fiber sequence
	\[ \L_{B/A}[-1] \to \L_{B/B(0)}[-1] \to \L_{B(0)/A} \underset{B(0)}{\otimes} B\]
where here and in what follows $\L_{B/B(0)} = \L_f$.
By assumption, $\L_{B/A}[-1]$ is perfect of Tor-amplitude $[-1,0]$. Since $B(0)$ is
the free algebra on $A$-module $M$, $\L_{B(0)/A} \simeq M \underset{A}{\otimes} B(0)$.
Hence, $\L_{B(0)/A} \underset{B(0)}{\otimes} B \simeq M \underset{A}{\otimes} B$,
which is a projective module over $B$.

By \cite[Proposition 7.2.4.23]{HA}, we deduce that $\L_{B/B(0)}[-1]$ is
perfect of Tor-amplitude $\leq 0$. 
Since the cofiber of $f: B(0) \to B$ is $1$-connective,
$\L_{B/B(0)}[-1]$ is connective
(see \cite[Corollary 7.4.3.2]{HA}).

Hence is a flat $B$-module by \cite[Remark 7.4.2.22]{HA}.
As a flat $B$-module, by \cite[Remark 7.2.4.22]{HA}
$\L_{B/B(0)}[-1]$ is a finitely generated, projective $B$-module. 
If $\L_{B/B(0)}[-1]$ is free, then our choice of $N$ in 
\cref{construction: standard form} yields an
isomorphism $\L_{B/B(0)} \simeq N[1] \underset{B(0)}{\otimes} B$.

If $\L_{B/B(0)}[-1]$ is not free, choose a localization of
$\pi_0 B$ over which $\pi_0 \L_{B/B(0)}[-1]$ is a free
module. This localization of $\pi_0 B$ is given by a
multiplicative subset $S \subset \pi_0 B$, a commutative ring,
and as such is a set of homogeneous elements
of the graded ring $\pi_*B$ which satisfies both the left
and the right Ore conditions. Denote by $C$ the
localization $B[S^{-1}]$, and observe that
it is trivially also 
locally of finite presentation over $A$.

Repeat the above discussion mutatis mutandis
with $C$ in place of $B$.
We see that $\L_{C/C(0)}[-1]$ is a finitely generated, projective
$C$-module such that $\pi_0 \L_{C/C(0)}[-1]$ is free.
Note that $\pi_0 C(0) \to \pi_0 C$ is surjective.
A careful accounting of the details
in the proof of \cite[Lemma 7.2.2.19]{HA}
reveals that $\L_{C/C(0)}[-1]$ is in the essential image
of the functor $\operatorname{hProj} C(0) \to \operatorname{hProj} C$.
Thus, there exists a projective $C(0)$-module $N$ such that
$N[1] \underset{C(0)}{\otimes} C \simeq \L_{C/C(0)}$.

Using $N$ as the choice of projective module surjecting onto
$(g_1, \ldots, g_m)$ in \cref{construction: standard form},
we obtain that the first morphism in the fiber sequence
\[\L_{B(1)/B(0)} \underset{B(1)}{\otimes} B \to
	\L_{B/B(0)} \to \L_{B/B(1)} \]
is an equivalence, so $\L_{B/B(1)}$ vanishes as desired.
\eproof

\brem
\label{remark: local form}
Note that in the proof of \cref{lemma: local form}, 
$C(0)$ may be taken to be $B(0)$
localized with respect to the preimage of $S$ under the
surjection $\pi_0 B(0) \to \pi_0 B$, which is also a
multiplicative subset of homogeneous elements.
In particular, this means we could
have required $N$ in the above proof to
be free after further localizing
$C$ and $C(0)$. This observation will be crucial
to the proof of
\cref{lemma: compare charts}.
\erem


\subsubsection{}
\label{sssec: quasi-smooth schemes are locally zero loci}
The previous lemma provided
us with a particular means
of constructing a local model
for a quasi-smooth derived scheme
over an arbitrary affine base.
A related means of describing 
quasi-smooth schemes locally is
given by a Kuranishi chart,
which expresses the scheme locally
as the derived zero locus of
a section of a vector bundle
over a smooth base.
We recall the definition
given in \cite{K21a}.

\bdef[{\cite[Definition 2.2]{K21a}}]
For a quasi-smooth derived scheme $\X$, a 
Kuranishi chart is a tuple
$(Z, U, E, s, \bm{\i})$ where $Z$ is an 
open subscheme of $X$, $U$ is a 
smooth scheme, $E$ is a vector bundle 
on $U$, $s$ is a section of $E$, and 
$\bm{\i}: \Z(s) \to \X$ is an open immersion
whose image is $Z$. A Kuranishi chart is 
said to be minimal at $p = \bm{i}(q) \in Z$ if the
differential $(ds)_q: T_qU \to E_q$ is zero. 
A Kuranishi chart is called good
if $U$ is affine and has global \'etale coordinates 
(i.e. regular functions $x_1, x_2, \ldots, x_n$
such that $d_{dR}x_1, d_{dR}x_2, \ldots, d_{dR}x_n$ 
form a basis of $\Omega^1_U$\footnotemark), 
and $E$ is a trivial vector 
bundle of a constant rank.
\footnotetext{So named because the 
``coordinate" map, $U \to \An$, 
defined by such a choice 
of functions is clearly \'etale,
and conversely any such \'etale 
map will give such coordinates.} 
\edefn

\brem 
Any smooth scheme will have 
\'etale coordinates Zariski locally
around each point (see e.g.\ \cite[Lemma 054L]{Stacks}), so
any Kuranishi chart locally gives a
good Kuranishi chart.
\erem

\brem
\label{remark: Kuranishi}
When $A$ is a classical, Noetherian ring
the condition that $B$ is locally finitely presented
over $A$ is equivalent to the condition
that $B$ is locally of finite type over $A$ in the sense of
\cite{GR17}. In this case, \cref{lemma: local form} 
furnishes a proof of the standard 
fact (see e.g.\ \cite[Proposition 2.1.10]{AG15})
that a quasi-smooth, derived
scheme $\Z$ locally of finite type over $k$ 
has a Zariski-local
presentation as the derived zero locus of a
section $s: U \to E$ where $U$ is smooth
and $E$ is a vector bundle over $U$.

Indeed, given a point
$x \in \Z$, we may let $U := \Spec C(0) = \An_k$ and
$s: \Spec C(0) \to \Spec(\Sym_{C(0)}N)$
be induced by the vertical map
in the pushout diagram, \labelcref{diagram: pushout}.
Applying the functor $\Spec$ 
to the pushout diagram, 
\labelcref{diagram: pushout}, 
we obtain the Kuranishi chart:
	\[\begin{tikzcd}
		{\Spec B} & U \\
		U & E
		\arrow["s", from=2-1, to=2-2]
		\arrow["0", from=1-2, to=2-2]
		\arrow[from=1-1, to=2-1]
		\arrow[from=1-1, to=1-2]
	\end{tikzcd}\]
In fact, by the Quillen-Suslin theorem,
$E$ is a trivial vector bundle over $U = \An_k$,
so this is a good Kuranishi chart.

Conversely, suppose that 
$\Z' := \Spec B \hook \Z$ is
an open immersion. 
Then, in a straightforward way,
a good Kuranishi chart $(Z', U, E, s, \bm{\i})$ whose
image is $\Spec \pi_0 B$ furnishes 
data of the sort used
algorithmically in 
\cref{construction: standard form}
to produce $B(1) \in \CAlg_k$. 

In the notation of that construction: 
\begin{enumerate}
\item $\pi_0 B = \O_{Z'}$;
\item $B(0) = \O_U$;
\item $N$ is the sheaf of sections of $E$;
\item $s$ furnishes the vertical map $\Sym_k(N) \to B(0)$ in diagram \labelcref{diagram: pushout};
\item $B(1) = \O_{\Z(s)}$; and
\item the open immersion $\bm{\i}$ induces the natural isomorphism $B(1) \to B$.
\end{enumerate}
\erem


\subsubsection{}
The following lemma is 
the central tool in our proof of 
\cref{theorem: global statement}. 
It will allow us to compare 
the isomorphisms obtained from
\cref{theorem: local statement} 
using local expressions of $\Z$
as the derived zero locus of different
functions on $X$. It is a variant of 
\cite[Proposition 2.3(ii)]{K21a} 
which in turn is a rephrasing of
\cite[Theorem 4.2]{BBJ19}.

\blem
\label{lemma: compare charts}
Let $\Z \to U$ be a closed immersion of a quasi-smooth derived
scheme $\Z$ into smooth $U$. 
Suppose that $U $ is the union of two opens, $U_1$ and $U_2$, such that
there exist Kuranishi charts, $(Z_i, U_i, E_i, s_i, \bm{\i}_i)$,
on $\Z$, where $E_i$ is constant rank $n$.
Denote $\Z_i = \Z \underset{U}{\times} U_i$
for $i=1,2$. 

Then around each point $x \in \Z_1 \underset{\Z}{\times} \Z_2$  
there exists a third \textit{good} Kuranishi chart
$(Z', U', E', s', \bm{\i}')$ such that $x \in \bm{\i}'(\Z(s'))$ 
and $E'$ is also constant rank $n$; 
open immersions $\eta_i: U' \to U_i$;
and isomorphisms $\tau_i: E' \to \eta_i^*E_i$ with the following properties:
	\begin{enumerate}[(i)]
		\item \label{property 1: compare charts}
		$\tau_i(s') = \eta_i^*(s_i)$
		\item \label{property 2: compare charts}
		The composition $\Z(s') \to \Z(s_i) \to U$ is equivalent to $\bm{\i}'$ where the
		first map is induced by $\eta_i$ and $\tau_i$, and the second map is $\bm{\i}_i$.
	\end{enumerate}
\elem

\bproof
$\Z_i \to \Z$ are clearly affine open immersions of
derived schemes, so $\Z_1 \underset{\Z}{\times} \Z_2$ 
is an open neighborhood
of the chosen point $x$ trivially contained
in $\Z_1 \underset{\Z}{\times} \Z_2$. Choose an affine open
subscheme $\Spec C \subset \Z_1 \underset{\Z}{\times} \Z_2$
containing $x$.

By \cref{remark: Kuranishi}
the data of a Kuranishi chart and
the data used in \cref{construction: standard form}
are equivalent, so the hypotheses of
the lemma ensure that $\Z(s_i) \simeq \Spec B_i(1)$
for derived rings $B_i(1)$ constructed in the manner
of \cref{construction: standard form}.
With this in mind, the preceding discussion yields the following
commutative diagram of derived schemes:
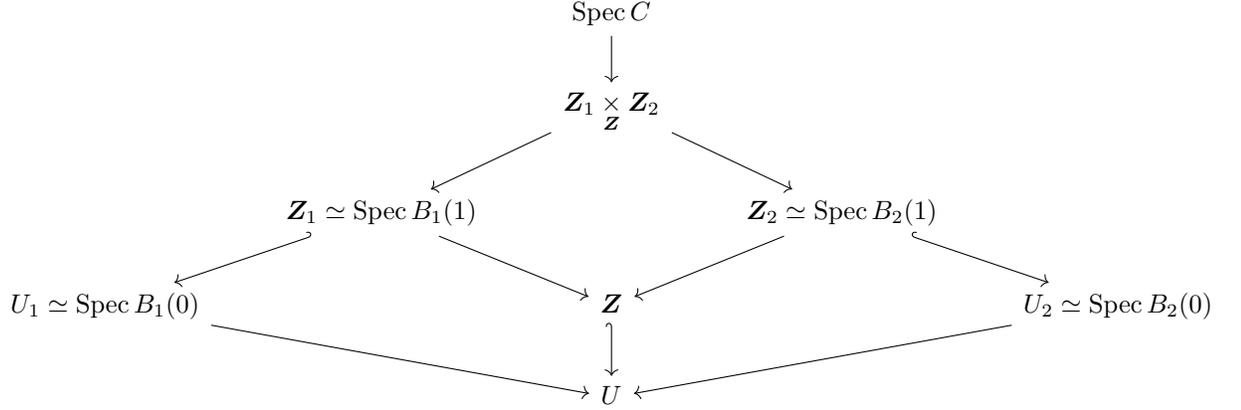
\[\begin{tikzcd}
	{\Spec C(1)} \\
	{\Z_1 \underset{\Z}{\times} \Z_2} & {\Z_2 \simeq \Spec B_2(1)} & {U_2 \simeq \Spec B_2(0)} \\
	{\Z_1 \simeq \Spec B_1(1)} & \Z \\
	{U_1 \simeq \Spec B_1(0)} && U
	\arrow[hook, from=3-2, to=4-3]
	\arrow[from=2-1, to=3-1]
	\arrow[from=2-1, to=2-2]
	\arrow[from=3-1, to=3-2]
	\arrow[from=2-2, to=3-2]
	\arrow[hook, from=3-1, to=4-1]
	\arrow[hook, from=2-2, to=2-3]
	\arrow[from=2-3, to=4-3]
	\arrow[from=4-1, to=4-3]
	\arrow["\lrcorner"{anchor=center, pos=0.125}, draw=none, from=2-1, to=3-2]
	\arrow[hook, from=1-1, to=2-1]
\end{tikzcd}\]

From the above diagram, we have a map,
$\Spec C \to U_1 \underset{U}{\times} U_2$,
which is composition of the open immersion
$\Spec C \to \Z_1 \underset{\Z}{\times} \Z_2$ followed
by the closed immersion $\Z_1 \underset{\Z}{\times} \Z_2
\hook U_1 \underset{U}{\times} U_2$.

Localizing if necessary, we may assume that the above
map factors as the composition of a closed
immersion $\Spec C \to U' = \Spec C(0)$ and an open affine
immersion $U' = \Spec C(0) \hook U_1 \underset{U}{\times} U_2$.
As an affine open subscheme of $U$, it is clear
that $U'$ is smooth.

It is clear from the 
set-up that $\Spec C$ is quasi-smooth, and
$C$ therefore satisfies the hypotheses of
\cref{construction: standard form}.
Using \cref{lemma: local form}
and localizing $C(0)$ and $C$ if necessary (see 
\cref{remark: local form}), 
we may use $C(0)$ (as notation suggests)
and the induced map $C(0) \to C$, to obtain
an object $C(1) \in \CAlg_k$ equivalent to $C$.
Thus, we have a map $\Spec C(1) \to \Spec B_i(1)$
and a specified factorization,
\[\Spec C(1) \to \Spec C(0) \to \Spec B_i(0),\]
for $i=1,2$.

Since Kuranishi charts are not intrinsic to
a quasi-smooth scheme, we need a more
explicit model for derived rings to prove
the desired statement. This more explicit
model is furnished by the model category
of cdgas over $k$ 
with the standard model
structure, $\cdga_k$. Since $k$ contains $\Q$,
the underlying $\infty$-category of
$\cdga_k$ is $\CAlg_k$.

Since $C(0)$ and $N_C$
are a discrete $k$-algebra
and discrete module over a discrete
ring, respectively, they are objects of ordinary
commutative algebra. View $C(0)$
and $\Sym_{C(0)}N_C
= \Sym_{C(0)}^\bullet N_C$ as
cdgas concentrated in degree $0$.
Then the pushout
	\[\begin{tikzcd}
		{\Sym^\bullet_{C(0)}N} & {C(0)} \\
		{C(0)} & {C(0) \underset{\Sym^\bullet_{C(0)N}}{\otimes} C(0)}
		\arrow[from=1-2, to=2-2]
		\arrow["f", from=2-1, to=2-2]
		\arrow["g"', from=1-1, to=2-1]
		\arrow["0", from=1-1, to=1-2]
	\end{tikzcd}\]
provides a cdga model for $C(1)$. We similarly
obtain cdga models for $B_i(1)$.

Assuming without loss of generality
that $N_{B_i}$ (see \cref{remark: local form})
is free, the model for $B_i(1)$
is also free over $B_i(0)$
as a graded commutative $k$-algebra,
generated by elements in degree $-1$.
As such, the structure morphism $f: B_i(0) \to B_i(1)$
is a cofibration in the model category
structure, so $B_i(1)$ is a cofibrant
object in $\cdga_{B_i(0)}$. 
Meanwhile,
the specified map $B_i(0) \to C(0)$
allows us to view the map $B_i(1) \to C(1)
\in \CAlg_k$ as a map in $\CAlg_{B_i(0)}$.
All objects in $\cdga_{B_i(0)}$ are fibrant, so 
\[\pi_0 \Hom_{\CAlg_{B_i(0)}}(B_i(1), C(1))
	\simeq \Hom_{\cdga_{B_i(0)}}(B_i(1), C(1)),\]
by the general yoga of model categories
and $\infty$-categorical localization.
The upshot is that 
the map of derived $k$-algebras,
$B_i(1) \to C(1)$, is induced
by a strict map of cdgas, for $i=1,2$.
In particular, since
$(B_i(1))^0 = B_i(0)$ and $(C(1))^0 = C(0)$,
we obtain a commutative diagram
(not necessarily Cartesian),
	\[\begin{tikzcd}
		\label{diagram: compare charts}
		{B_i(1)} & {C(1)} \\
		{B_i(0)} & {C(0)}
		\arrow[from=2-2, to=1-2]
		\arrow[from=2-1, to=2-2]
		\arrow[from=2-1, to=1-1]
		\arrow[from=1-1, to=1-2]
	\end{tikzcd}\]
in the category $\cdga_k$, where
the vertical arrows are the
inclusions of the degree $0$ parts of
the cochain complexes.
Passing to $\CAlg_k$, the underlying
$\infty$-category of $\cdga_k$,
this obtains a homotopy commutative
diagram (i.e. the only
kind of commutative diagram
in $\CAlg_k$).

There is a natural base change morphism,
\[ \L_{B_i(1)/B_i(0)} \underset{B_i(0)}{\otimes} C(0)
	\to \L_{C(1)/C(0)},\]
arising from the commutativity of this diagram,
where we have viewed $\L_{B_i(1)/B_i(0)}$
as an object in $\Mod_{B_i(0)}$ and
$\L_{C(1)/C(0)}$ as an object in
$\Mod_{C(0)}$.
Recall from the proof of
\cref{lemma: local form} that 
$\L_{B_i(1)/B_(0)}$ and $\L_{C(1)/C(0)}$
are projective modules over $B_i(0)$ 
and $C(0)$, respectively, concentrated
in homotopy degree $1$. We labeled them
$N_{B_i(1)}[1]$ and $N_{C(1)}[1]$, respectively.
As such, the base change morphism
yields an isomorphism,
\[N_{B_i(1)} \underset{B_i(0)}{\otimes} C(0) 
	\simeq N_{C(1)}.\]

With this last isomorphism, we
have unwittingly proven
\cref{lemma: compare charts}.
Explicitly,
\begin{itemize}

\item $(Z', U', E', s', \bm{\i}') :=
(\Spec \pi_0 C(1), \Spec C(0), 
\Spec_{U'}(\Sym^\bullet_{C(0)} N_{C(1)}),
\Spec(C(0) \to C(1)), 
\Spec C(1) \to \Z_1 \underset{\Z}{\times} \Z_2 \to \Z)$;

\item $\eta_i := \Spec(B_i(0) \to C(0))$ for $i=1,2$;

\item $\tau_i$ is induced by the isomorphism
$N_{B_i(1)} \underset{B_i(0)}{\otimes} C(0) 
\simeq N_{C(1)}$, which are the dual sheaf of 
sections of $\eta_i^*E_i$ and $E'$, respectively.

\end{itemize}

The desired properties and compatibilities
of these objects and maps follow clearly from
their construction above. Furthermore, by
shrinking $U'$ around $x$ we may take 
$(Z', U', E', s', \bm{\i}')$ to be good.
\eproof


\subsection{Gluing local statements}

\subsubsection{}
Continue to suppose that $\Z$ 
is quasi-smooth and 
$i: \Z \hook X$ is a 
closed immersion
into a smooth scheme $X$. 
By \cref{lemma: quasi-smooth}, 
$i$ is quasi-smooth.
Choose a Zariski open cover 
of $X$, $\{U_\a \hook X\}_{\a \in \mathcal{I}}$, 
such that on each open, 
$U_\a \hook X$, $i$ has the form
of a derived zero locus, 
$\Z(f_\a)$, per 
\cref{lemma: local model}:
	\[\begin{tikzcd}
		{\bm{Z} \times_{U_\a} X} & {U_\a} \\
		{*} & {\mathbb{A}^n}
		\arrow["{f_\a}", from=1-2, to=2-2]
		\arrow["0"', from=2-1, to=2-2]
		\arrow["i", hook,  from=1-1, to=1-2]
		\arrow[from=1-1, to=2-1]
		\arrow["\lrcorner"{anchor=center, pos=0.125}, draw=none, from=1-1, to=2-2]
	\end{tikzcd}.\]
Note that the target 
of $f_\a$ will have 
the same rank for any 
choice of index, $\a$.
Let $\a_0$
denote a fixed index value.
Observe that $\N_{\Z(f_{\a_0})/U_{\a_0}} 
\simeq \O_{\Z(f_{\a_0})}^{\oplus n}$,
so the normal bundle 
is trivial of rank $n$.
To $f_{\a_0}$ we associate 
the corresponding regular 
function $\widetilde{f_{\a_0}}$ on 
$U_{\a_0} \times {\An}^{\vee}$.
This is the setting of
\cref{theorem: local statement},
so we obtain the isomorphism,
\begin{equation}
\label{iso: a}
g_{\a_0}: \varphi_{\widetilde{f_{\a_0}}}(\pr_1^*\const_{U_{\a_0}}[\dim X+n]) 
\xrightarrow{\simeq} \mu_{\Z(f_{\a_0})}(\const_{U_{\a_0}}[\dim X+n]).
\end{equation}

\subsubsection{}
The following lemmas show 
that the left-\ and right-hand
sides of $g_{\a_0}$ are the restrictions to 
$\T^*[-1]\Z \times_{\Z} \Z(f_{\a_0})$
of the perverse sheaves, $\varphi_{\T^*[-1]\Z}$ 
and $\mu_{\Z}[n](\const_X[\dim X])$, respectively. 

\blem
\label{lemma: cotangent of open}
Suppose that $j: \Z_1 \hook \Z_2$ is an open immersion
of derived Artin stacks. 
Then $\T^*[-1]\Z_2 \times_{\Z_2} \Z_1 \simeq \T^*[-1]\Z_1$
\elem

\bproof
Observe that $\mathbb{L}_j \simeq 0$, so
the natural map $\O_{Z_1} \otimes_{\O_{\Z_2}} 
\mathbb{L}_{\Z_2} = j^*\mathbb{L}_{\Z_2} 
\to \mathbb{L}_{\Z_1}$ is an isomorphism. 
The result follows after shifting, dualizing,
and applying the total space construction.
\eproof

\blem
\label{lemma: restriction of phi and mu}
Suppose that $j: \Z_1 \hook \Z_2$ is an open immersion
of derived Artin stacks. Let $\T_j: \T^*[-1]\Z_1 \to \T^*[-1]\Z_2$
be the map of stacks induced by projection on the first
factor of $\T^*[-1]\Z_1 \times_{\Z_1} \Z_2$ via
\cref{lemma: cotangent of open}.
Then
\begin{align}
	\label{iso: restrict phi} \varphi_{\T^*[-1]\Z_1} 	&\simeq T_j^*\varphi_{\T^*[-1]\Z_2} \, and \\
	\label{iso: restrict mu} \mu_{\Z_1} 			&\simeq T_j^*\mu_{\Z_2}.
\end{align}
\elem

\bproof
By \cite[Corollary 6.11]{BBDJSS15},
the perverse sheaf $\varphi_{\T^*[-1]\Z_1}$ 
is defined uniquely up to canonical isomorphism
as the perverse sheaf possessing the characteristics
(i) and (ii) of \cite[Theorem 6.9]{BBDJSS15}. 
Thus, to prove \labelcref{iso: restrict phi}, it
suffices to show that $T_j^*\varphi_{\T^*[-1]\Z_2}$
possesses these characteristics. 

Observe that if $\mathscr{R}$
is a critical chart for the canonical
d-critical locus associated to $(\T^*[-1]\Z_1,
\omega_{\T^*[-1]\Z_1})$, it is also a
critical chart for the d-critical
locus $T^*[-1]\Z_2$, viewing
$j_R: R \hook T^*[-1]\Z_1$ as an
open subset of $T^*[-1]\Z_2$ via
composition with the open embedding 
$T_j^*$, denoted $j'_R: R \hook T^*[-1]\Z_2$.
As such, ${\varphi_{\T^*[-1]\Z_2}|}_R 
:= {j'_R}^*\varphi_{\T^*[-1]\Z_2} 
\simeq {(T_j^*\varphi_{\T^*[-1]\Z_2})|}_R$.
Both characteristics (i) and (ii)
now follow from the corresponding
characteristics of $\varphi_{\T^*[-1]\Z_2}$.

The isomorphism \labelcref{iso: restrict mu}
is obtained as the composition of the
smooth base change isomorphism along
$j$ for quasi-smooth specialization 
(\cref{lemma: smooth base change for specialization}) 
and the general base change isomorphism for
the Fourier--Sato transform 
(\cref{lemma: base change for FS transform}).
\eproof

\subsubsection{}
Letting $\a$ range over 
all index values, we obtain 
a collection of isomorphisms, 
$\{\varphi_{\widetilde{f_\a}}
(\pr_1^*\const_{U_{\a_0}}[\dim X+n]) \simeq 
\mu_{\Z(f_\a)}(\const_{U_{\a_0}}[\dim X+n])\}_{\a \in \mathcal{I}}$,
whose sources and targets glue
to the perverse sheaves 
$\varphi_{\T^*[-1]\Z}$
and $\mu_{\Z}[n](\const_X[\dim X])$ by 
\cref{lemma: cotangent of open}
and \cref{lemma: 
restriction of phi and mu}.
We finally prove
\cref{theorem: global statement}.

\bproof[Proof of \cref{theorem: global statement}]
\label{proof of global statement}
Let all notation and terms 
be as above.
Since $\Perv(T^*[-1]\Z)$ is a 
groupoid-valued stack in the 
\'etale topology of $T^*[-1]\Z$\footnotemark, 
\footnotetext{\cite[Theorem 2.7]{BBDJSS15}}
it suffices to show that the 
isomorphisms $g_\a$ agree 
on overlaps $\{U_{\a\b}\}$ for 
any values of $\a$ and $\b$.
We establish some 
alternative notation to ease
the exposition. 
Denote $\Z_\a := \Z(f_\a)$ and
$\Z_{\a\b} := \Z_\a \times_{\Z} 
\Z_\b$. Let $\bm{i}^\a_{\a\b}:
\Z_{\a\b} \to \Z_\a \to \Z$ 
denote the obvious inclusion. 
By abuse of notation, we also
use $\bm{i}^\a_{\a\b}$ to denote 
the corresponding map
$\T^*[-1]\Z_{\a\b} \to 
\T^*[-1]\Z_\a \to \T^*[-1]\Z$.
Now choose arbitrary $\a$ and $\b$.
We proceed to show that the following
equality holds:
\begin{equation}
\label{equation: compatibility} 
(i^\a_{\a\b})^*g_\a = (i^\b_{\a\b})^*g_\b,
\end{equation}
where $g_\a$ and $g_\b$ are seen as maps
in $\Perv(T^*[-1]\Z)$ using extension by $0$. 

Note that the Zariski open 
$U_\a \subset X$, and the pull-back diagram, 
	\[\begin{tikzcd}
		{\bm{Z}(f_\a)} & {U_\a} \\
		{*} & {\mathbb{A}^n}
		\arrow["{f_\a}", from=1-2, to=2-2]
		\arrow["0"', from=2-1, to=2-2]
		\arrow["i", hook,  from=1-1, to=1-2]
		\arrow[from=1-1, to=2-1]
		\arrow["\lrcorner"{anchor=center, pos=0.125}, draw=none, from=1-1, to=2-2]
	\end{tikzcd},\]
comprise a good Kuranishi chart for $\Z$,
which we denote by $\K_\a$.
Choose an arbitrary point $x \in \Z_{\a\b}$.
By \cref{lemma: compare charts}
there exists an third good Kuranishi chart 
$\K'_x$ around $x$ which openly embeds
into each of $\K_\a$ and $\K_\b$.
The collection of all such charts as
$x$ ranges over all points 
in $\Z_{\a\b}$ forms a
Zariski open cover of $\Z_{\a\b}$.
Denote by $\bm{i}'_x$ both the open
immersion $\Z'_x \hook \Z_{\a\b}$
and the induced open immersion
$\T^*[-1]\Z'_x \hook \T^*[-1]\Z_{\a\b}$. 
Since the $\{\Z'_x\}$ form an open
cover of $\Z$, in order to show
\labelcref{equation: compatibility}
it suffices to show the equality
\[(i'_x)^*(i^\a_{\a\b})^*g_\a 
	= (i'_x)^*(i^\b_{\a\b})^*g_\b,\]
for arbitrary $x$.

The equality,
$i^\a_{\a\b} \circ i'_x = \i'$,
is merely a reformulation
of \cref{lemma: compare charts},
property (\labelcref{property 2: compare charts}).
We obtain the induced equality
on the level of $-1$-shifted cotangent
bundles, denoting (again by abuse of 
notation) the induced map
$T^*[-1]\Z'_x \hook T^*[-1]\Z$
by $\i'$. As such,
	\[\begin{tikzcd}
		{(i'_x)^*(i^\a_{\a\b})^*g_\a = (i'_x)^*(i^\b_{\a\b})^*g_\b} & {\i'^*g_\a = \i'^*g_\b},
		\arrow[Leftrightarrow, from=1-1, to=1-2]
	\end{tikzcd}\]
and the latter equation is true
since all of the constructions
in section \cref{sec: the local equivalence}
are self-evidently compatible with open
restrictions. That is, the isomorphism,
$\varphi_{\widetilde{f|_U}} 
\xrightarrow{\simeq} 
\mu_{\Z({f|}_U)} \in \Shv(\Z({f|}_U) 
\times V^{\vee})$,
obtained in the proof
of \cref{theorem:
local statement} for a Zariski open, 
$U \subset X$, is naturally the restriction of
the isomorphism $\varphi_{\widetilde{f}} 
\xrightarrow{\simeq} \mu_{\Z(f)}$
along $Z({f|}_U) \times V^{\vee} 
\simeq {Z(f)|}_U \times V^{\vee} 
\hook Z(f) \times V^{\vee}$.
Indeed, this follows from the following
canonical isomorphism of sections:
${\F|}^*_U(V) \simeq \F(U \cap V)$.

Since the choices of $\a$ and $\b$
in \labelcref{equation: compatibility}
were arbitrary, this concludes the proof.
\eproof

\subsubsection{}
Since the canonical sheaf 
$\varphi_{\T^*[-1]\Z}$
is intrinsic to $\Z$,
\cref{theorem: global statement}
immediately implies the following
corollary.

\begin{cor}
\label{cor: mu is intrinsic}
Suppose $\Z$ is a quasi-smooth
derived scheme and $\Z \hook X$ is
a closed immersion into a smooth
scheme $X$. Then
$\mu_{\Z/X}[\vircodim(\Z, X)](\const_{X}[\dim X]) 
\in \Perv(\T^*[-1]\Z)$, is
independent of the closed 
immersion $\Z \hook X$.
In particular, if $\Z \hook X'$ is another
closed immersion, there is a canonical
isomorphism,
\begin{equation}
\mu_{\Z/X}[\vircodim(\Z,X)](\const_{X}[\dim X]) 
	\simeq \mu_{\Z/X'}[\vircodim(\Z,X')](\const_{X'}[\dim X']).
\end{equation}
\end{cor}

In light of \cref{cor: mu
is intrinsic}, we make the following
definition for any quasi-smooth derived
scheme which admits a closed immersion
into a smooth scheme. 

\bdef
The \emph{microlocal homology sheaf} on $\Z$
is defined to be:
\[\mu_{\Z} := \mu_{\Z/X}[\vircodim(\Z,X)](\const_{X}[\dim X]),\]
for any choice of closed immersion,
$\Z \hook X$, into a smooth scheme, $X$.
\edefn

\brem
It should be possible to modify 
the given proof of \cref{theorem: global statement}
to prove a variant in which $\Z$ is not assumed to
admit a closed immersion into
a smooth scheme, $X$. 
Using the derived 
microlocalization functor of \cite{AAK}\footnotemark
\footnotetext{See \cref{rem: Adeel's work}.},
one might obtain such a variant by replacing
the right-hand side of the isomorphism
in \cref{theorem: global statement} by
$\mu_{\Z/\ast}(e)$,
which is manifestly independent of
any embedding $\Z$ might admit. 
In the case that $\Z$
happens to admit an embedding into a smooth
scheme $X$, $\mu_{\Z/\ast}(e)$ and $\mu_{\Z}$
would necessarily agree.
Such a generalization of our present theorem is the subject
of ongoing work in progress by Adeel Khan and Tasuki Kinjo.

\erem


\section{Microlocal homology}
\label{sec: microlocal chains}
Throughout this section, let $e = \mathbb{C}$.
Consider the derived
pullback square,
\[\begin{tikzcd}
	\Z(f) & X \\
	\ast & V
	\arrow["0", hook, from=2-1, to=2-2]
	\arrow[from=1-1, to=2-1]
	\arrow["f", from=1-2, to=2-2]
	\arrow["i", hook, from=1-1, to=1-2]
	\arrow["\lrcorner"{anchor=center, pos=0.125}, draw=none, from=1-1, to=2-2],
\end{tikzcd}\]
where $X$ is smooth of
pure dimension $m$, $V$ is an 
$n$-dimensional $\C$-vector 
space viewed as a scheme,
and $\Z(f)$ is the derived
zero fiber of a propert map $f$. Let $Z(f)$
denote the underlying classical
scheme of $\Z(f)$---equivalently,
the \emph{classical} zero fiber
of $f$. In a graduate course
taught during the Fall 2013
quarter at UC Berkeley, David
Nadler introduced 
the following concept.

\bdef[cf.\ \cite{Y13}]
\label{def: microlocal homology}
Let $\Lambda \subset V^{\vee}$
be a closed conic subset. Then the 
\emph{$\Lambda$-microlocal homology} of $f$
with support $\Lambda$ is
defined to be
\[H_{\bullet}^{\Lambda}(f) := \Gamma_{\Lambda}
	(V^{\vee}; \mu_{0/V}(f_*\underline{\C}_{X}[m])),\]
where $\Gamma_{\Lambda}(V^{\vee};-)$ is
the global sections of the ``sections with support in
$\Lambda$" functor, using the 
terminology of \cite{KS90}---alternatively, 
the functor ${pt}_*
(\Lambda \hook V^{\vee})^!$. 
\edefn

\brem
In \cref{def: microlocal homology},
one can take $\mu_{0/V}$ to
be the classical
microlocalization of \cite{KS90}
or the quasi-smooth 
microlocalization---they
coincide here.
\erem

The microlocal homology is
meant to interpolate between
the singular cohomology and
Borel--Moore homology of $Z(f)$
(with the complex analytic topology)
as $\Lambda$ ranges through
the poset of conic subsets of $V^{\vee}$.
Indeed, supposing that $X$ is oriented,
a calculation using the standard
properties of $\mu_{0/V}$ yields
the equivalences,
\begin{align}
H_{\bullet}^{\{0\}}(f) 			&\simeq \Gamma(Z(f); \underline{\C}_{Z(f)})[m-n] \\
H_{\bullet}^{V^{\vee}}(f) 			&\simeq \Gamma(Z(f); \omega_{Z(f)}),
\end{align}
meaning that global sections
of the microlocal homology
recover, up to a shift, 
singular cochains on $Z(f)$ when
$\Lambda = 0$ and Borel--Moore
chains on $Z(f)$ when 
$\Lambda = V^{\vee}$.
A natural line of inquiry is
the following.

\begin{quest}
\label{question: nice description}
What is a nice description of
the intermediate complexes
given by taking $0 \subsetneq \Lambda
\subsetneq V^{\vee}$?
\end{quest}

\subsection{Geometric chains}
As discussed in the introduction, 
on an oriented smooth manifold, singular
chains and singular cochains are,
in a sense, the same type of
mathematical object: Poincar\'e
duality furnishes an isomorphism
between Borel--Moore homology 
and singular cohomology groups.

Though $Z(f)$ may in general
have singularities, it is a tame
space in the sense that $X$ and
$V$ can always be Whitney stratified
in such a way that $f$ is a map of
stratified spaces and $Z(f)$ is a union
of strata. On such a space, the
complex of singular chains and
singular cochains are also not totally
different species of mathematical
objects: both may be represented
by \emph{geometric chains} and
\emph{geometric cochains}, which
satisfy a kind of Poincar\'e duality. 
Both geometric chains and cochains
are, roughly speaking, formal linear 
combinations of (suitable) subsets 
of $Z(f)$, together with some additional
data (e.g.\ (co)orientation); the precise
definition may be found in \cite{G81}. A geometric
$k$-chain is one whose formal terms
are $k$-dimensional subsets of $Z(f)$,
and a geometric $k$-\emph{co}chain is one whose
formal terms are $k$-\emph{co}dimensional
subsets of $Z(f)$. Like in the
case of a smooth manifold, there
is a means of taking cocycles
to cycles by taking the intersection
of the ``support" of a cocycle
(the union of the subsets making
up the formal sum) with the support
of a cycle---giving the cap product
in singular homology.

The animating idea behind the ``geometric cycle" 
description of homology and 
cohomology---to view chains
and cochains as concrete subsets of our space
(with some extra data and decoration)---proved
enormously fruitful to Mark Goresky 
and Robert MacPherson in their development
of intersection homology. 
In their search for a robust homology
theory for singular spaces, Goresky and
MacPherson decided to look at the chain
theories obtained by considering only
\emph{some} allowable geometric chains. If we
work with $Y \subset X$, the closure of
a single stratum inside the Whitney stratified
manifold, $X$, allowable
geometric chains are 
those whose intersections
with the various strata of $Y$
have dimensions dictated by a chosen
perversity function\footnotemark.
\footnotetext{An integer-valued function on the
collection of singular strata of $Y$. 
See \cite[\S1.1]{M90}.} 
There are two extreme
perversity functions, the zero perversity,
given by the zero function, and
the so-called ``top" perversity function, 
given by the assignment,
\[Y_{\a} \mapsto \codim Y_{\a} -2,\]
for any stratum $Y_{\a} \subset Y$,
and where, if $Y$ is of pure real dimension $n$, 
$\codim Y_{\a} := n - \dim Y_{\a}$. 
By \cite[\S4.3]{GM80}, the 
chain theory determined by using allowable
geometric chains for the top
perversity is singular homology, and
the chain theory determined 
by allowable chains for the 
zero perversity is singular 
cohomology. Thus, intersection
homology interpolates between
singular homology and cohomology,
giving different geometric homology
theories as the perversity function varies.

The central lesson of the 
Goresky-MacPherson construction
of intersection homology is that
considering collections 
of chains with fixed 
geometric properties leads to 
interesting theories of homology.

\subsection{Currents}

Returning to the microlocal
homology, let us rewrite the
Borel--Moore homology and
the singular cohomology of $Z(f)$
using resolutions by currents 
on $X$. If $\mathscr{D}
^{\bullet}_X$ denotes the sheaf of
currents\footnotemark on $X$ (i.e.\ distributional
differential forms), then
\begin{align}
\underline{\C}_{Z(f)} 			&\simeq i^*\mathscr{D}^{\bullet}_X  \label{currents resolution 1} \\
\omega_{Z(f)}[-m] 	&\simeq i^*\mathscr{D}^{\bullet}_{[Z(f)]} \label{currents resolution 2},
\end{align}
	\footnotetext{The theory of currents has been
	developed mostly in the context of geometric
	measure theory, where they play the role of
	manifolds in differential geometry. For details 
	of the theory, the reader is referred to \cite{Fed69}, 
	the standard reference on the subject.}
where $\mathscr{D}^{\bullet}_{[Z(f)]}$
denotes the sheaf of currents on $X$ supported
on $Z(f) \subset X$.
Since $X$ is smooth,
$\Omega^{\bullet}_X$
is equivalent to 
$\mathscr{D}^{\bullet}_X$\footnotemark,
so \labelcref{currents resolution 1} 
follows from the
algebraic de Rham theorem.
	\footnotetext{See pg.\ 385 of \cite{GH78} 
	for a proof of the
	so-called ``smoothing of cohomology."}
On the other hand,
\labelcref{currents resolution 2}
follows from the following 
slightly more general lemma.

\blem
\label{lemma: distributional dR for BM homology}
Suppose that $i: Z \hook X$ is
the inclusion of a closed subset
into a smooth, oriented manifold $X$
of dimension $m$. Then,
\[\omega_Z[-m] \simeq i^*\mathscr{D}^{\bullet}_{[Z]}.\]
\elem

\bproof
Let $j: U \hook X$ denote the
inclusion of the open subset
complementary to $Z$. By the
standard adjunction triangle,
we obtain $i_*i^!\underline{\C}_{X}
\simeq i_*i^!(\omega_X[-m])
\simeq i_*\omega_{Z}[-m]$ as
the fiber of the map, $\underline{\C}_X
\to j_*\underline{\C}_U$. Since $U$ and
$X$ are smooth, $\underline{\C}_X$ and
$\underline{\C}_U$ have standard resolutions
by the sheaves of currents,
$\mathscr{D}^{\bullet}_X$ and
$\mathscr{D}^{\bullet}_U$,
respectively. The sheaf of
$p$-currents, $\mathscr{D}^p_U$
is fine for each $p$, so the
push-forward $j_*\underline{\C}_U$ has
the obvious resolution given
by the pushforward, 
$j_*\mathscr{D}^{\bullet}_U$.
Now the map, $\underline{\C}_X \to j_*\underline{\C}_U$,
induces the obvious map on
resolutions $\mathscr{D}^{\bullet}_X
\to j_*\mathscr{D}^{\bullet}_U$
given by domain restriction 
of distributions, and its fiber is
clearly the sheaf of currents 
supported on $Z$.
\eproof
Using such resolutions, we may
express the microlocal homology
at $\Lambda = \{0, V^{\vee}\}$
using the pullback of certain currents on $X$:
\begin{align}
H_{\bullet}^{\{0\}}(f) 			&\simeq \Gamma(Z(f); i^*\mathscr{D}^{\bullet}_X)[m-n] \\
H_{\bullet}^{V^{\vee}}(f) 			&\simeq \Gamma(Z(f); i^*\mathscr{D}^{\bullet}_{[Z(f)]}[m]).
\end{align}

While we have taken the perspective
that currents are merely differential
forms with distributional coefficients,
in most standard accounts of currents
(e.g.\ \cite{Fed69}) they are 
defined as the continuous dual
to the topological vector space of compactly
supported, smooth differential forms.
The latter perspective suggests a connection
between currents and geometric chains\footnotemark
	\footnotetext{Indeed, the author learned from 
	\cite{DP07} that currents were studied 
	by Georges de Rham in order to ``show the 
	connections between differential forms and 
	singular chains" (\cite[pg. 168]{DP07}).} 
that can be made precise. For example,
integration against a piece-wise smooth, oriented
$p$-chain in $X$ gives a $p$-current\footnotemark.
	\footnotetext{Beware notation conventions 
	in the literature for currents.} 
In particular, any non-degenerate $p$-simplex 
determined by a general choice of $p$ points
in $X$ gives a $p$-current. In fact, geometric
measure theorists have defined the theories
of ``integral flat" homology and ``integral
rectifiable" homology based on special classes
of currents called ``integral flat chains" and
``integral currents," respectively\footnotemark. 
	\footnotetext{The interesting article, \cite{DP07}, compares 
	these homology theories with the theories 
	of \v{C}ech and singular homology.}

In light of this discussion, 
the philosophy behind 
the intersection homology
of Goresky-MacPherson 
suggests that we address 
\cref{question: nice description} by
fixing certain conditions indexed by $\Lambda$
and by considering only those
currents which satisfy these conditions.  
Inspired by this line of thought, 
the author's graduate thesis 
advisor, David Ben-Zvi, made the
following rough conjecture.

\begin{conj}[Ben-Zvi c.\ 2016]
\label{conj: David's conjecture v1}
There should exist a resolution
of $H_{\bullet}^{\Lambda}(f)$ by the pullback of
currents on $X$ with prescribed 
regularities (e.g. wave front sets contained in
a prescribed conic, subset of the cotangent
bundle) dictated in some way by $\Lambda$.
\end{conj}

Of course, we can view currents as
differential forms with distributional
coefficients as we did initially.
From this perspective, \cref{conj: David's
conjecture v1} states that 
there should be a $\Lambda$-indexed 
family of ``distributional de Rham
theorems" for the microlocal
homology that interpolates
between the usual de Rham theorem
when $\Lambda = 0$ and the one
proven in \cref{lemma: distributional
dR for BM homology}.

\subsection{Singular support of Borel--Moore chains}
On a smooth manifold, the 
\emph{cochain} theories obtained by using 
currents of any specified
regularity---i.e.\ whether we consider
all currents, only currents with Lipschitz
regular coefficients, only currents
with smooth coefficients, etc.\--are 
equivalent. While we can model
the de Rham complex with currents,
the regularity (i.e.\ wave front sets) of
individual currents is not meaningful
since every current can be represented
by a smooth form in its cohomology
class. To the extent that currents
can be viewed as chains, this is
saying that any chain on a smooth
manifold always has a smooth
chain representative in homology.
This suggests that, for a singular 
space $Z \subset X$ sitting inside 
of an ambient smooth $X$, when
we model homology classes of
chains on $Z$ using currents
on $X$, there should be an
obstruction to representing these
currents by smooth forms which
model the same homology class
of $Z$. Furthermore, this obstruction
should behave a like a theory of
singular support, assigning to each
chain a conic subset of some vector
space or vector bundle such 
that an assignment of the zero
section implies smoothness 
of (a representative of) the chain. 

When $Z:=Z(f)$, the microlocal
homology provides a systematic
way of producing sheaves of currents
with various regularity conditions
(assuming some kind of positive
resolution to \cref{conj: David's 
conjecture v1}). At one extreme,
$\Lambda = 0$ produces only
currents with smooth coefficients.
At the other extreme, $\Lambda = V^{\vee}$
produces \emph{all} currents. The
intermediate values for $\Lambda$ produce
classes of currents somewhere in
between these two, and as $\Lambda$
becomes larger, so too 
does the collection of currents
it classifies. This suggests that
$\Lambda \subset V^{\vee}$
controls the behavior of chains
in a manner similar to 
the theory of singular support
we seek. A working definition
of the singular support of a current,
$T$, might now be: the minimal closed, 
conic subset, $\Lambda$, such that $T
\in M^f_{\Lambda}$.
This train of thought leads
to another conjecture of David
Ben-Zvi, which reformulates
some of the ideas behind 
\cref{conj: David's conjecture v1}.

\begin{conj}[Ben-Zvi c.\ 2020]
\label{conj: David's conjecture v2}
There is a good notion of 
singular support for 
chains on $Z(f)$ which
measures the failure of a
chain to be represented by
``smooth chains" on $X$.
\end{conj}

\subsubsection{Main motivation}

Clearly, in the remarks
and discussion leading up
to \cref{conj: David's conjecture v2}, 
the central object has 
been, not $f$, but
rather the closed immersion
$i: Z(f) \hook X$. Yet, the failure
of a chain in $Z(f)$ to
be represented by smooth
chains in $X$ is ultimately
the failure of $Z(f)$ to be smooth. 
Indeed, if a chain in $Z(f)$ avoids the
singular locus, it is always
representable by a smooth
chain in $X$.

Using elementary homological
algebra and the standard sequence
on cotangent complexes,
we obtain the following 
closed immersion of 
schemes over $Z(f)$,
\begin{equation}
\label{equation: inclusion}
T^*[-1]\Z(f) \hook Z(f) \times V^{\vee}.
\end{equation}
Ideally, the singular support of
a chain on $Z(f)$ would live
inside of a vector bundle
over $Z(f)$ itself, so rather
than consider $\Lambda \subset
V^{\vee}$ we instead consider
$Z(f) \times \Lambda \subset
Z(f) \times V^{\vee}$. The following
conjecture now represents a important
step toward showing a good theory
of singular support of chains.

\begin{conj}
\label{conj: microlocal homology depends on scheme of singularities}
Let $\widetilde{\Lambda} := 
T^*[-1]\Z(f) \cap (Z(f) 
\times \Lambda)$, for $\Lambda 
\subset V^{\vee}$. Let $\Lambda, 
\Lambda' \subset V^{\vee}$ be such that $\widetilde{\Lambda} = 
\widetilde{\Lambda'}$. Then
\[H_{\bullet}^{\Lambda}(f) \simeq H_{\bullet}^{\Lambda'}(f).\]
\end{conj}

Roughly speaking, this conjecture 
states that the microlocal homology
should only depend on the part of $\Lambda$
that intersects $T^*[-1]\Z(f)$
(the ``scheme of singularities" of $\Z(f)$
in the language of \cite{AG15}) 
when made into a trivial family over
$Z(f)$,  viewed as a subfamily of the
conormal bundle to $\Z(f)$. It suggests
that the microlocal homology
is actually something intrinsic to
the quasi-smooth derived scheme
$\Z(f)$ rather than to the map
$f$ itself. 


\Cref{theorem: global statement}
allows us to make the following definition,
and answer \cref{conj: microlocal homology
depends on scheme of singularities} in the
affirmative in the case when $f$ is proper. 

\bdef
\label{def: intrinsic microlocal homology}
Given a quasi-smooth derived
scheme $\Z$ which admits
a closed immersion into a
smooth scheme $X$, the 
\emph{$\Lambda$-microlocal
homology sheaf of $\Z$} is defined to
be the sheaf,
\begin{equation}
\label{def: intrinsic microlocal homology}
\mu_{\Z}^{\Lambda} := 
	\Gamma_{\Lambda}\varphi_{\T^*[-1]\Z},
\end{equation}
for $\Lambda \subset T^*[-1]\Z$ 
closed and conical.
\edefn

When $\Z = \Z(f) \hook X$ is the derived
zero locus of a proper map
$f: X \to V$, the following theorem
justifies our terminology by
showing that the global sections of the
microlocal homology of $\Z(f)$ recovers 
the microlocal homology
of $f$.

\bthm
\label{theorem: equivalence of microlocal homologies}
Let $\Lambda \subset V^{\vee}$ be
a closed conic subset, and suppose $f: X \to V$
is a proper map from a smooth scheme $X$ to
a finite dimensional vector space $V$.
Then there exists an equivalence,
\[\Gamma(T^*[-1]\Z(f); \mu_{\Z(f)}^{\widetilde{\Lambda}}) 
	\simeq H_{\bullet}^{\Lambda}(f)[\vircodim \, \Z(f)].\]
\ethm

\bproof
Since the inclusion $i: \Z(f) 
\hook X$ is quasi-smooth, it
has a well-defined normal bundle,
as described in \cref{sec: deformation to
the normal bundle}. Furthermore, since
the diagram,
\[\begin{tikzcd}
	{\Z(f)} & X \\
	\ast & {V^{\vee}}
	\arrow["0", from=2-1, to=2-2]
	\arrow[from=1-1, to=1-2]
	\arrow["f", from=1-2, to=2-2]
	\arrow["{f'}"', from=1-1, to=2-1]
	\arrow["\lrcorner"{anchor=center, pos=0.125}, draw=none, from=1-1, to=2-2]
\end{tikzcd}\]
is derived Cartesian, $f$ induces a
map of bundles over $Z(f)$, which
we call $T_f: \simeq N_{\Z(f)/X} \to 
{f'}^*N_{0/V^{\vee}} \simeq Z(f) \times V$---upon 
the trivialization $N_{\Z(f)/X} \simeq Z(f) \times V$,
$T_f$ becomes the identity.
Composition of $T_f$ with the projection,
$\pr: Z(f) \times V \to V$, gives the
 map of normal bundles,
$N_f: N_{\Z(f)/X} \to N_{0/V}$,
described in \cref{sec:
deformation to the normal bundle}.

Let $T^{\vee}_f$ and $\pi$ denote the
maps on dual bundles induced
by $T_f$ and $\pr$, respectively.
Note the Cartesian diagram,
\begin{equation}
\label{diagram: lambda}
\begin{tikzcd}
	{Z(f) \times \Lambda} & {Z(f) \times V^{\vee}} \\
	\Lambda & {V^{\vee}}
	\arrow["{\pi'}", from=1-1, to=2-1]
	\arrow["i", hook, from=2-1, to=2-2]
	\arrow["\pi", from=1-2, to=2-2]
	\arrow["{i'}", hook, from=1-1, to=1-2]
	\arrow["\lrcorner"{anchor=center, pos=0.125}, draw=none, from=1-1, to=2-2].
\end{tikzcd}
\end{equation}
The following chain
of equivalences proves the theorem.
\begin{align}
pt_*\mu_{\Z(f)}^{\widetilde{\Lambda}}[-\vircodim \, \Z(f)] 						&\simeq pt_*{i'}^*\mu_{\Z(f)/X}(\const_X) \label{big boy} \\
																&\simeq pt_*\pi'_*{i'}^!\mu_{\Z(f)/X}(\const_X) \\
																&\simeq pt_*i^!\pi_*\mu_{\Z(f)/X}(\const_X) \label{base change} \\
																&\simeq pt_*i^!\pi_*{T^{\vee}_f}^!\mu_{\Z(f)/X}(\const_X) \label{tee vee} \\
																&\simeq pt_*i^!\pi_*{T^{\vee}_f}^!\Four(\Sp_{\Z(f)/X}(\const_X)) \\
																&\simeq pt_*i^!\Four(\pr_*{T_f}_*\Sp_{\Z(f)/X}(\const_X)) \\
																&\simeq pt_*i^!\Four({N_f}_*\Sp_{\Z(f)/X}(\const_X)) \\
																&\simeq pt_*i^!\mu_{\Z(f)/X}(f_*\const_X) \label{base change 2} \\
																&=: H_{\bullet}^{\Lambda}(f).
\end{align}
\Cref{theorem: global statement}
implies the first equivalence, 
\labelcref{big boy}; 
base change along diagram
\labelcref{diagram: lambda}
gives \labelcref{base change};
\labelcref{tee vee} follows from 
$T^{\vee}_f$ being the identity;
and \labelcref{base change 2} follows from
\cref{lemma: proper base change
for specialization}.
\eproof

\begin{cor}
\cref{conj: microlocal homology depends
on scheme of singularities} is true in the 
case when $f$ is proper.
\end{cor}

\Cref{theorem: equivalence of microlocal 
homologies} also allows us to extract from 
$\mu_{\Z(f)}^{\Lambda}$ a plausible singular
support theory for Borel--Moore chains when $f$ is proper,
thereby making progress toward resolving 
\cref{conj: David's conjecture v2} in
an important special case.

\begin{spec}[``Singular support for Borel--Moore chains"]
Let $f$ be as in \cref{theorem: equivalence of microlocal homologies}.
The singular support for an element $[\sigma] \in H_{\bullet}^{\BM}(Z(f))$
is the minimal subset $\Lambda \subset T^*[-1]Z(f)$ such that
$[\sigma]$ lies in the image of the natural map,
\begin{align*}
\Gamma\left(\mu_{\Z(f)}^{\Lambda}[-\vircodim \, \Z(f)]\right) 	&\to \Gamma\left(\mu_{\Z(f)}^{T^*[-1]\Z(f)}[-\vircodim \, \Z(f)]\right) \\
												&\simeq H_{\bullet}^{\BM}(Z(f)),
\end{align*}
where where have denoted global sections
by $\Gamma$ in order to ease the notation.
\end{spec}

\brem
In order to show that this proposed definition
of singular support for chains is reasonable, we
would need to show, for example, that a minimal
subset $\Lambda \subset T^*[-1]Z(f)$ exists. We
plan to address this in future work.
\erem

\appendix


\section{Sheaf theory}
\label{sec: sheaf theory}
\subsection{Sheaves valued in $\Mod_e$}

The main reference for sheaves of
spectra is \cite{SAG}. Further details
or justifications for unproven technical
assertions in this subsection can
be found there.

\subsubsection{}
Suppose that $X$ is a complex
analytic space, locally of finite type.
By abuse of notation, $X$ will denote
both the complex analytic space and
its underlying topological space.
Denote by $\U(X)$ the partially
ordered collection of open sets of
$X$, viewed as a category.

\subsubsection{The category of sheaves}
 
Let $\Mod_e$ be the
$\infty$-category of 
$e$-module objects in
the category of
spectra\footnotemark 
and consider the category 
of $\Mod_e$\-valued
presheaves on $X$, 
\[\PrShv(X; e) := \Fun(\U(X)^{\op},
\Mod_e).\] 

\footnotetext{See 
\cite[\S3.3]{HA} for a general discussion of
modules over an algebra object.}

\bdef
The category of $e$\-sheaves
on $X$, denoted by $\Shv(X; e)$, is
defined as the full subcategory 
of $\Fun(\U(X)^{\op},
\Mod_e)$ spanned by those functors 
$\F$ which satisfy the following descent
condition for open covers:
\begin{itemize}
\item[$(\ast)$] 
Let $U \in \U(X)$ 
and let $\{U_{\a} \to U\}$ be a
covering of $U$, then the canonical
map,
\[\F(\U) \to \varprojlim_{V}\F(V)\]
is an equivalence in $\Mod_e$, where
$V$ ranges over all opens contained
in $U$ whose inclusion factors through 
some $U_{\a}$.
\end{itemize}
\edefn

\brem
The above definition makes 
sense for any topological
space, not just those underlying
a complex analytic space.
\erem

We may also define 
more generally the category of 
sheaves on $X$ valued in
an arbitrary $\infty$-category 
$\mathscr{C}$, denoted 
$\Shv(X; \mathscr{C})$, 
using the same descent 
condition (see 
\cite[Definition 1.1.2.1]{SAG}).
An important example for us
will be the category of 
spectral sheaves
on $X$, $\Shv(X; \Sp)$.
The stalk at a point
and the support of a sheaf have 
the same definitions
as in the setting of ordinary categories.

\bdef
Suppose that $\F \in \Shv(X; \Cee)$,
where $\Cee$ is a stable $\infty$-category.
The stalk at the point $x \in X$,
denoted $\F_x$, is defined to be
\[\F_x := \varinjlim_{\{U \in \U(X)|x \in U\}} \F(U).\]
The support of $\F$, denoted
$\supp \F$, is defined to be
\[\supp \F := \overline{\{x \in X| \F_x \neq 0\}}\]
The stalk and support of a section
$s \in \F(U)$ are defined similarly.
\edefn

\brem
The category $\Shv(X;e)$ has
an alternative description as the
category of $\const$\-module
objects in spectral sheaves, where
$\const$ is the constant sheaf
on $X^{\an}$ with stalk $e$, viewed
as a discrete $\mathbb{E}
_{\infty}$\-ring.
\erem

\subsubsection{$t$-structures}

If $\mathscr{C}$ is given as the
stabilization of a presentable
$\infty$-category $\mathscr{C}'$, 
there is a standard $t$-structure
on $\mathscr{C}$ defined by the
property that $\mathscr{C}_{<0}$
is the fiber of the ``infinite loop
space" functor $\mathscr{C} \simeq 
\Sp(\mathscr{C}') \to \mathscr{C}'$.
The category of $\mathscr{C}$-valued 
sheaves on $X$ is itself the
stabilization of $\Shv(X; \mathscr{C}')$
(see \cite[Remark 1.3.2.2]{SAG}), so
$\Shv(X; \mathscr{C})$ inherits a standard
$t$-structure. 

\cite[Proposition 1.3.2.7]{SAG} states
the resultant $t$-structure on spectral
sheaves on $X$ (or more generally, 
on any $\infty$-topos), is compatible with 
filtered colimits (i.e. a filtered colimit of
coconnective objects is coconnective) 
and right-complete. 
The $n$th homotopy sheaf of 
an object $\F \in \Shv(X; \Sp)$,
$\pi_n(\F)$, may be identified (
\cite[Example 1.3.2.4]{SAG}) with the
sheaf of abelian groups given 
by sheafifying the presheaf given by
\[U \mapsto \pi_n(\F(U)).\]

Because $\const$ is connective
as a sheaf of spectral rings, 
the category of $e$\-sheaves, 
$\Shv(X; e)$, admits a $t$-structure 
defined by letting $(\Shv(X; e))_
{\leq 0}$ be the inverse image 
of $(\Shv(X; \Sp))_{\leq 0} \subset 
\Shv(X; \Sp)$ under the forgetful 
functor $\Shv(X;e) \to \Shv(X; \Sp)$.
Like the $t$-structure on $\Shv(X; \Sp)$,
it is right complete and compatible
with filtered colimits by 
\cite[Proposition 2.1.1.1]{SAG}. The heart of
this $t$-structure can be identified
with the abelian category of 
$\const$\-modules in sheaves of sets
(\cite[Remark 2.1.2.1]{SAG}), meaning
the homotopy sheaves of an object
$\F \in \Shv(X;e)$ have the additional
structure of $\const$\-modules.



\subsection{Hyperdescent}
\label{ssec: hyperdescent}

The classical results of
sheaves on locally compact
spaces rely on the property
that a map in the derived category 
(of complexes of sheaves) is an
equivalence if and only if the
induced map on stalks is an
equivalence at each point; and
if and only if it induces 
isomorphisms on 
cohomology sheaves.
While the latter property
essentially defines the
derived category,
on an arbitrary $\infty$-topos,
it is not always the case
that the category of sheaves
valued in an arbitrary 
stable $\infty$-category
enjoys these properties. 

\subsubsection{}
By reducing to the case 
of sheaves of spaces
on $X$, we show that 
\blem
\label{lemma: checking isos}
In the
category $\Shv(X;e)$, a
map is an equivalence:
\begin{enumerate}[(i)]
\item
iff it is $\infty$-connective (i.e. 
induces isomorphisms 
on all homotopy sheaves).
\item 
iff the induced map
on stalks is an equivalence
at each point in $X$.
\end{enumerate}
\elem

Both properties are consequences
of the $\infty$-topos of sheaves
of spaces on $X$ being
hypercomplete.\footnotemark 

\subsubsection{}
For a reasonable class
of topological spaces---those
which are suitably `finite' in
some sense---$\Shv(-; \mathcal{S})$ 
is automatically
hypercomplete. Fortunately,
among them are finite
dimensional complex
analytic spaces locally of finite
type, as we now show
in the next two lemmas.
\footnotetext{See 
\cite[\S6.5]{HTT} for a thorough discussion
of hypercompleteness.}

\blem
\label{lemma: covering dimension of a subset}
Any closed 
subspace $Z \subset S$
of a topological space $S$ 
has covering dimension
less than or equal to that of $S$.
\elem

\bproof
Take any open cover $\{U_i\}$
of the subspace $Z$. 
Each $U_i = U'_i \cap Z$
for some open $U'_i \subset S$.
Complete $\{U'_i\}$ to an 
open cover of $S$ by adding the
complement $S \setminus Z$ to the
collection. Let $n$ be
the covering dimension of $S$,
and take a refinement, $\{V'_{\a}\}$,
of the completed cover of $S$
such that the $(n+2)$\-fold intersection
of pairwise disjoint members
of $\{V'_{\a}\}$ are empty. Now take collection
$\{V'_{\a} \cap Z\}$. This is clearly an open
cover of $Z$, and $(n+2)$\-fold
intersections of pairwise disjoint members
of $\{V'_{\a} \cap Z\}$ are clearly empty.
Finally, it is a refinement of
$\{U_i\}$: take any nonempty 
element of the cover, $V'_{\a} \cap Z$.
Necessarily, $V'_{\a} \subset U'_j$ for some $j$,
since otherwise it would be contained
in $S \setminus Z$, and $V'_{\a} 
\cap Z = \varnothing$.
Thus, $V'_{\a} \cap Z \subset U_j$.
\eproof

\blem
\label{lemma: shv(X) is hypercomplete}
Suppose $X$ is a finite
dimensional complex analytic
space, locally of finite type. 
Then the $\infty$-topos 
$\Shv(X; \mathcal{S})$ is
hypercomplete.
\elem

\bproof
By \cite[Corollary 7.2.1.12]{HTT},
it suffices to show that $\Shv(X; \mathcal{S})$
is locally of homotopy dimension $\leq n$
for some integer $n$.
We recall that this means
that there exists a collection 
of objects $\{\F_{\a}\}$ in $\Shv(X; \mathcal{S})$
which generate $\Shv(X; \mathcal{S})$ under
colimits, such that $\Shv(X; \mathcal{S})_{/\F_{\a}}$
is of homotopy dimension $\leq n$.

An obvious collection of objects which
generates $\Shv(X; \mathcal{S})$ is 
the collection of delta sheaves,
$\{\underline{\ast}_{U}\}$,
where $U$ ranges over $\U(X)$.
Observe that 
$\Shv(X; \mathcal{S})_{/\underline{\ast}_{U}}
\simeq \Shv(U; \mathcal{S})$ via the adjoint
pair $(j_{U !}, j_U^*)$ induced
by the inclusion $j_U$, so it
suffices to show that $\Shv(U; \mathcal{S})$ is
of homotopy dimension $\leq n$ for any
open subset $U$.

In fact, it suffices to show that
$\Shv(U_i; \mathcal{S})$ 
is of homotopy dimension
$\leq n$ for some chosen open cover
$\{U_i\}$ of $X$ because if
$U \subset U'$, the homotopy
dimension of $\Shv(U; 
\mathcal{S})$ is less than
or equal to that of $\Shv(U'; 
\mathcal{S})$. Indeed,
if $\F \in \Shv(U'; \mathcal{S})$ 
is $n$-connective,
$\F' := (U \hook U')_!\F$ 
is $n$-connective
since $(U \hook U')_!$ is left exact.
If  $\F'$ admits a global section 
$\underline{\ast}|_{U'} \to \F'$, then
$\underline{\ast}|_{U} = 
(U \hook U')^*\underline{\ast}|_{U'} 
\to (U \hook U')^*\F' \simeq \F$ is a
global section of $\F$.

Note that $X$ locally takes the form
$\O_{\C^n}/(f_1, \ldots, f_m)$
for some natural numbers $n \leq \dim X$
and $m$ . Choose an open
cover $\{U_i\}$ of $X$ on
each member of which $X$
has such a presentation.
Each open subset $U_i$ of our cover,
embeds into $\mathbb{C}^{n_i}$
as the intersection of an
open subset $V \subset \C^{n_i}$ 
and the closed analytic subset given by
\[\{f_1 = \cdots = f_{m_i} = 0\}.\]
By \cite[Theorem 7.2.3.6]{HTT}, 
$\Shv(U_i; \mathcal{S})$ is of
homotopy dimension $\leq \dim X$
if $U_i$ has covering dimension
$\leq \dim X$, but
this follows immediately from
\cref{lemma: covering
dimension of a subset}, since
any open subset of $\C^{n_i}$
is a smooth manifold of dimension
$n < \dim X < \infty$, and the
covering dimension of smooth
manifolds is their usual dimension.
\eproof

\subsubsection{}
We now prove
\cref{lemma: checking isos}(i) and (ii). 
In each case, we only
prove the ``if" direction since the
``only if" direction is obvious.

\bproof[Proof of \cref{lemma: checking isos}]

Suppose $\F \xrightarrow{\a} \F' \in \Shv(X;e)$ 
were an $\infty$-connective 
morphism, i.e. $\pi_i\F \xrightarrow{
\simeq} \pi_i\F'$ for every 
integer $i$. Since the forget functor
$\oblv: \Shv(X;e) \to \Shv(X; \Sp)$
is conservative and preserves
small colimits (\cite[Proposition 2.1.0.3]{SAG}), 
$\oblv(\a)$ is also $\infty$-connective; 
and it suffices to show that
$\oblv(\a)$ is an equivalence.
By abuse of notation we
refer to $\oblv(\F)$ as $\F$. 
Now the induced map
$\Omega^{\infty}\F 
\to \Omega^{\infty}\F'$
is $\infty$-connective for each
integer $n$ by design,
so it is an equivalence 
by hypercompletion using
\cref{lemma: shv(X) is
hypercomplete}. Combining
\cite[Proposition 1.3.3.3]{SAG}
and \cite[Remark 1.1.3.3]{SAG},
we obtain (i).

Now suppose $\F \xrightarrow{\a} \F' \in 
\Shv(X; e)$ were a morphism
such that the induced maps of stalks,
$\F_x \to \F'_x,$
were equivalences for 
all $x \in X$.
Recall the definition of
the stalk at $x$ as the
filtered colimit of sections
over open sets containing
$x$:
\[\F_x := \varinjlim_{x \in U} \F(U).\]

Because the forgetful functor
to spectral sheaves is conservative
and preserves small colimits, 
$\oblv(\a)$ also induces
isomorphisms on stalks, and
it suffices to prove that that $\oblv(\a)$
is an equivalence.
As remarked earlier, the
$t$-structure on spectra
is compatible with filtered colimits.
This is equivalent via the fiber
sequence (for $S \in \Sp$),
\[\tau_{\geq 1}S \to S \to \tau_{\leq 0}S,\]  
to the truncation
functor $\tau_{\leq 0}$ preserving
filtered colimits. The functor $\tau_{\geq 0}$
is a left adjoint, so it preserves all
colimits; therefore, $\pi_0 S :=
\tau_{\geq 0}\tau_{\leq 0}S$ preserves
filtered colimits, as well as $\pi_i$ for
all integers $i$.

In particular, we obtain that
\[\varinjlim_{x \in U} \pi_i(\F(U)) 
	\xrightarrow{\simeq} \pi_i(\varinjlim_{x 
		\in U} \F(U)) =: \pi_i(\F_x).\]
The left-hand side of this
map is, by definition, the stalk
of the presheaf whose sheafification
is $\pi_i \F$. Sheafification preserves
stalks as a left adjoint, so
the left-hand side is also $(\pi_i \F)_x$.

Thus, we obtain that the maps
\[\pi_i \F \xrightarrow{\pi_i(\a)} \pi_i \F'\]
induce isomorphisms
on stalks, for all integers
$i$. Since each $\pi_i \F$
is a sheaf of sets on a topological
space, $\pi_{\a}$ is an isomorphism
for each $i$. By property (i), proven
above, $\oblv(\a)$ is an
equivalence, and we have
proven (ii).
\eproof

\subsubsection{}
Later in this paper, we will
need to consider sheaves
on closed subspaces of manifolds,
$Z \subset M$.
The following lemma is an immediate
consequence of \cref{lemma:
covering dimension of a subset},
and implies that the conclusions of
\cref{lemma: checking isos}
hold for $\Shv(Z;e)$.
\blem
\label{lemma: hypercompleteness
for closed subsets}
Suppose that $Z$ is a closed subspace
of $\R^n$ for some value of $n$.
Then the $\infty$-topos $\Shv(Z; 
\mathcal{S})$ is hypercomplete.
\elem


\subsection{Bounded variations}

We may also consider sheaves
valued in the category of $e$\-modules
objects in one of the bounded 
variations of spectra---eventually
connective, eventually coconnective,
or bounded spectra: $\Mod^+_e$,
$\Mod^{-}_e$, or $\Mod^b_e$ 
We denote the resultant categories
by $\Shv^+(X; e)$, $\Shv^{-}(X;e)$,
and $\Shv^b(X; e)$, respectively.
The results of the previous
sections all hold mutatis mutandis
for these variants.
The most important category
for us will be $\Shv^b(X;e)$,
which we call the bounded derived 
$\infty$-category of $e$\-sheaves
on $X$. This terminology is 
justified on the ground that 
the homotopy category of 
$\Shv^*(X; e)$ (with its
standard $t$-structure) is
the usual bounded derived
category $D^*(X; e)$,
where $* = +,-,$ or $b$.

\brem
\label{rem: sheaves on schemes}
Suppose that $X$ is a classical scheme,
locally of finite type over $k$. 
Associated to $X$ is the
topological space $X^{\an}$,
defined by endowing
the set of closed
points of $X$ with the
complex analytic topology.
One may further define the
analytification of
$X$ by introducing a particular sheaf
of rings $\O_{X^{\an}}$
on $X^{\an}$ in such a way
that the ringed space $(X^{\an}, 
\O_{X^{\an}})$ is a complex
analytic space, though we will
not need details of this construction
other than to note that the resulting
complex analytic space is also
locally of finite type\footnotemark.
\footnotetext{Details of the complex
analytic topology and 
the analytification of $X$ may 
be found in e.g. \cite{N07}}
In this case, we denote
the corresponding categories
of $e$\-sheaves on $X^{\an}$
by $\Shv^*(X;e)$, as well as refer
to them as categories of sheaves
$X$ rather than on
$X^{\an}$. If $X$ is finite 
dimensional, so is its 
associated complex
analytic space. In that
case, \cref{lemma: checking isos}
holds for $\Shv^*(X;e)$.
\erem

\bdef
If $\X$ is a derived $k$-scheme
whose underlying classical scheme
$X$ is locally of finite type, we define
\[\Shv(\X;e) := \Shv(X;e)\]
\edefn


\subsection{Six functor formalism}
\label{ssec: six functor formalism}

In the sheaf theory context
we have just defined, there
is the usual Grothendieck six functor
formalism.

\subsubsection{}
Let $f:X \to Y$ be a map
of locally compact topological 
spaces (not necessarily complex
analytic spaces). Then 
there exist the
following functors,
\begin{enumerate}
\item \label{f_*} $f_*: \Shv(X;e) \to \Shv(Y;e)$
\item \label{f_!} $f_! : \Shv(X;e) \to \Shv(Y;e)$
\item \label{f^*} $f^*: \Shv(Y;e) \to \Shv(X;e)$
\item \label{f^!} $f^!: \Shv(Y;e) \to \Shv(X;e)$.
\end{enumerate}
There are, additionally,
two functors
\begin{enumerate}
\item[$(5)$] $- \boxtimes - : \Shv(X;e) \times \Shv(Y;e) \to \Shv(X \times Y;e)$
\item[$(6)$] $\sHom(-,-): \Shv(X;e) \times \Shv(X;e) \to \Shv(X;e)$,
\end{enumerate}
which round out the list
of the six functors.

\begin{notn}
If $S \hook X$ is the
inclusion of a subspace $S$
and $\F \in \Shv(X;e)$, we sometimes
use ``${\F|}^*_S$" or ``${\F|}_S$"
to denote the $*$-pullback
$(S \hook X)^*\F$, and ``${\F|}^!_S$"
to denote the $!$-pullback $(S \hook X)^!\F$.
\end{notn}

\subsubsection{$f^*$ and $f_*$}
We briefly describe the
existence of functors \labelcref{f_*}
and \labelcref{f^*}.
The pushforward $\underline{f}_*$ of
presheaves is defined by
precomposition with the 
functor of $\infty$-sites,
$\U(Y) \to \U(X)$, induced by
the map $f$. It admits a left
adjoint, given by left Kan extension
along $f$, which we denote by
$\underline{f}^*$. By 
\cite[Lemma 2.13]{PYY16}, $\underline{f}_*$ 
preserves sheaves
(like one would expect
from the classical definition),
so we define,
\begin{align*}
f_* &:= \underline{f}_* \circ \iota_{X} \\
f^* &:= L_{X} \underline{f}^* \circ \iota_{Y},
\end{align*}
where $\iota_X$ and $\iota_Y$
denotes the inclusions of
sheaves into presheaves
and where $L_X$ is sheafification.
By \cite[Lemma 2.14]{PYY16},
\[f^*: \Shv(Y; e) \rightleftarrows \Shv(X;e): f_*\]
are also adjoints. See \cite[
\S2.4]{PYY16} for more details about
the construction and properties
of these functors.

\subsubsection{$f_!$ and $f^!$}
We briefly sketch the
construction of functors
\labelcref{f_!} and \labelcref{f^!}.
For an individual sheaf
$\F \in \Shv(X; e)$,
$f_!\F$ is defined as the
subsheaf of $f_*\F$ given 
by the formula,
\[f_!\F(U) := \{s \in \F(f^{-1}(U))| \,
	\supp s \to Y \, \textnormal{is proper}\}.\]
Given a map $\F \to \F'
\in \Shv(X;e)$, we obtain
the map $f_*\F \to f_*\F'$
by functoriality; it is 
easy to see that the image
of the restriction of this
map to the subsheaf $f_!\F$
is the subsheaf $f_!\F'$.
From the functoriality
of $f_*$, we therefore
obtain the $\infty$-functor
$f_!$. The functor $f^!$ is
defined as the right adjoint
to $f_!$. Its existence can be
shown either by imitating
the construction in \cite[
\S3.1]{KS90} in an $\infty$-setting
or by inferring its existence
from the Adjoint Functor 
Theorem 
\cite[Corollary 5.5.2.9]{HA} after showing
$f_!$ preserves small colimits.

\subsubsection{Tensor product}
Using the first 
four functors, we 
may also define a tensor
product of sheaves using the functor 
$- \boxtimes -$.
It may be defined using 
the diagonal map 
$X \xrightarrow{\Delta} X \times X$
in one of two inequivalent ways.

\bdef
The $*$\-tensor product is defined as,
\[- \overset{*}{\otimes} - := \Delta^*(- \boxtimes -).\]
The $!$\-tensor product is defined as,
\[- \overset{!}{\otimes} - := \Delta^!(- \boxtimes -).\]
\edefn

Each these functors endows $\Shv(X;e)$
with a symmetric monoidal structure,
which are generally not equivalent to each other.
The practical difference between the two tensor products is
that $f^*$ is symmetric monoidal
for $\overset{*}{\otimes}$, while
$f^!$ is symmetric monoidal for
$\overset{!}{\otimes}$. The $*$\-tensor
product coincides with the canonical
tensor product of module spectra 
(see \cite[\S2.1]{SAG}). For simplicity,
we use the undecorated symbol, $\otimes$,
to refer to the $*$\-tensor product, and
refer to it is simply as the ``tensor product
of sheaves."

\brem
While the exterior product,
$- \boxtimes -$, is the more primordial
functor, the tensor product
is usually taken as one of the Grothendieck
six functors instead.
\erem

\subsubsection{}
The functor $\sHom$ should not be 
confused with the space
of maps in the $\infty$-category
$\Shv(X;e)$, which is denoted
by $\Hom$. In general, $\Hom(-,-)$
is always the space of maps between
two objects in a given category 
while $\sHom(-,-)$ is always a sheaf. 

The relation between
the two for the category $\Shv(X,e)$ is 
as follows: global sections of
$\sHom(\F, \mathcal{G})$ is 
$\Hom(\F, \mathcal{G})$.
$\sHom$ is the internal hom functor
for $\Shv(X;e)$.

\subsubsection{Adjunction triangles}
Given a complementary pair 
of open and closed subsets 
$Z \xhookrightarrow{i} X 
\xhookleftarrow{j} U$, there exists 
the following two canonical 
fiber sequences, which we
call the ``adjunction triangles":
\begin{align}
\label{adj1} i_!i^!\F \to \F \to j_*j^*\F \\
\label{adj2} j_!j^!\F \to \F \to i_*i^*\F.
\end{align}


\subsection{Accessing classical results in the $\infty$-categorical setting}

The canonical reference for the theory of
sheaves on manifolds (or more generally,
locally compact spaces) is the tome 
\cite{KS90}. It and other standard 
references (e.g. \cite{Dim04}) on
the subject were written during a time when
the theory of $\infty$-categories and sheaves
valued therein were either not established
or not adopted in widespread use. Rather,
these references work in the (left, right,
or bounded) derived category of complexes
of sheaves.

Since then, the treatment of $\infty$-categories
in the literature has matured immensely.
Much recent work in the subject of microlocal
sheaf theory and symplectic geometry takes
place in the setting of sheaves valued, not
just in certain stable $\infty$-categories,
but in the category of $\infty$-categories
itself. The toolbox of (microlocal) sheaf theory,
however, continues to be accessed by 
referring to the work of Kashiwara, Schapira.
Usually these references are accompanied
by a short note assuring the reader that
the same definitions, theory, and results
hold in the $\infty$-categorical setting.
This is one such note. 

\Cref{ssec: hyperdescent}
established the ability to check
equivalences in the category $\Shv(X;e)$
by checking on stalks or homotopy sheaves.
\Cref{ssec: six functor formalism} 
established the existence of the six
functors in our setting. We make the 
meta-claim that the results of \cite{KS90}
from this ability and these six functors.
As such, in the rest of this paper, we
make reference to results in \cite{KS90} 
and other classical sources 
with clear conscience and without
further comment, whenever the 
corresponding $\infty$-categorical 
result is obvious.


\subsection{Constructible sheaves}

In this subsection we define
the category of constructible
sheaves, and recall its relevant
properties.

\bdef
Let $M$ be an arbitrary object
of $\Mod^b_e$, viewed as a
sheaf on the point. Then the
constant sheaf on $X$ with
value $M$, denoted by
$\underline{M}_X$, is defined as
\[(X^{\an} \to \ast)^*M.\]
\edefn

\bdef
A sheaf $\F \in \Shv(X;e)$ is
locally constant if there
exists a module $M \in \Mod^b_e$ 
and an open cover of $X^{\an}$,
$\{U_{\a}\}$, such that 
$\F|_{U_{\a}} \simeq 
\underline{M}_{X_{\a}}$.
\edefn

\bdef
\label{definition: constructible sheaves}
A sheaf $\F \in \Shv(X; e)$
is constructible if there
exists a locally finite stratification of $X$
by locally closed complex analytic
subsets $X_{\a}$
having the property that, for each index $\a$,
there exists a perfect $e$\-module
$M_{\a}$ such that $\F|_{X_{\a}}$
is locally constant with value $M_{\a}$.
\edefn

\begin{warning}
When $X$ is a $k$-scheme
locally of finite type, it is typical\footnotemark 
to use a narrower definition of
constructible sheaf, considering 
only stratifications
which are algebraic, i.e.
in which the strata 
are algebraic subsets,
rather than arbitrary analytic
subsets. We do not adopt
this convention here; we instead
refer to the latter as ``algebraically
constructible" sheaves.
We denote the category of
such by $\Shv^{\alg}(X)$.
\end{warning}
\footnotetext{See \cite[Remark 4.1.2(i)]{Dim04}.}

\bdef
We use $\Shv(X)$ to denote 
the full subcategory
of $\Shv^b(X;e)$ spanned by
constructible objects
and call it the bounded derived
category of constructible sheaves
on $X$. 
\edefn

\brem
$\Shv(X)$ is also a stable
$\infty$-category.
\erem



\subsubsection{Constructibility and
the six functors}

The functors
\[- \otimes -, \sHom(-,-),\]
preserve constructibility.
Consider a map of
complex analytic spaces,
$g: (X, \O_{X}) 
\to (Y, \O_{Y})$.
In general, the functors,
\[g^*, g^!\]
each preserve constructibility.
If, additionally, $g$ is
proper on $\supp(\F)$, 
then the sheaves
\[f_*\F, f_!\F\]
are also constructible. In particular
if $g$ is a proper map,
all six functors preserve $\Shv(X)$.
The reference for these
statements is \cite[Theorem 4.1.5]{Dim04}.

\subsubsection{Verdier duality}
Consider the terminal map
from a scheme, $\pt: X \to \ast$.

\bdef
The dualizing sheaf on X is defined as:
\[\omega_X := \pt^!e\]
\edefn

The dualizing sheaf is used to define
the following important functor.

\bdef
The functor of Verdier duality,
$\mathbb{D}: \Shv(X;e) \to \Shv(X;e)^{\op}$,
is defined by the formula,
\[\mathbb{D}\F := \sHom(\F,\omega_X).\]
\edefn

Below are the most immediate
and relevant properties of Verdier
duality.

\begin{enumerate}[(i)]
\item
Verdier duality has the property that
\[\mathbb{D} \circ \mathbb{D} \simeq \id_{\Shv(X;e)}.\]
It is therefore a sort of anti-involution.

\item
It follows immediately from the definition
that Verdier duality also reverses shifts:
\[\mathbb{D}(\F[n]) \simeq (\mathbb{D}\F)[-n],\]
 for all $\F \in \Shv(X;e)$ and $n \in \mathbb{Z}$.

\item{\cite[Theorem 4.1.16]{Dim04}}
Verdier duality preserves constructibility,
i.e. $\mathbb{D}: \Shv(X) \to \Shv(X)^{\op}$.

\item{\cite[Corollary 4.1.17]{Dim04}}
When $f: X \to Y$ is map of schemes 
and $\F \in \Shv^{\alg}(X)$
is an algebraically constructible sheaf,
\[f_!(\mathbb{D}\F) \simeq \mathbb{D}(f_*\F).\]
If $f: X \to Y$ is a map of complex analytic
spaces and $\F \in \Shv(Y)$,
\[f^*(\mathbb{D}\F) \simeq \mathbb{D}(f^!\F).\]
\end{enumerate}

\subsubsection{Perverse sheaves}
 
Recall that a $t$-structure on 
a stable $\infty$-category, 
$\Cee$, is a $t$-structure on its
homotopy category, $\h\Cee$.
$\h\Shv(X)$ is the usual (i.e. $1$\-categorical) 
derived category of constructible sheaves
defined, for example, in \cite[\S4.1]{Dim04},
which we denote by $D_c^b(X; e)$. 
Thus, we define the perverse 
$t$-structure on $\Shv(X)$ to be
the usual perverse $t$-structure on
$D_c^b(X;e)$ for the middle perversity.

\bdef
The category of perverse sheaves
on $X$, denoted by $\Perv(X)$, is
the heart of the perverse $t$-structure
on $\Shv(X)$.
\edefn


\subsection{Nearby and vanishing cycles}

The six-functor formalism
allows us to define the functors
of nearby and vanishing cycles
on $\Shv(X)$.

We choose the exponential map,
$\exp: \widetilde{\G_m} \to \G_m$,
as a model for the universal cover
of the punctured complex plane.
Let $X$ be a complex analytic
space and $f: X \to \Aone$ be a 
complex analytic function. We have
the following diagram:
\begin{equation}
\label{diagram: definition of nearby cycles}
\begin{tikzcd}
	{X_0} & X & {X^*} & {\widetilde{X^*}} \\
	0 & {\Aone} & {\G_m} & {\widetilde{\G_m}}
	\arrow[hook, from=2-1, to=2-2]
	\arrow["i", hook, from=1-1, to=1-2]
	\arrow[hook', from=2-3, to=2-2]
	\arrow["j"', hook', from=1-3, to=1-2]
	\arrow["{\widetilde{\pi}}"', from=1-4, to=1-3]
	\arrow["\exp"', from=2-4, to=2-3]
	\arrow["f"', from=1-2, to=2-2]
	\arrow[from=1-1, to=2-1]
	\arrow[from=1-3, to=2-3]
	\arrow["{\widetilde{f}}"', from=1-4, to=2-4]
	\arrow["\lrcorner"{anchor=center, pos=0.125, rotate=-90}, draw=none, from=1-3, to=2-2]
	\arrow["\lrcorner"{anchor=center, pos=0.125}, draw=none, from=1-1, to=2-2]
	\arrow["\lrcorner"{anchor=center, pos=0.125, rotate=-90}, draw=none, from=1-4, to=2-3]
\end{tikzcd}
\end{equation}

\bdef
\label{definition: nearby cycles}
The nearby cycles 
functor with respect to $f$,
denoted by 
$\psi_f: \Shv(X;e) \to \Shv(X_0;e)$,
is defined by the formula,
\[\psi_f \F := i^*(j \circ 
	\widetilde{\pi})_*(j \circ \widetilde{\pi})^*\F.\]
\edefn

There is an obvious
natural map,
\begin{equation}
\label{equation: nearby cycles comparison map}
i^*\F \xrightarrow{c} \psi_f \F,
\end{equation}
induced by the unit
of the $(j \circ 
\widetilde{\pi})^* \dashv
(j \circ \widetilde{\pi})_*$
adjunction.

\bdef
\label{definition: vanishing cycles}
The vanishing cycles
functor with respect to $f$, denoted
by $\varphi_f: \Shv(X;e) \to 
\Shv(X_0;e)$, is defined to be 
the fiber of the natural map $c$,
\begin{equation}
\label{fiber sequence: nearby/vanishing cycles}
\varphi_f \F \to i^*\F \to \psi_f \F.
\end{equation}
\edefn

\brem
Our definition of the vanishing
cycles functor agrees with the
one defined in \cite[\S8.6]{KS90},
and differs from the ones in
\cite[\S4.2]{Dim04} and \cite[\S3]{M16} 
by a shift of $[-1]$.
\erem

\subsubsection{}
We record a few
important properties
of the nearby and vanishing
cycles functors.
\begin{enumerate}[(i)]
\item{\cite[pg. 103]{Dim04}} $\psi_f$ and $\psi_f$ preserve
constructibility. In the case when
$X$ is a scheme, and $f$ is an algebraic
function, $\psi_f$ and $\psi_f$
preserve algebraic constructibility.
\item{\cite[Theorem 5.2.21]{Dim04}} $\psi_f[-1]$ and $\varphi_f$
are $t$-exact with respect to the
perverse $t$-structure.
\item{\cite[Proposition 4.2.10]{Dim04}}
$\psi_f[-1]$ and $\varphi_f$ commute
with Verdier duality.
\end{enumerate}

\brem
\label{rem: vanishing cycles name convention}
It is conventional to call the perverse
sheaf obtained by applying $\varphi_f$
to the perverse constant sheaf
``the vanishing cycles along $f$,"
and denote it simply by $\varphi_f
:= \varphi_f(\const[\dim X])$.
\erem

\subsubsection{Monodromy automorphism}
\label{sssec: monodromy automorphism}
The nearby cycles functor comes
equipped with a monodromy action
of $\mathbb{Z}$, natural in $\F$
and generated by
the monodromy operator
\[\psi_f \F \xrightarrow{T} \psi_f \F,\]
which is induced by the action of
the fundamental group
of $\G_m$ on the cover
$\widetilde{X^*}$.
If $i^*\F$ is viewed
as having a trivial monodromy
action by $\mathbb{Z}$, $\varphi_f \F$
also obtains a monodromy action
as the fiber of the (trivially
$\mathbb{Z}$\-equivariant) 
map $i^*\F \xrightarrow{c} 
\psi_f \F$ in the category
of constructible sheaves valued
in $e[t, t^{-1}]$.  

\subsubsection{Nearby cycles at a fixed angle}
In practice, there are equivalent
expressions for $\psi_f$ and $\varphi_f$
called the ``nearby and vanishing cycles
at fixed angle" which are 
helpful for computation. In order
to state them, define:
\begin{align*}
H^{\geq}_{\theta} 	&:= f^{-1}(e^{i\theta}\Aone_{\geq 0}) \\
H^{>}_{\theta}		&:= f^{-1}(e^{i\theta}\Aone_{>0}) \\
\ell^{\geq}_{\theta} 	&:= f^{-1}(e^{i\theta}\ell_{\geq 0}) \\
\ell^{>}_{\theta}	&:= f^{-1}(e^{i\theta}\ell_{>0}),
\end{align*}
and denote the inclusions
into $X$ of each by:
\begin{align*}
&I_{\theta}: H^{\geq}_{\theta} \hook X \\
&J_{\theta}: H^{>}_{\theta} \hook X \\
&i_{\theta}: \ell^{\geq}_{\theta} \hook X \\
&j_{\theta}: \ell^{>}_{\theta} \hook X.
\end{align*}

The following equivalences hold\footnotemark:
\footnotetext{\cite[\S3]{M16}}
\begin{align}
\psi_f \F 			&\simeq \psi^{\theta}_f \F	\\	
				&:= {({J_{\theta+ \pi}}_*J_{\theta + \pi}^*\F)|}^*_{f^{-1}(0)} \label{nearby cycles at fixed angle no.1} \\
				&\simeq {({j_{\theta}}_*j_{\theta}^*\F)|}^*_{f^{-1}(0)}  \label{nearby cycles at fixed angle no.2} \\
\varphi_f \F 		&\simeq \varphi^{\theta}_f \F \\ 		
				&:= {(I_{\theta + \pi}^!\F)|}^*_{f^{-1}(0)} \label{vanishing cycles at fixed angle no.1} \\
				&\simeq {(i_{\theta}^!\F)|}^!_{f^{-1}(0)} \label{vanishing cycles at fixed angle no.2}
\end{align}

\brem
One can show that there
are natural isomorphisms
\begin{align*}
T_{\theta}^{\psi}: \psi^{0}_f \F 	&\xrightarrow{\simeq} \psi^{\theta}_f \F \\
T_{\theta}^{\varphi}: \varphi^0_f \F 	&\xrightarrow{\simeq} \varphi^{\theta}_f \F
\end{align*}
which correspond to rotation
counterclockwise over the
origin in $\Aone$. Then
$T_{2\pi}^{\psi}$ and 
$T_{2\pi}^{\varphi}$ correspond
to the natural monodromy
automorphisms on $\psi_f$
and $\varphi_f$ described in
\cref{sssec: monodromy automorphism}.
See \cite[\S2.3]{K21b} for the
explicit construction of these
isomorphisms (as well as
proofs of many of the assertions
in this section).
\erem

\brem
In this paper we 
only utilize expressions
\labelcref{nearby cycles at fixed angle no.1},
\labelcref{nearby cycles at fixed angle no.2},
and \labelcref{vanishing cycles at fixed angle no.1}
for the values $\theta = 0, \pi$.
\erem


\subsection{Monodromic sheaves}

We begin this section by recalling
the definition of a monodromic sheaf,
as first coined by Verdier in \cite[Def. 3.1]{V83b}.

Let $X$ be a $k$-scheme of finite type,
and $C$ a cone over $X$. Let $\theta:
\G_m \times C \to C$ denote the scaling
action on $C$, and for $p \in C$, let $W_p:
\G_m \to C = \theta \circ (p \hook C)$ 
denote the orbit of $p$ under this action.

\bdef
Let $\F \in \Shv(C;e)$ 
be a sheaf on $C$. 
$\F$ is called monodromic
if $W_p^*\F$ is locally constant for all
values of $p$.
\edefn


\begin{notn}
We denote the $\infty$-category of
monodromic sheaves on $C$ by 
$\Shv_{\G_m}(C;e)$.
\end{notn}

In order to enumerate the
properties of monodromic 
sheaves relevant to this work,
we define the notions of equivariant
maps and equivariant sheaves.

\bdef
\label{def: equivariant map}
Suppose that $X$ and $Y$ are
schemes with left $G$\-actions, where
$G$ is a group scheme. A map
$f: X \to Y$ is equivariant if 
the following diagram commutes:
\[\begin{tikzcd}
	{G \times X} & {G \times Y} \\
	X & Y
	\arrow["f", from=2-1, to=2-2]
	\arrow["{\act_Y}"', from=1-2, to=2-2]
	\arrow["{\act_X}"', from=1-1, to=2-1]
	\arrow["{\id_G \times f}", from=1-1, to=1-2].
\end{tikzcd}\]
\edefn

\bdef
\label{def: equivariant sheaf}
Suppose that $X$ is a scheme
with $G$\-action, $\act: G \times X
\to X$. A $G$-equivariant sheaf
is a sheaf on the stack quotient $X/G$,
(i.e.\ an object $\F \in \Shv(X/G;e)$).
Alternatively, the category of $G$-equivariant
sheaves on $X$
is the limit of the usual cosimplicial diagram,
\[\begin{tikzcd}
	{\Shv(X;e)} & {\Shv(X \times_{X/G} X;e)} & {\cdots}
	\arrow[shift right=2, from=1-1, to=1-2]
	\arrow[shift left=2, from=1-1, to=1-2]
	\arrow[from=1-2, to=1-1]
	\arrow[shift right=4, from=1-2, to=1-3]
	\arrow[from=1-2, to=1-3]
	\arrow[shift left=4, from=1-2, to=1-3]
	\arrow[shift left=2, from=1-3, to=1-2]
	\arrow[shift right=2, from=1-3, to=1-2],
\end{tikzcd}\]
which we denote by $\Shv_G^{eq}(X;e)$.
\edefn

We draw attention to the following
well-known property of monodromic
sheaves.

\begin{property}
Denote by $0_{X}: X \hook C$ the
inclusion of $X$ as the base of $C$,
and by $\pi: C \to X$ the structure
map of $C$ as a scheme over $X$.
If $\F \in \Shv_{\G_m}(C;e)$, then
\begin{align*}
\pi_*\F & \simeq 0_{X}^*\F \\
\pi_!\F & \simeq 0_{X}^!\F
\end{align*}
\end{property}


\subsection{Conic sheaves}
If we had taken $k = \R$,
and left $C$ be any locally
compact topological space 
with a $\G_m = \R^+$\-action,
the above definition of 
monodromic sheaf coincides 
with what \cite[\S3.7]{KS90}
calls conic sheaves. 
Monodromic sheaves as 
we've defined them are
obviously conic in this sense, and
herein we will refer to such sheaves 
as either conic or monodromic 
according to whichever property 
we consider in the moment.

\subsubsection{}
We denote the category of
conic sheaves by $\Shv_{\R^+}(-)$.

\subsubsection{}
Conic sheaves on a real, finite
dimensional vector space 
can be seen
as sheaves on 
the much coarser topology
consisting of open conic subsets.

\bdef
Suppose $V$ is a real,
finite dimensional
vector space considered
with the Euclidean topology. 
The collection of open
conic subsets of $V$ is a topology
which we call the conic
topology on $V$. We denote
the corresponding topological
space by $V_{\con}$.
\edefn

There is an obvious continuous map
$\varphi_{\con}: V \to V_{\con}$. Denote
by $\U_{\R^+}(V)$ the collection
of conic open subsets of $V$, considered
as elements of the topology on $V$.

\begin{prop}
\label{prop: conic sheaves are sheaves on conic topology}
The functor $\varphi_{\con}^*: \Shv(V_{\con})
\to \Shv(V)$ induces an equivalence of categories
between $\Shv(V_{\con})$ and $\Shv_{\R^+}(V)$
\end{prop}

\bproof
It suffices to show that the unit and
counit of the $\varphi_{\con}^* \dashv {\varphi_{\con}}_*$
adjunction are equivalences.
Suppose that $\F \in \Shv_{\R^+}(V)$,
and let $x \in V$ be an arbitrary point.
Then
\begin{align*}
(\varphi_{\con}^*{\varphi_{\con}}_*\F)_x 							&\to \F_x \\
({\varphi_{\con}}_*\F)_x 									&\to \F_x \\
\varinjlim_{\{x \in U \subset V_{\con}\}} \Gamma(U; {\varphi_{\con}}_*\F) 	&\to \F_x \\
\varinjlim_{\{x \in U|U \in \U_{\R^+}(V)\}} \Gamma(U; \F)		&\xrightarrow{\simeq} \F_x,
\end{align*}
where the last map is
an equivalence because $\F$ is conic.
Checking that the unit is
an equivalence is done similarly.
\eproof

\subsubsection{}
The terminology and 
results of this section
for conic sheaves on real finite dimensional
vector spaces extend in an
obvious way to sheaves on finite dimensional
real vector bundles which are conic 
with respect to the scaling action
on fibers.
\Cref{prop: conic sheaves 
are sheaves on conic topology} will be
used in section \S3 to specify a map
of conic sheaves by
specifying maps of sections over
each element of a basis
for the conic topology on
a vector bundle.

\brem
\label{rem: clarifying definition of convex}
For sake of clarity later
in the paper, we define
a convex subset of a
real vector bundle to be a subset
whose intersection with any
fiber is convex.
\erem


\subsection{Fourier--Sato transform}
Let $X$ be a complex analytic space,
and $E \to X$ a complex vector bundle
over $X$. Denote by $E^{\vee} \to X$
the (complex) dual vector bundle.
Define the closed, subanalytic subset
$P \subset E \oplus E^{\vee}$ by
\[P = \{(v,w) \in E \oplus E^{\vee}|
	\Real w(v) \geq 0 \}.\]
Consider the following diagram:
\[\begin{tikzcd}
	& P \\
	& {E \oplus E^{\vee}} \\
	{E^{\vee}} && E
	\arrow["{\pr_2}"', from=2-2, to=3-1]
	\arrow["{\pr_1}", from=2-2, to=3-3]
	\arrow["{i_P}", hook, from=1-2, to=2-2].
\end{tikzcd}\]

\bdef{(\cite[\S 5]{V83a}, \cite[\S 3.7]{KS90})}
\label{definition: Fourier--Sato transform}
The Fourier--Sato transform is the functor
$\Four: \Shv_{\G_m}(E) \to \Shv_{\G_m}(E^{\vee})$
defined by the formula,
\begin{align*}
\Four(\F) &:= \pr_{2*}i_{P!}i_P^!\pr_1^*\F \\
		&\simeq \pr_{2*}\Gamma_P(\pr_1^*\F).
\end{align*}
\edefn
The Fourier--Sato transform is an equivalence of categories (\cite[Theorem 3.7.9]{KS90}).

\subsubsection{}
The above definition contains 
the hidden claim that $\Four(\F)$
is monodromic, but this is easily
seen from the formula as follows.

\bproof[Proof that $\Four(\F)$ is monodromic]
Denote by $\overline{E}^{\vee}$,
the vector bundle whose fibers are
the complex conjugates of the fibers
of $E^{\vee}$. Note that $\overline{E}^{\vee}$
and $E^{\vee}$ have the same underlying
topological space, $|E^{\vee}|$. The orbit of
$p \in |E^{\vee}|$ under 
the complex conjugate action of $\G_m$
is the same as 
under the original action of $\G_m$. 
Thus, $\Shv_{\G_m}(E^{\vee})$ and
$\Shv_{\G_m}(\overline{E}^{\vee})$ are
identical subcategories of $\Shv(|E^{\vee}|)$.

Consider the bundle $E \oplus 
\overline{E}^{\vee}$ with the diagonal
$\G_m$\-action (i.e. its action as a vector
bundle over $X$). The subset $P$ is invariant
under this action, meaning we obtain the 
Cartesian square,
\[\begin{tikzcd}
	{P \times \G_m} & P \\
	{E \oplus \overline{E}^{\vee} \times \G_m} & {E \oplus \overline{E}^{\vee}}
	\arrow["{i_P}", hook, from=1-2, to=2-2]
	\arrow["{i_P \times \id_{\G_m}}"', hook, from=1-1, to=2-1]
	\arrow["{\act^P_{\G_m} := \act^{E \oplus \overline{E}^{\vee}}_{\G_m}|_{P}}", from=1-1, to=1-2]
	\arrow["{\act^{E \oplus \overline{E}^{\vee}}_{\G_m}}", from=2-1, to=2-2],
\end{tikzcd}\]
demonstrating that $i_P$ 
is $\G_m$\-equivariant.
At the same time, $\pr_1$ and $\pr_2$ are
clearly equivariant as maps of vector bundles,
so we obtain that $\Four(\F) := 
\pr_{2*}i_{P!}i_P^!\pr_1^*\F$
is monodromic.
\eproof

\subsubsection{Naming conventions}
The definition of $\Four$ first
appears in \cite[\S 5]{V83a} under the
name ``geometric Fourier transform"---as 
opposed to the Fourier transform of 
$\mathscr{D}$-modules, which 
Verdier calls the ``analytic Fourier
transform." In keeping with modern convention
(e.g. \cite[\S2.4]{K21b}), 
however, we refer the functor
$\Four$ as the Fourier--Sato transform, a name
which first appears in \cite[\S 3.7]{KS90}\footnotemark.
\footnotetext{The functor which appears
in \cite[\S 3.7]{KS90} is technically different
from that defined above, being
defined for conic sheaves on
real vector bundles, but its usage
later in that work indicates that the authors
were aware of the obvious generalization
to complex vector bundles.}

\subsubsection{}
There is an alternate definition of
the Fourier--Sato transform using
the subset 
\[P' := \{(v,w) \in E \oplus
	E^{\vee}| \Real w(v) \leq 0 \},\]
as well as an inverse functor for the 
Fourier--Sato transform exist, but
we will not need these in this paper.
A wonderful, detailed account of
these functors and their various
compatibilities is found in \cite[\S2.4]{K21b},
from which we are borrowing
much of the notation in this
section.

\subsubsection{}
We will need the following two
properties of the 
Fourier--Sato transform.

\blem
\label{lemma: base change for FS transform}
Suppose that $f: X' \to X$ is a map
of complex analytic spaces, and let
$E \to X$ be a complex vector bundle.
Denoting $E' := E \times_{X} X'$,
we obtain the induced maps $f_{E}:
E' \to E$ and $f_{E^{\vee}}: E'^{\vee}
\to E^{\vee}$. The following natural
transformation is an equivalence:
\[f_{E^{\vee}}^* \Four_{E} \xrightarrow{\simeq}
	\Four_{E'}f_{E}^*.\]
\elem

\bproof
See \cite[Proposition 3.7.13]{KS90}.
\eproof

\blem
\label{lemma: FST and perversity}
Suppose that $E \to X$ is a rank $n$ 
(i.e. $\dim_{\mathbb{C}}=n$) complex
vector bundle. Then, 
	\[\Four[n]: \Perv(E) \to \Perv(E^{\vee}).\]
\elem

\bproof
See \cite[Proposition 10.3.18]{KS90}.
\eproof


\section{Deformation to the normal bundle}
\label{sec: deformation to the normal bundle}

In this section, we recall the definition of
a quasi-smooth morphism and the construction
in \cite{KR19} of deformation to the normal 
bundle along quasi-smooth closed immersions.
We then prove several general facts about their
construction needed in order to set up the theory
of quasi-smooth specialization and microlocalization.


\subsection{Quasi-smooth morphisms}
\subsubsection{}
In the body of the paper, we will need the notion
of a quasi-smooth closed immersion:
the correct notion of a regular
embedding in the derived setting. 
Such morphisms $\Z \to \X$ 
are Zariski locally
modeled on the derived zero locus
of a map from $\X$ to affine 
space of some finite rank. They are 
highly structured, having, 
for example, a well-behaved
notion of normal bundle.

\subsubsection{}
\bdef
A map $f: \bm{X} \to \bm{Y}$ of derived schemes is quasi-smooth if the 
relative cotangent complex $\mathbb{L}_f$ is perfect of Tor-amplitude $[-1, 0]$.
$\mathbb{L}_f$ and $\mathbb{L}_{\bm{X}/\bm{Y}}$ will be used interchangeably 
when the map $f$ is clear from context.
\edefn

As a special case of the above definition, we say that a derived 
scheme $\bm{X}$ is quasi-smooth if $\mathbb{L}_{\bm{X}/k}$ 
is perfect of Tor-amplitude $[-1, 0]$.

\blem
\label{lemma: quasi-smooth}
Assume $X$ is smooth, and suppose 
$i: \bm{Z} \hook X$ is a closed immersion
of derived schemes. Then $i$ is 
quasi-smooth iff $\bm{Z}$ is quasi-smooth.
\elem

\bproof
Since $X$ is smooth, $i^*\mathbb{L}_X$ is locally free. 
$\mathbb{L}_{\bm{Z}}$ is perfect by assumption. 
Consider the standard sequence on 
cotangent complexes, $i^*\mathbb{L}_X \xrightarrow{i} 
\mathbb{L}_{\bm{Z}} \to \mathbb{L}_i$.
If any two of the three terms in a distinguished 
triangle in $\QCoh(\bm{Z})$
are perfect, the third is as well \cite[Lemma 066R]{Stacks}, so 
if either of $\mathbb{L}_i$ 
or $\mathbb{L}_{\bm{Z}}$ is perfect, the other is as well. 
Meanwhile, $i^*\mathbb{L}_X$ is flat, so has Tor amplitude $[0]$. 
By \cite[Lemma 0655]{Stacks}, this implies that if 
either $\mathbb{L}_i$ or $\mathbb{L}_{\bm{Z}}$ 
has Tor-amplitude in $[-1,0]$, the other does as well.
\eproof

The following lemma gives a local model for 
quasi-smooth closed immersions.

\blem 
\label{lemma: local model}
Suppose that $i: \Z \hook \X$ is a quasi-smooth closed immersion 
of derived schemes.
Then Zariski-locally on $\X$, $i$ fits into a Cartesian diagram
	\[\begin{tikzcd}
		{\bm{Z}} & {\bm{X}} \\
		{*} & {\mathbb{A}^n}
		\arrow["0", from=2-1, to=2-2]
		\arrow[from=1-2, to=2-2]
		\arrow["i", hook, from=1-1, to=1-2]
		\arrow[from=1-1, to=2-1]
		\arrow["\lrcorner"{anchor=center, pos=0.125}, draw=none, from=1-1, to=2-2]
	\end{tikzcd}.\]
Moreover, $\mathbb{L}_i[-1]$ is then locally free
of finite rank. Denote this sheaf by $\N_{\Z/\X}$.
\elem

\bproof 
See \cite[Proposition 2.3.8]{KR19}, which in fact
shows the converse is also true.
\eproof



\bdef{(\cite[\S 2.3.11]{KR19})}
\label{definition: virtual codimension}
The virtual codimension of $i$,
defined Zariski-locally on $\Z$, is the
rank of the $\O_{\Z}$\-module $\N_{\Z/\X}$.
We denote it by $\vircodim(\Z, \X)$, or
$\vircodim \, \Z$ when there is no ambiguity.
\edefn

\bdef
\label{definition: normal bundle}
The normal bundle to $\Z$ inside $\X$
is the scheme, 
\[\bm{N}_{\Z/\X} := 
	\Spec_{\Z}(\Sym_{\O_{\Z}}^\bullet(\N_{\Z/\X})).\]
Similarly, the conormal bundle to $\Z$ inside $\X$,
denoted $\bm{N}^*_{\Z/\X}$, is the total space 
of the sheaf dual to $\N_{\Z/\X}$. 
\edefn

\subsubsection{}
Obviously the rank of $\bm{N}_{\Z/\X}$ is $\vircodim \, \Z$.


\subsection{Definition of deformation to the normal bundle}
\label{ssec: definition of deformation to the normal bundle}
The deformation to the normal
bundle described in \cite{KR19} 
is obtained by removing a divisor
from a blow-up, just as in the
 classical construction---only
the precise notions of ``blow-up"
and ``divisor" differ.
The bulk of \cite{KR19}
is spent developing the theory
of blow-ups along quasi-smooth
closed immersions of derived stacks.
Roughly, quasi-smooth blow-ups
are defined by a universal
property that generalizes the
classical one characterizing the
blow-up of $X$ along $Z$ as the
terminal object in schemes over
$X$ whose fiber over $Z$ is a Cartier
divisor. The exact definition is
given in \cite[\S4]{KR19}.

Once the authors of \textit{op. cit.} 
obtain a good theory
of quasi-smooth blow-ups,
they define the deformation to the
normal bundle of $\Z \hook \X$,
which we denote $\bm{D}_{\Z/\X}$,
as the open substack of 
$\Bl_{\Z \times \{0\}}(\X \times \Aone)$
given by removing the immersed copy of
$\Bl_{\Z}\X$ inside it. We recall the
theorem of \cite{KR19}, establishing
the existence of deformation
to the normal bundle.

\bthm[{\cite[Theorem 4.1.13]{KR19}}]
\label{theorem: deformation}
Let $i: \Z \hook \X$ be a quasi-smooth 
closed immersion of derived
stacks with normal bundle 
$\bm{N}_{\Z/\X} = \Spec_{\Z}(\Sym_{\O_{\Z}}
(\N_{\Z/\X}))$. Then there exists a
canonical factorization of 
$i \times \id$, the deformation to the 
normal bundle:
\[i \times \id: \Z \times \Aone \overset{j}{\longhook} \bm{D}_{\Z/\X} \overset{\pi}{\longrightarrow} \X \times \Aone\]
satisfying the following properties
\begin{enumerate}[(i)]
	\item The factorization is stable under arbitrary derived base change along $\X$.
	\item $j$ is a quasi-smooth closed immersion.
	\item $\pi$ is quasi-smooth and quasi-projective.
	\item Restricting to $\mathbb{G}_m = \Aone \setminus \{0\}$ we obtain
		\[i \times \id: \Z \times \mathbb{G}_m
		\overset{j_\mathrm{gen}}{\longhook} \X \times \G_m
		\overset{\pi_\mathrm{gen}}{\longrightarrow} \X \times \G_m\]
		where $j_{\mathrm{gen}} = i \times \id$ 
		and $\pi_\mathrm{gen} = \id$.
	\item Restricting to $\{0\}$ we obtain
		\[i: \Z \overset{j_0}{\longhook} \bm{N}_{\Z/\X} \overset{\pi_0}{\longrightarrow} \X\]
		where $j_0$ is the zero-section and 
		$\pi_0$ the composition of the projection 
		$\bm{N}_{\Z/\X} \to \Z$
		and $i: \Z \hook \X$.
\end{enumerate}
\ethm

\subsubsection{}
An immediate consequence of
\cref{theorem: deformation}
is the existence of the following
commutative diagram of derived 
stacks, almost identical to diagram
\labelcref{diagram: deformation 1}
above: 
\begin{equation}
\begin{tikzcd}
\label{diagram: deformation 2}
	{\bm{N}_{\Z/\X}} & {\bm{D}_{\Z/\X}} & {\X \times \G_m} & {\X} \\
	{\{0\}} & \Aone & {\G_m}
	\arrow[hook, from=2-1, to=2-2]
	\arrow[from=2-3, to=2-2]
	\arrow[from=1-1, to=2-1]
	\arrow[from=1-3, to=2-3]
	\arrow["{\bm{t}}", from=1-2, to=2-2]
	\arrow["{\bm{pr}_1}", from=1-3, to=1-4]
	\arrow["{\bm{j}_{\neq 0}}"', from=1-3, to=1-2]
	\arrow["{\bm{s}}", hook, from=1-1, to=1-2].
\end{tikzcd}
\end{equation}

\subsubsection{}
If $i: \Z \hook \X$ is a quasi-smooth closed immersion 
of derived schemes, $\bm{D}_{\Z/\X}$
is (representable by) a derived scheme by
\cite[Theorem 4.1.5(i)]{KR19}.


\subsection{Formal properties of deformation to the normal bundle}

\subsubsection{}
Familiarity with the notation, 
terminology, and concepts of 
\cite{KR19} is required to understand
this subsection.

\subsubsection{}
The following lemma says,
roughly speaking, that a map of derived
schemes induces a map of their
blow-ups along compatible 
quasi-smooth, closed immersions.

\blem
\label{lemma: blowup map}
Suppose that $\bm{f}: \X' \to \X$ is a
map of derived schemes, and that 
$\Z \hook \X$, $\Z' \hook \X'$ 
are quasi-smooth
closed immersions such that
$\Z' \subset \X' \times_{\X} \Z$;
$Z' \simeq X' \times_{X} Z$; and the
canonical morphism $\bm{f}_{\Z}^*
\N_{\Z/\X} \to \N_{\Z'/\X'}$ is surjective
on $\pi_0$.
Then there exists a natural
map $\Bl_{\bm{f}}: \Bl_{\Z'}\X'
\to \Bl_{\Z}\X$ such that the following
diagram commutes:
\begin{equation}
\label{diagram: blowup map 1}
\begin{tikzcd}
	{\P_{\Z'}(\N_{\Z'/\X'})} & {\Bl_{\Z'}\X'} \\
	{\P_{\Z}(\N_{\Z/\X})} & {\Bl_{\Z}\X} \\
	\Z & \X
	\arrow["{\P(\bm{N}_{\bm{f}})}", from=1-1, to=2-1]
	\arrow["{\Bl_{\bm{f}}}"', from=1-2, to=2-2]
	\arrow[hook, from=1-1, to=1-2]
	\arrow[hook, from=2-1, to=2-2]
	\arrow["{\bm{\pi}_{\X}}"', from=2-2, to=3-2]
	\arrow["{\bm{\tau}}", from=2-1, to=3-1]
	\arrow[hook, from=3-1, to=3-2]
	\arrow["{\bm{f}|_{\Z'} \circ \bm{\tau}'}"'{pos=0.47}, curve={height=40pt}, from=1-1, to=3-1]
	\arrow["{\bm{f}\circ \bm{\pi}_{\X'} }"{pos=0.5}, curve={height=-20pt}, from=1-2, to=3-2].
\end{tikzcd}
\end{equation}
If, in fact, $\Z' \simeq \X' \times_{\X} \Z$, then
we furthermore have a commutative diagram:
\begin{equation}
\label{diagram: blowup map 2}
\begin{tikzcd}
	{\Bl_{\Z'}\X'} & {\X' \setminus \Z'} \\
	{\Bl_{\Z}\X} & {\X \setminus \Z} \\
	\X & {\X \setminus \Z}
	\arrow["{\bm{f}|_{\X' \setminus \Z'}}"{pos=0.6}, from=1-2, to=2-2]
	\arrow[from=2-2, to=2-1]
	\arrow[from=1-2, to=1-1]
	\arrow[from=3-2, to=3-1]
	\arrow["{\bm{\pi}_{\X}}", from=2-2, to=3-2]
	\arrow["{\bm{f}\circ \bm{\pi}_{\X'} }"'{pos=0.5}, curve={height=20pt}, from=1-1, to=3-1]
	\arrow["{\bm{f}|_{\X \setminus \Z'} \circ \bm{\pi}_{\X'}}"{pos=0.5}, curve={height=-50pt}, from=1-2, to=3-2]
	\arrow["{\bm{\pi}_{\X}}", from=2-1, to=3-1]
	\arrow["{\Bl_{\bm{f}}}", from=1-1, to=2-1].
\end{tikzcd}
\end{equation}
Moreover, all the squares in diagrams
\labelcref{diagram: blowup map 1}
and \labelcref{diagram: blowup map 2}
are Cartesian on underlying
classical schemes.
\elem

\bproof
Given a derived scheme $\bm{S} \in \DSch_{/\X'}$,
the mapping space $\Hom_{\DSch_{/\X'}}
(\bm{S}, \Bl_{\Z'}\X'$ is the space
of virtual Cartier divisors on $\bm{S}$ 
lying over $(\X', \Z')$. Recall that a point
in this space consists of a commutative
square,
\[\begin{tikzcd}
	{\bm{D}} & {\bm{S}} \\
	{\Z'} & {\X'}
	\arrow["{\bm{h}'}", from=1-2, to=2-2]
	\arrow["{\bm{i}'}", hook, from=2-1, to=2-2]
	\arrow["{\bm{i}_{\bm{D}}}", hook, from=1-1, to=1-2]
	\arrow["{\bm{g}'}", from=1-1, to=2-1],
\end{tikzcd}\]
satisfying the following properties:
\begin{enumerate}[(a)]
\item $\bm{i}_{\bm{D}}: \bm{D} \hook \bm{S}$
	exhibits $\bm{D}$ as a virtual Cartier divisor on $\bm{S}$
	(i.e. a quasi-smooth closed immersion whose normal
	bundle has rank $1$).
\item The underlying square of classical schemes is Cartesian.
\item The canonical base change morphism
	\[\bm{g}'^*\N_{\Z'/\X'} \to \N_{\bm{D}/\bm{S}} \]
	is surjective on $\pi_0$.
\end{enumerate}
Denote the space of such by
$\VDiv(\bm{S}/(\X', \Z'))$.

Given a point $\bm{i}_{\bm{D}} \in 
\VDiv(\bm{S}/(\X', \Z'))$,
we show that it naturally obtains a virtual Cartier
divisor on $\bm{S}$ living over $(\Z, \X)$. 
Indeed, by hypothesis, $\bm{i}_{\bm{D}}$
fits into the commutative diagram,
\[\begin{tikzcd}
	{\bm{D}} & {\bm{S}} \\
	{\Z'} & {\X'} \\
	\Z & \X
	\arrow["{\bm{g}'}", from=1-1, to=2-1]
	\arrow["{\bm{h}'}", from=1-2, to=2-2]
	\arrow["{\bm{f}_{\Z}}", from=2-1, to=3-1]
	\arrow["{\bm{f}}", from=2-2, to=3-2]
	\arrow["{\bm{i}}", hook, from=3-1, to=3-2]
	\arrow["{\bm{i}'}", hook, from=2-1, to=2-2]
	\arrow["{\bm{i}_{\bm{D}}}", from=1-1, to=1-2].
\end{tikzcd}\]
Denote the composition of the outer left and right
vertical arrows by $\bm{g}$ and $\bm{h}$ respectively.
The claim is that the outer commutative
square of this diagram satisfies the conditions
of a virtual Cartier divisor on $\bm{S}$ over
$(\X, \Z)$.

Indeed, $\bm{i}$ was taken to be a virtual
Cartier divisor on $\bm{S}$, so condition (a)
is satisfied. Condition (b) is satisfied because
of the hypothesis that $Z' \simeq Z \times_{X} X'$.
Condition (c) is satisfied because the base change
morphism $\bm{g}^*\N_{\Z/\X} \to 
\N_{\bm{D}/\bm{S}}$ factors as
\[\bm{g}'^*\bm{f}_{\Z}^*\N_{\Z/\X}
	\to \bm{g}'^*\N_{\Z'/\X'} \to \N_{\bm{D}/\bm{S}}\]
where the first morphism is induced by
$\bm{f}_{\Z}^*\N_{\Z/\X} \to \N_{\Z'/\X'}$,
which is surjective by assumption.
Surjectivity on $\pi_0$ then follows from
that of $\bm{g}'^*\N_{\Z'/\X'} 
\to \N_{\bm{D}/\bm{S}}$ and the
fact that $\bm{g}'^*$ preserves
surjectivity on $\pi_0$ (this is just
the classical statement that 
pull-back of quasi-coherent sheaves
is right exact).

The assignment $\VDiv(\bm{S}/(\X', \Z'))
\to \VDiv(\bm{S}/(\X, \Z))$ described above
is obviously functorial in $\bm{S}$, so by
Yoneda we obtain a map 
$\Bl_{\Z'}\X' \to \Bl_{\Z}\X \times_{\X} \X'
\in \DSch_{/\X'}$. Composition
with the projection onto $\Bl_{\Z}\X$
yields the desired map $\Bl_{\bm{f}}$
of derived schemes.

A map $\bm{S} \to \bm{S}' \in \DSch_{/\X}$
induces a map on virtual Cartier divisors
over $(\X, \Z)$ by precomposition with
points in $\Hom_{\DSch_{/\X}}
(\bm{S'}, \Bl_{\Z}\X)$. Explicitly,
the image of an element 
$\bm{i}_{\bm{D}'} \in
\VDiv(\bm{S}'/(\X, \Z))$ under precomposition 
by $\bm{S} \to \bm{S}'$ is given by an
element $\bm{i}_{\bm{D}} \in \VDiv(\bm{S}/(\X,\Z))$
whose diagram 
\[\begin{tikzcd}
	{\bm{D}} & {\bm{S}} \\
	\Z & \X
	\arrow[from=1-1, to=2-1]
	\arrow["{\bm{i}_{\bm{D}}}", hook, from=1-1, to=1-2]
	\arrow[from=1-2, to=2-2]
	\arrow[hook, from=2-1, to=2-2]
\end{tikzcd}\]
factors as
\[\begin{tikzcd}
	{\bm{D}} & {\bm{S}} \\
	{\bm{D}'} & {\bm{S}'} \\
	\Z & \X
	\arrow["{\bm{i}_{\bm{D}}}", hook, from=1-1, to=1-2]
	\arrow["{\bm{i}_{\bm{D}'}}", hook, from=2-1, to=2-2]
	\arrow[hook, from=3-1, to=3-2]
	\arrow[from=2-2, to=3-2]
	\arrow[from=2-1, to=3-1]
	\arrow[from=1-2, to=2-2]
	\arrow[from=1-1, to=2-1].
\end{tikzcd}\] 

Recall that the closed immersion
$\P_{\Z}(\N_{\Z/\X}) \hook \Bl_{\Z}\X$
is the virtual Cartier divisor lying over
$(\X, \Z)$ given by the identity 
$\Bl_{\Z}\X \overset{\id}{\to} \Bl_{\Z}\X$
(see proof of \cite[Theorem 4.1.4(iv)]{KR19}).
Note that the map $\Bl_{\bm{f}}:
\Bl_{\Z'}\X' \to \Bl_{\Z}\X$ factors trivially
through the identity:
\[\Bl_{\Z'}\X' \to \Bl_{\Z}\X \overset{\id}{\to}
	\Bl_{\Z}\X.\]
As such, we see that the virtual Cartier
divisor given by $\Bl_{\bm{f}}$,
$\bm{D}_{\Bl_{\bm{f}}} \to \Bl_{\Z'}\X'$, 
factors as
\[\begin{tikzcd}
	{\bm{D}_{\Bl_{\bm{f}}}} & {\Bl_{\Z'}\X'} \\
	{\P_{\Z}(\N_{\Z/\X})} & {\Bl_{\Z}\X} \\
	\Z & \X
	\arrow[from=2-1, to=3-1]
	\arrow[hook, from=2-1, to=2-2]
	\arrow[from=2-2, to=3-2]
	\arrow[hook, from=3-1, to=3-2]
	\arrow["{\Bl_{\bm{f}}}", from=1-2, to=2-2]
	\arrow[from=1-1, to=2-1]
	\arrow[hook, from=1-1, to=1-2].
\end{tikzcd}\]
It remains to be shown that 
$\bm{D}_{\Bl_{\bm{f}}} \hook \Bl_{\Z'}\X'$
is the divisor $\P_{\Z'}(\N_{\Z'/\X'} \hook \Bl_{\Z'/\X'}$.
This follows easily by noting that, as illustrated here,
\[\begin{tikzcd}[row sep=0]
	{\Hom_{\DSch_{/\X'}}(\Bl_{\Z'}\X', \Bl_{\Z'}\X')} & {\Hom_{\DSch_{/\X}}(\Bl_{\Z'}\X', \Bl_{\Z}\X)} \\ 
	\rotatebox{90}{$\in$} & \rotatebox{90}{$\in$} \\
	{\id} & {\Bl_{\bm{f}}}
	\arrow["{\Bl_{\bm{f}} \circ -}", from=1-1, to=1-2]
	\arrow[maps to, from=3-1, to=3-2],		
\end{tikzcd}\]
the divisor $\Bl_{\bm{f}}$ is the image
of the divisor $\id_{\Bl_{\Z'}\X'}$ 
under the map which
views elements of $\VDiv(\Bl_{\Z'}\X'/(\X', \Z'))$
as elements of $\VDiv(\Bl_{\Z'}\X'/(\X, \Z))$.

Let $\bm{\pi}: \P_{\Z'}(\N_{\Z'/\X'})
\to \Z$ and 
$\bm{\tau}': \P_{\Z'}(\N_{\Z'/\X'}) 
\to \Z'$ denote the structure maps of
$\P_{\Z'}(\N_{\Z'/\X'})$ as a divisor
lying over $(\X, \Z)$ and
$(\X', \Z')$, respectively.
Unraveling the definitions
used in the argument above,
observe that the natural map 
$\P_{\Z'}(\N_{\Z'/\X'}) \to 
\P_{\Z}(\N_{\Z/\X})$, which we
label by ``$?$", is classified
by the data of the canonical, 
surjective (on $\pi_0$) morphism
\[\bm{\pi}^*\N_{\Z/\X} \to
	\N_{\P_{\Z'}(\N_{\Z'/\X'})/\Bl_{\Z'}\X'}\]
associated to $\P_{\Z'}(\N_{\Z'/\X'}) \hook
\Bl_{\Z'}\X'$ as a virtual divisor over $(\X, \Z)$.
This morphism in turn factors as
\[(\bm{\pi}^* \simeq) \bm{\tau}'^*\bm{f}_{\Z}^*\N_{\Z/\X} \to
	\bm{\tau}'^*\N_{\Z'/\X'} \twoheadrightarrow
		\N_{\P_{\Z'}(\N_{\Z'/\X'})/\Bl_{\Z'}\X'}.\]
This factorization then classifies a commutative
triangle in $\DSch_{/\X}$,
\[\begin{tikzcd}
	{\P_{\Z'}(\N_{\Z'/\X'})} & {\P_{\Z'}(\N_{\Z'/\X'})} \\
	& {\P_{\Z}(\N_{\Z/\X}) }
	\arrow["{\P(\bm{N}_{\bm{f}})}", from=1-2, to=2-2]
	\arrow["\id", from=1-1, to=1-2]
	\arrow["{?}"', from=1-1, to=2-2],
\end{tikzcd}\]
by the moduli-theoretic description of
$\P_{\Z}(\N_{\Z/\X})$. Thus 
$? = \P(\N_{\bm{f}})$. This 
yields diagram
\labelcref{diagram: blowup map 1}.

To see the existence of diagram
\cref{diagram: blowup map 2}
observe that
the map $?: \Bl_{\Z'}\X' \setminus 
\P_{\Z'}(\N_{\Z'/\X'}) \to \Bl_{\Z}\X
\setminus \P_{\Z}(\N_{\Z/\X})$
classifies, in part, the following
commutative triangle,
\[\begin{tikzcd}
	{\Bl_{\Z'}\X' \setminus \P_{\Z'}(\N_{\Z'/\X'})} & {\Bl_{\Z}\X \setminus \P_{\Z}(\N_{\Z/\X})} \\
	{\X \setminus \Z}
	\arrow["\simeq", from=1-2, to=2-1]
	\arrow["{\bm{f}|_{\X' \setminus \Z'} \circ \bm{\pi}_{\X'}}"', from=1-1, to=2-1]
	\arrow["{?}", from=1-1, to=1-2].
\end{tikzcd}\]
Identifying $\Bl_{\Z'}\X' \setminus 
\P_{\Z'}(\N_{\Z'/\X'})$ with 
$\X' \setminus \Z'$ and $\Bl_{\Z}\X
\setminus \P_{\Z}(\N_{\Z/\X})$ with
$\X \setminus \Z$ using 
\cite[Theorem 4.1.5(v)]{KR19}, 
we obtain that 
$? = \bm{f}|_{\X' \setminus \Z'}$
as desired.

The claim that all underlying
squares of classical schemes in
diagrams \labelcref{diagram: blowup
map 1} and \labelcref{diagram: blowup
map 2} are Cartesian follows trivially
from the definitions in \cite{KR19}
and our assumptions.
\eproof

\brem
Note that, in contrast with
the theory of classical blow-ups,
unless $\Z \hook \X$ 
itself has virtual codimension $1$,
the square, 
\[\begin{tikzcd}
	{\P_{\Z}(\N_{\Z/\X})} & {\Bl_{\Z}\X} \\
	\Z & \X
	\arrow[from=1-2, to=2-2]
	\arrow[hook, from=2-1, to=2-2]
	\arrow[hook, from=1-1, to=1-2]
	\arrow[from=1-1, to=2-1],
\end{tikzcd}\]
will never be Cartesian as a diagram
of derived schemes, so in general 
the squares in
diagram \labelcref{diagram: blowup map 1}
are only Cartesian on the level of
underlying classical schemes.
\erem

An immediate corollary 
of \cref{lemma: blowup map}
is the following.

\begin{cor}
\label{cor: deformation map}
Suppose that $\bm{f}: \X' \to \X$ is a
map of derived schemes, and 
that $\Z \hook \X$, $\Z' \hook \X'$ 
are quasi-smooth
closed immersions such that
$\Z' \simeq \X' \times_{\X} \Z$.
Then there exists a natural
map $\D_{\f}: \bm{D}_{\Z'/\X'} \to
\bm{D}_{\Z/\X}$
making the following diagram
commute:
\begin{equation}
\begin{tikzcd}
\label{diagram: map of deformations}
	{\bm{N}_{\Z'/\X'}} & {\bm{D}_{\Z'/\X'}} & {\X' \times \G_m} & {\X'} \\
	{\bm{N}_{\Z/\X}} & {\bm{D}_{\Z/\X}} & {\X \times \G_m} & \X \\
	{\{0\}} & \Aone & {\G_m}
	\arrow["{\bm{s}'}", hook, from=1-1, to=1-2]
	\arrow["{\bm{s}}", hook, from=2-1, to=2-2]
	\arrow["{\bm{N}_{\bm{f}}}", from=1-1, to=2-1]
	\arrow["{\D_{\f}}"', from=1-2, to=2-2]
	\arrow["{\bm{f} \times \id_{\G_m}}", from=1-3, to=2-3]
	\arrow["{\bm{f}}", from=1-4, to=2-4]
	\arrow["{\bm{pr}_2}", from=2-3, to=3-3]
	\arrow["{\bm{t}}"', from=2-2, to=3-2]
	\arrow[from=2-1, to=3-1]
	\arrow[hook, from=3-1, to=3-2]
	\arrow[from=3-3, to=3-2]
	\arrow["{\bm{j}_{\neq 0}}"', from=2-3, to=2-2]
	\arrow["{\bm{j}_{\neq 0}'}"', from=1-3, to=1-2]
	\arrow["{\bm{pr}_1'}", from=1-3, to=1-4]
	\arrow["{\bm{pr}_1}", from=2-3, to=2-4]
	\arrow["{\bm{t}'}"{pos=0.3}, curve={height=-20pt}, from=1-2, to=3-2].
\end{tikzcd}
\end{equation}
\end{cor}

\bproof
Using \cref{lemma: blowup map} and
\cite[Theorem 4.1.5(iii)]{KR19} for the
inclusions $\X = \X \times \{0\} \hook
\X \times \Aone$, $\X' = \X' \times \{0\}
\hook \X' \times \Aone$ we obtain the
following map of blow-ups,
\begin{equation}
\label{diagram: deformation map 1}
\begin{tikzcd}
	{\Bl_{\Z' \times \{0\}}(\X' \times \{0\})} & {\Bl_{\Z' \times \{0\}}(\X' \times \Aone)} \\
	{\Bl_{\Z \times \{0\}}(\X \times \{0\})} & {\Bl_{\Z \times \{0\}}(\X \times \Aone)}
	\arrow["{\Bl_{\bm{f} \times \id_{\Aone}}}", from=1-2, to=2-2]
	\arrow["{\Bl_{\bm{f}}}", from=1-1, to=2-1]
	\arrow["{\Bl_{\bm{i}'}}", hook, from=1-1, to=1-2]
	\arrow["{\Bl_{\bm{i}}}", hook, from=2-1, to=2-2],
\end{tikzcd}
\end{equation}
whose commutativity follows easily
from unwinding the definition of each map.
Roughly speaking, each map in the
above diagram is obtained by viewing
virtual Cartier divisors lying over one pair
as Cartier divisors lying over another pair.
Clearly, the two ways specified in
the above diagram of viewing a
virtual Cartier divisor on $\bm{S}$ lying
over $(\X', \Z')$ as one lying over
$(\X \times \Aone, \Z \times \{0\})$
are compatible.

\begin{claim}
The diagram
\labelcref{diagram: deformation map 1}
is Cartesian. 
\end{claim}

\bproof
Indeed, suppose we are
given a commutative diagram,
\[\begin{tikzcd}
	{\bm{S}} & {\Bl_{\Z' \times \{0\}}(\X' \times \Aone)} \\
	{\Bl_{\Z \times \{0\}}(\X \times \{0\})} & {\Bl_{\Z \times \{0\}}(\X \times \Aone)}
	\arrow["\b"', from=1-1, to=2-1]
	\arrow["{\Bl_{\bm{i}}}", from=2-1, to=2-2]
	\arrow["{\Bl_{\bm{f}}}", from=1-2, to=2-2]
	\arrow["\a", from=1-1, to=1-2].
\end{tikzcd}\]
This classifies the data of an
element $\bm{i}_{\bm{D}_{\Z}}:
\bm{D}_{\Z} \to \bm{S} \in
\VDiv(\bm{S}/(\X, \Z))$, an element
$\bm{i}_{\bm{D}_{\Z' \times \{0\}}}:
\bm{D}_{\Z \times \{0\}} \to \bm{S} \in
\VDiv(\bm{S}/(\X \times \Aone, \Z' \times \{0\}))$,
and an equivalence $\D_{\Z} \simeq 
\D_{\Z' \times \{0\}}$ as divisors lying 
over $(\X \times \Aone, \Z \times \0)$.
Alternatively, using \cite[Remark 4.1.3]{KR19}, 
an equivalence of divisors lying over
$(\X \times \Aone, \Z \times \0)$ is
an equivalence $\D_{\Z} \simeq 
\D_{\Z' \times \{0\}}$ of schemes
over $\bm{S}_{\Z \times \{0\}} := \bm{S}
\times_{\X \times \Aone} \Z \times \{0\}$\footnotemark.
	\footnotetext{The product $\bm{S}
	\times_{\X \times \Aone} \Z \times \{0\}$ is defined unambiguously because 
	$\Bl_{\bm{i}} \circ \b$ is canonically homotopic to $\Bl_{\bm{f}} \circ \a$
	as specified by the commutative diagram.}
More precisely we have the following
commutative diagram,
\[\begin{tikzcd}
	{\bm{D} :=\bm{D}_{\Z} \simeq \bm{D}_{\Z' \times \0}} & {\bm{S}_{\Z' \times \0}} \\
	{\bm{S}_{\Z}} & {\bm{S}_{\Z \times \0}}
	\arrow[from=1-1, to=2-1]
	\arrow[from=2-1, to=2-2]
	\arrow[from=1-2, to=2-2]
	\arrow[from=1-1, to=1-2]
\end{tikzcd}\]
which is classified by a map
$\bm{D} \to \bm{S}_{\Z} 
\times_{\bm{S}_{\Z \times \0}}
\bm{S}_{\Z' \times \0}$.
It is now a routine exercise to
check that $\bm{S}_{\Z} 
\times_{\bm{S}_{\Z \times \0}}
\bm{S}_{\Z' \times \0} \simeq
\bm{S}_{\Z'}$, yielding a map
$\bm{D} \to \bm{S}_{\Z'}$ which
classifies a virtual Cartier divisor
on $\bm{S}$ over $(\X' \times \0,
\Z' \times \0)$. Thus, we obtain a map
$\bm{S} \to \Bl_{\Z' \times \0}(\X' \times \0)$,
verifying the universal property for
products.
\eproof

The claim having been shown,
we see that the restriction of
$\Bl_{\bm{f} \times \id_{\Aone}}$ to 
$\bm{D}_{\Z'}\X' := \Bl_{\Z' \times \0}(\X' \times \Aone)
\setminus \Bl_{\Z' \times \0}(\X' \times \0)$
induces a map
	\[\D_{\bm{f}}: \D_{\Z'}\X' \to \D_{\Z}\X. \]
The lemma now follows from
\cref{lemma: blowup map} and 
\cite[Theorem 4.1.13]{KR19}, which describes the complements
of $\Bl_{\Z'}\X'$ and $\Bl_{\Z}\X$ inside each of
the terms in diagrams \labelcref{diagram:
blowup map 1}) and 
\labelcref{diagram: blowup map 2}.
\eproof
 
\subsubsection{}
Consistent with our notation
conventions, denote by
$D_{\Z/\X}$ the underlying classical
scheme of $\bm{D}_{\Z/\X}$, and by
$N_{\Z/\X}$ the underlying
classical scheme of $\bm{N}_{\Z/\X}$.
The following commutative diagram 
of schemes is obtained as the classical 
truncation of diagram \labelcref{diagram: deformation 2}:
\begin{equation}
\label{diagram: classical deformation 2}
\begin{tikzcd}
	{N_{\Z/\X}} & {D_{\Z/\X}} & {X \times \G_m} & X \\
	{\{0\}} & \Aone & {\G_m}
	\arrow[hook, from=2-1, to=2-2]
	\arrow[from=2-3, to=2-2]
	\arrow[from=1-1, to=2-1]
	\arrow[from=1-3, to=2-3]
	\arrow["{\pr_1}", from=1-3, to=1-4]
	\arrow["{j_{\neq 0}}"', from=1-3, to=1-2]
	\arrow["{s}", hook, from=1-1, to=1-2]
	\arrow["t"', from=1-2, to=2-2]
	\arrow["p", curve={height=-25pt}, from=1-2, to=1-4],
\end{tikzcd}
\end{equation}
and the following diagram of
schemes is obtained as the classical
truncation of diagram
\labelcref{diagram: map of deformations}:
\begin{equation}
\label{diagram: classical map of deformations}
\begin{tikzcd}
	{N_{\Z'/\X'}} & {D_{\Z'/\X'}} & {X' \times \G_m} & {X'} \\
	{N_{\Z/\X}} & {D_{\Z/\X}} & {X \times \G_m} & X \\
	{\{0\}} & \Aone & {\G_m}
	\arrow["{s'}", hook, from=1-1, to=1-2]
	\arrow["s", hook, from=2-1, to=2-2]
	\arrow["{N_{\bm{f}}}", from=1-1, to=2-1]
	\arrow["{D_{\f}}"', from=1-2, to=2-2]
	\arrow["{f \times \id_{\G_m}}", from=1-3, to=2-3]
	\arrow["f", from=1-4, to=2-4]
	\arrow["{\pr_2}", from=2-3, to=3-3]
	\arrow["t"', from=2-2, to=3-2]
	\arrow[from=2-1, to=3-1]
	\arrow[hook, from=3-1, to=3-2]
	\arrow[from=3-3, to=3-2]
	\arrow["{j_{\neq 0}}"', from=2-3, to=2-2]
	\arrow["{j_{\neq 0}'}"', from=1-3, to=1-2]
	\arrow["{\pr_1'}", from=1-3, to=1-4]
	\arrow["{\pr_1}", from=2-3, to=2-4]
	\arrow["{t'}"{pos=0.3}, curve={height=-20pt}, from=1-2, to=3-2].
\end{tikzcd}
\end{equation}

\brem
We draw attention
to the fact that all the squares 
in diagrams
\labelcref{diagram: classical map of
deformations}
and \labelcref{diagram: classical 
deformation 2} are Cartesian.
\erem


\subsubsection{$\G_m$\-action}
The classical construction of deformation
to the normal cone has a natural 
$\G_m$\-action induced by the grading 
on the Rees algebra. While there is no
obvious grading on the structure sheaf
of $\D_{\Z/\X}$ (indeed, we have not
even considered any explicit model for the
blow-up as locally-ringed space), there
is a natural $\G_m$\-action obtained by
application of \cref{lemma:
blowup map}.

Let $\bm{m}: \X \times \Aone \times \G_m
\to \X \times \Aone$ be the obvious $\G_m$\-action
map given by scaling the second factor.
The following square is Cartesian:
\[\begin{tikzcd}
	{\Z \times \0 \times \G_m} & {\X \times \Aone \times \G_m} \\
	{\Z \times \0} & {\X \times \Aone}
	\arrow["{\bm{m}}", from=1-2, to=2-2]
	\arrow[hook, from=2-1, to=2-2]
	\arrow[hook, from=1-1, to=1-2]
	\arrow["{\pr_{12}}", from=1-1, to=2-1]
	\arrow["\lrcorner"{anchor=center, pos=0.125}, draw=none, from=1-1, to=2-2]
\end{tikzcd}\]
We therefore obtain the map of blow-ups,
\[\Bl_{\bm{m}}:
	\Bl_{\Z \times \0 \times \G_m}(\X \times \Aone \times \G_m)
		\to \Bl_{\Z \times \0}(\X \times \Aone).\] 
Identifying $\Bl_{\Z \times \0 \times \G_m}(\X \times \Aone \times \G_m)$
and $\Bl_{\Z \times \0}(\X \times \Aone) \times \G_m$
via \cite[Theorem 4.1.5(ii)]{KR19},
$\Bl_{\bm{m}}$ gives a $\G_m$\-action
on $\Bl_{\Z \times \0}(\X \times \0)$.

For reasons similar to those
in the proof of
\cref{cor: deformation map}, the
following diagram is Cartesian:
\[\begin{tikzcd}
	{\Bl_{\Z \times \0}(\X \times \0) \times \G_m} & {\Bl_{\Z \times \0}(\X \times \Aone) \times \G_m} \\
	{\Bl_{\Z \times \0}(\X \times \0)} & {\Bl_{\Z \times \0}(\X \times \Aone)}
	\arrow["{\pr_{12}}", from=1-1, to=2-1]
	\arrow["{\Bl_{\bm{m}}}", from=1-2, to=2-2]
	\arrow[hook, from=2-1, to=2-2]
	\arrow[hook, from=1-1, to=1-2]
\end{tikzcd}\]
Therefore, $\Bl_{\bm{m}}$ induces the map
\[\bm{act}_{\G_m}: \D_{\Z/\X} \times \G_m \to \D_{\Z/\X},\]
endowing $\D_{\Z/\X}$ with a $\G_m$\-action.
We denote by $\act: D_{\Z/\X} \times \G_m \to
D_{\Z/\X}$ the corresponding action on 
underlying classical schemes.

\brem
\Cref{lemma: blowup map}
ensures that the action on $\X \times \G_m \subset
\D_{\Z/\X}$ is given by expected scaling on the second 
factor, while the action on $\N_{\Z/\X}$ is
given the canonical scaling of fibers,
so $\bm{act}_{\G_m}$ is indeed
a reasonable notion of $\G_m$\-action
on $\D_{\Z/\X}$.
\erem


\printbibliography

\end{document}